\documentclass[12pt, a4paper]{amsart}

\usepackage{graphicx}
\usepackage{mathtools}
\usepackage{amssymb}
\usepackage{amsmath, amsthm}
\usepackage{amscd}

\usepackage{enumitem}
\setlist{nosep}

\usepackage[all]{xy}

\usepackage{hyperref}
\hypersetup{hidelinks, hypertexnames=false}
\hypersetup{
    hypertexnames=false,
    colorlinks,
    linkcolor={red!30!black},
    citecolor={blue!50!black},
    urlcolor={blue!80!black}
}

\usepackage[margin=2cm]{geometry}

\usepackage[dvipsnames]{xcolor}
\usepackage{cancel, soul}
\usepackage{pdfcomment}

\allowdisplaybreaks

\newenvironment{axiom}[1]
{
\par\smallskip\noindent
{\bf Axiom #1}\begin{it}
}
{
\end{it}\par\smallskip
}

\newtheorem{theorem}[subsection]{Theorem}
\newtheorem{lemma}[subsection]{Lemma}
\newtheorem{sublemma}[subsubsection]{Lemma}
\newtheorem{proposition}[subsection]{Proposition}
\newtheorem{corollary}[subsection]{Corollary}

\newtheorem{claim}[subsection]{Claim}

\newtheorem{remark}[subsection]{Remark}
\newtheorem{example}[subsection]{Example}

\newtheorem{definition}[subsection]{Definition}

\makeatletter
\@addtoreset{equation}{section}
\@addtoreset{figure}{section}
\@addtoreset{table}{section}
\makeatother

\newtheorem{theorem*}{Theorem}

\makeatletter
\newcommand\testshape{family=\f@family; series=\f@series; shape=\f@shape.}
\def\myemphInternal#1{\if n\f@shape%
\begingroup\itshape #1\endgroup\/%
\else\begingroup\sffamily #1\endgroup%
\fi}
\def\myemph{\futurelet\testchar\MaybeOptArgmyemph}
\def\MaybeOptArgmyemph{\ifx[\testchar \let\next\OptArgmyemph
                 \else \let\next\NoOptArgmyemph \fi \next}
\def\OptArgmyemph[#1]#2{\index{#1}\myemphInternal{#2}}
\def\NoOptArgmyemph#1{\myemphInternal{#1}}
\makeatother

\newcommand\bC{{\mathbb C}}

\newcommand\bN{{\mathbb N}}

\newcommand\bR{{\mathbb R}}

\newcommand\bZ{{\mathbb Z}}

\newcommand\FF{{\mathcal F}}

\newcommand\isot{\mathrm{isot}}

\newcommand\cov[1]{\tilde{#1}}
\newcommand\eps{\varepsilon}

\newcommand\id{\mathrm{id}}          
\newcommand\Int{\mathrm{Int}}        
\newcommand\Fix[1]{\mathrm{Fix}(#1)} 
\newcommand\supp{\mathrm{supp\,}}    
\newcommand\wrm[1]{\wr_{#1}} 
\newcommand\ptnum[1]{|#1|}
\newcommand\MOD[2]{#1\,(\mathrm{mod}\,#2)}

\newcommand\AxCrPt{{\text{\rm(L)}}}
\newcommand\AxBd{{\text{\rm(B)}}}

\newcommand\Bman{B}

\newcommand\Mman{M}
\newcommand\Nman{N}
\newcommand\Pman{P}
\newcommand\Qman{Q}

\newcommand\Uman{U}
\newcommand\Vman{V}
\newcommand\Wman{W}
\newcommand\Xman{X}
\newcommand\Yman{Y}
\newcommand\Zman{Z}

\newcommand\tPman{\widetilde{\Pman}}

\newcommand\Circle{S^1}            
\newcommand\UInt{[0,1]}

\newcommand\Disk{D^2}              
\newcommand\Cylinder{A} 

\newcommand\Torus{T^2}             

\newcommand\PrjPlane{\bR{P}^2}    

\newcommand\Cyli[1]{\Cylinder_{#1}}

\newcommand\GL{\mathrm{GL}}

\newcommand\SO{\mathrm{SO}}

\newcommand\Orb{\mathcal{O}}        
\newcommand\Stab{\mathcal{S}}       
\newcommand\Diff{\mathcal{D}}       
\newcommand\Homeo{\mathcal{H}}      

\newcommand\DiffId{\Diff_{\id}}     
\newcommand\StabId{\Stab_{\id}}     

\newcommand\Cinfty{\mathcal{C}^{\infty}}
\newcommand\Cr[1]{\mathcal{C}^{#1}}
\newcommand\Ci[2]{\mathcal{C}^{\infty}(#1,#2)}               

\newcommand\func{f}
\newcommand\gfunc{g}
\newcommand\dif{h}
\newcommand\hdif{\cov{\dif}}
\newcommand\gdif{g}
\newcommand\hgdif{\cov{\gdif}}
\newcommand\qdif{q}

\newcommand\tfunc{\tilde{\func}}

\newcommand\DiffM{\Diff(\Mman)}

\newcommand\DiffIdM{\DiffId(\Mman)}

\newcommand\DiffMX{\Diff(\Mman, \Xman)}
\newcommand\DiffIdMX{\DiffId(\Mman, \Xman)}

\newcommand\Morse{\mathrm{Morse}}

\newcommand{\nb}{\mathrm{nb}}
\newcommand\DiffNbh{\Diff_{{\nb}}}     

\newcommand\Stabilizer[1]{\Stab(#1)}             
\newcommand\StabilizerPlus[1]{\Stab^{+}(#1)}     
\newcommand\StabilizerId[1]{\StabId(#1)}         
\newcommand\StabilizerIsotId[1]{\Stab'(#1)}      
\newcommand\StabilizerNbh[1]{\Stab_{{\nb}}(#1)}  
\newcommand\StabilizerNbhIsotId[1]{\Stab'_{{\nb}}(#1)}  

\newcommand\Orbit[1]{\Orb(#1)}                   
\newcommand\OrbitPathComp[2]{\Orb_{#2}(#1)}      

\newcommand\SingularSet[1]{\Sigma_{#1}}             

\newcommand\AutKRGraphStab[1]{\mathbf{G}}             

\newcommand\fStab{\Stabilizer{\func}}             
\newcommand\fStabId{\StabilizerId{\func}}         
\newcommand\fStabIsotId{\StabilizerIsotId{\func}}  

\newcommand\fOrbComp{\OrbitPathComp{\func}{\func}}  

\newcommand\fSing{\SingularSet{\func}}                

\newcommand\regU[1]{U_{#1}}
\newcommand\regN[1]{R_{#1}}
\newcommand\canN[1]{N_{#1}}

\newcommand\hMman{\cov{\Mman}}  

\newcommand\crLev{K}

\newcommand\dCyl{Q}
\newcommand\hdCyl{\cov{\dCyl}}

\newcommand\regNK{R_{K}}

\newcommand\fld{F}
\newcommand\flow{\mathbf{F}}
\newcommand{\hflow}{\cov{\flow}}

\newcommand\orb{\omega}

\newcommand\LStab{\mathcal{L}}

\newcommand\FolStab{\Delta} 
\newcommand\FolStabilizer[1]{\FolStab(#1)}
\newcommand\FolStabilizerNbh[1]{\FolStab_{{\nb}}(#1)}
\newcommand\FolStabilizerIsotId[1]{\FolStab'(#1)}
\newcommand\FolStabilizerNbhIsotId[1]{\FolStab'_{{\nb}}(#1)}

\newcommand\GKR{\mathbf{G}}
\newcommand\GrpKR[1]{\GKR(#1)}
\newcommand\GrpKRIsotId[1]{\GKR'(#1)}
\newcommand\GrpKRNbh[1]{\GKR_{{\nb}}(#1)}
\newcommand\GrpKRNbhIsotId[1]{\GKR'_{{\nb}}(#1)}

\newcommand\GH{\Gamma}
\newcommand\GHom[1]{\GH(#1)}
\newcommand\GHomIsotId[1]{\GH'(#1)}

\newcommand\halpha{\cov{\alpha}}

\newcommand\hXman{{\Xman}^{\partial}}

\newcommand\hZman{{\Zman}^{\partial}}
\newcommand\hUman{\Uman^{\partial}}

\newcommand\FSet{A}

\newcommand\HPlane{\mathbb{H}}

\newcommand\threepar[3]{{#1}_{#2,#3}}

\newcommand\Yi[2]{\threepar{\Yman}{#1}{#2}}

\newcommand\monoArrow{\lhook\joinrel\rightarrow}
\newcommand\longmonoArrow{\lhook\joinrel\longrightarrow}
\newcommand\xmonoArrow[1]{\lhook\joinrel\xrightarrow{~#1~}}

\newcommand\epiArrow{\rightarrow\!\!\!\!\!\to}
\newcommand\longepiArrow{\longrightarrow\!\!\!\!\!\to}
\newcommand\xepiArrow[1]{\xrightarrow{#1}\!\!\!\!\!\to}

\newcommand\DAFunc[1]{\Theta(#1)}
\newcommand\DFunc{\DAFunc{\fld}}

\newcommand\Sh[1]{\varphi_{#1}}

\newcommand\mprod{\prod} 
\newcommand\myprod{\mathop{\prod}}

\def\cnt{k}

\newcommand{\tMman}{\tilde{\Mman}}

\newcommand{\tXman}{\tilde{\Xman}}

\newcommand\jIncl{j}
\newcommand\jInclZ{\jIncl_0}

\newcommand\KRGraphf{\Gamma_{\func}}
\newcommand\EKRGraphf{\widehat{\Gamma}_{\func}}
\newcommand\PF[1]{\widehat{#1}}

\newcommand\crsdl{x}
\newcommand\cra{y}
\newcommand\crb{z}

\newcommand\funcSeq{\mathbf{b}}
\newcommand\seqStab[1]{\funcSeq(#1)}
\newcommand\seqStabIsotId[1]{\funcSeq'(#1)}
\newcommand\seqStabNbh[1]{\funcSeq_{\nb}(#1)}
\newcommand\seqStabNbhIsotId[1]{\funcSeq'_{\nb}(#1)}

\newcommand\DSG{{\funcSeq'}\!} 
\newcommand\sDSG[1]{\DSG\left(#1\right)}
\newcommand\dDSG[2]{\DSG(#1|_{#2},\partial#2)}
\newcommand\seqZ[1]{\mathbf{z}_{#1}}

\newcommand\seqWrm[2]{#1 \wr \seqZ{#2}}
\newcommand\seqTriv{\seqZ{0}} 
\newcommand\aSeq{\mathbf{q}}
\newcommand\uSeq{\mathbf{u}}
\newcommand\vSeq{\mathbf{v}}

\newcommand\kSeq{\mathbf{k}}
\newcommand\lSeq{\mathbf{l}}

\newcommand\classZ{\mathcal{Z}}
\newcommand\classStab{\mathcal{B}}
\newcommand\classGrp{\mathcal{P}}
\newcommand\classZBP{\classZ\classStab\classGrp}
\newcommand\ZBP{\classZ\classStab\classGrp_{min}}

\newcommand\JacIdeal{\added{\mathcal{I}_0(\func_z)}}
\newcommand\lfrm{\Phi}
\newcommand\zfrm[1]{\lfrm(#1)}

\newcommand\bZman{\mathbf{Z}}
\newcommand\bZmanX{\bZman^{fix}}
\newcommand\bZmanY{\bZman^{reg}}

\newcommand\XFix{\Xman}
\newcommand\XFixA{\XFix_0}
\newcommand\XFixB{\XFix_1}

\newcommand\hb{\eta}

\newcommand\ia{\alpha}
\newcommand\fc{\rho}

\newcommand\xA{A} 
\newcommand\xB{B} 
\newcommand\xC{C} 

\newcommand\kA{K} 
\newcommand\kB{L} 
\newcommand\kC{M} 

\newcommand\zA[1]{K_{#1}}
\newcommand\zB[1]{L_{#1}}
\newcommand\zC[1]{M_{#1}}

\newcommand\xv[1]{b_{#1}} 
\newcommand\xw[1]{c_{#1}} 

\newcommand\RB{\Xman}
\newcommand\cc[1]{\mathrm{(a#1)}}
\newcommand\ccInvm{\cc{1}}
\newcommand\ccSfR{\cc{2}}
\newcommand\ccSfX{\cc{3}}
\newcommand\ccTildeEta{\cc{4}}
\newcommand\ccIntervals{\cc{5}}
\newcommand\ccEdges{\cc{6}}

\newcommand\xc[1]{$\mathrm{(c#1)}$}
\newcommand\ppx{x}
\newcommand\ppy{y}
\newcommand\ppz{z}

\newcommand\rIso{\alpha}

\newcommand\UA{\Uman_{\Cylinder}}
\newcommand\UB{\Uman_{\Bman}}

\newcommand\lift{l}
\newcommand\lf[1]{\lift(#1)}

\newcommand\aDisk{D}

\newcommand\ggA{\mathsf{A}}
\newcommand\ggB{\mathsf{B}}
\newcommand\ggC{\mathsf{C}}

\title[Deformations of functions on surfaces]{Deformations of functions on surfaces by isotopic to the identity diffeomorphisms}
\author{Sergiy Maksymenko}
\email{maks@imath.kiev.ua}
\address{Institute of Mathematics of NAS of Ukraine, Tereshchenkivska st. 3, Kyiv, 01024 Ukraine}

\begin{document}

\providecommand\removed[1]{}
\providecommand\fremoved[1]{}
\providecommand\hidden[1]{}
\providecommand\added[1]{#1}
\providecommand\changed[2]{#2}
\providecommand\aremoved[1]{}
\providecommand\jadded[2]{#2}
\providecommand\jremoved[2]{}
\providecommand\jchanged[3]{#3}
\providecommand\xadded[4]{#4}
\providecommand\xremoved[4]{}
\providecommand\xchanged[5]{#5}

\keywords{surface, isotopy, Morse function, solvable group, wreath product}
\subjclass[2010]{
57S05, 
20E22, 
58B05 
}

\begin{abstract}
Let $M$ be a compact surface and $P$ be either $\mathbb{R}$ or $\Circle$.
For a smooth map $f:M\to P$ and a closed subset $V\subset M$\xadded{-0.8}{0}{1.1}{,} denote by $\mathcal{S}(f,V)$ the group of diffeomorphisms $h$ of $M$ preserving $f$, i.e.\! satisfying the relation $f\circ h = f$, and fixed on $V$.
Let also $\mathcal{S}'(f,V)$ be its subgroup consisting of diffeomorphisms isotopic relatively $V$ to the identity map $\mathrm{id}_{M}$ via isotopies that are not necessarily $f$-preserving.
The groups $\pi_0 \mathcal{S}(f,V)$ and $\pi_0 \mathcal{S}'(f,V)$ can be regarded as analogues of mapping class group for $f$-preserving diffeomorphisms.
The paper describes precise algebraic structure of groups $\pi_0 \mathcal{S}'(f,V)$ and some of their subgroups and quotients for a large class of smooth maps $f:M\to P$ containing all Morse maps, where $M$ is orientable and distinct from $2$-sphere and $2$-torus.
In particular, it is shown that for certain subsets $V$ ``adapted'' in some sense with $f$, the groups $\pi_0 \mathcal{S}'(f,V)$ are solvable and Bieberbach.
\end{abstract}

\maketitle

\section{Introduction}\label{sec:introduction}

Let $\Mman$ be a smooth $(\Cinfty)$ compact not necessarily connected surface, $\Pman$\jremoved{1.2}{be} either the real line $\bR$ or the circle $\Circle$, $\func:\Mman\to\Pman$\jremoved{1.3}{be} a smooth map, \xadded{-1.5}{2}{3.5}{and $\fSing$ the set of critical points of $\func$.}
A diffeomorphism $\dif:\Mman\to\Mman$ is said to be \myemph{$\func$-preserving} whenever $\func\circ\dif=\func$, which is equivalent to the assumption that $\dif$ leaves invariant each \myemph{level-set}, $\func^{-1}(c)$, $c\in\Pman$, of $\func$.

\begin{figure}[ht]
\centering
\footnotesize
\begin{tabular}{cccc}
\includegraphics[height=1.5cm]{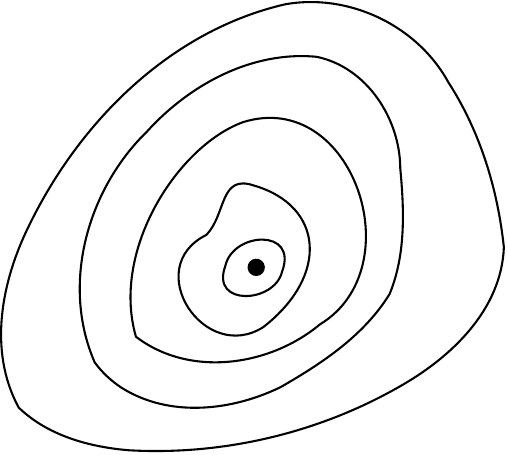}  &
\includegraphics[height=1.5cm]{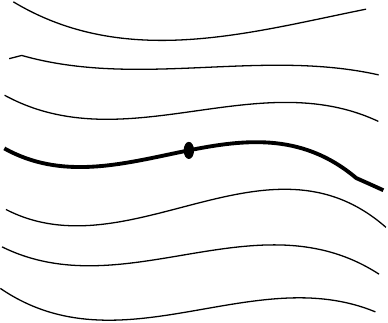}   &
\includegraphics[height=1.5cm]{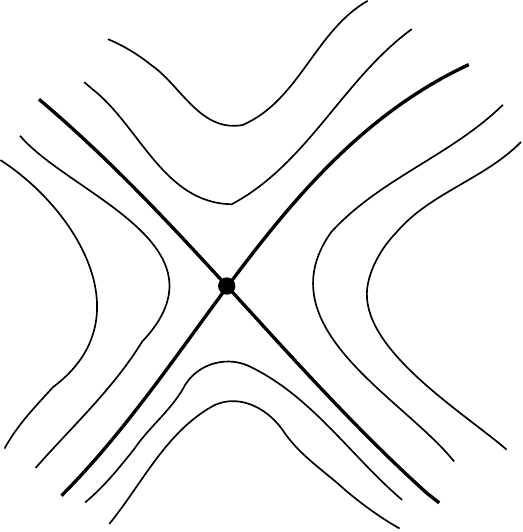}    &
\includegraphics[height=1.5cm]{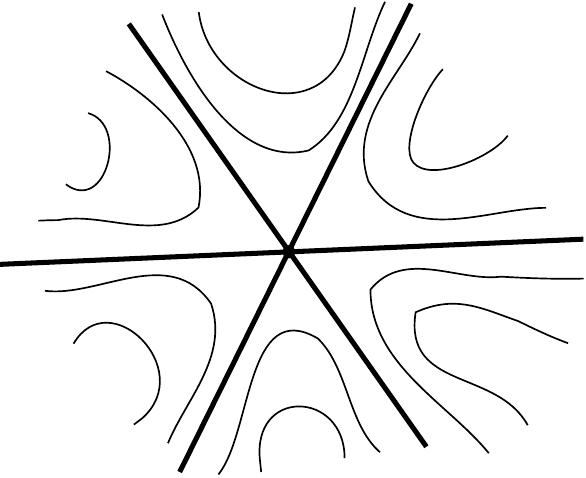}  \\
a) local extreme &
b) quasi-saddle, &
c) saddle, &
d) saddle, \\
& $m=1$ &
$m=2$ &
$m=3$ \\
$A^{+}_{2k-1}: \pm(x^{2k}+y^2)$ & $A_{2k}: \pm(x^{2k+1} + y^2)$   & $A^{-}_{2k-1}: \pm(x^{2k}-y^2)$ & $D_{2k}^{-}: x^2 y - y^{2k-1}$ \\
                                & $D^{+}_{2k}: x^2y + y^{2k-1}$   & $D_{2k+1}: \pm(x^2y+y^{2k})$     & \\
                                & $E_6 : \pm(x^3 + y^4)$          & $E_7: x^3+xy^3$                 & \\
                                & $E_8: x^3+y^5$                  &                           &
\end{tabular}
\caption{Topological structure of level-sets of homogeneous polynomials without multiple factors and simple singularities $A-D-E$.}\label{fig:isol_crit_pt}
\end{figure}

Denote by $\FF(\Mman,\Pman)$ the subspace of $C^{\infty}(\Mman,\Pman)$ consisting of maps $\func$ satisfying the following two axioms:
\begin{axiom}{\AxBd}
The map $\func$ takes a constant value at each connected component of $\partial\Mman$ and has no critical points in $\partial\Mman$.
\end{axiom}
\begin{axiom}{\AxCrPt}
For every critical point $z$ of $\func$\jadded{1.5}{,} the germ of $\func$ at $z$ is smoothly equivalent%
\footnote{Suppose $A_1, A_2, B_1,B_2$ are smooth manifolds, $a_1\in A_1$ and $a_2\in A_2$ are two points, and let $f_1:A_1\to B_1$ and $f_2:A_2\to B_2$ are smooth maps.
Then the germ of $f_1$ at $a_1$ is \myemph{smoothly} (resp. \myemph{topologically}) \myemph{equivalent} to the germ of $f_2$ at $a_2$ whenever there exist germs of diffeomorphisms (resp. homeomorphisms) $h:(A_1,a_1) \to (A_2, a_2)$ and $\phi:(B_1,f_1(a_1)) \to (B_2, f_2(b_2))$ such that $\phi\circ f_1 = f_2 \circ h$.
If in this case $\phi$ is the germ of the identity map, that is $f_1 = f_2 \circ h$, then the germ of $f_1$ at $a_1$ is \myemph{smoothly} (resp. \myemph{topologically}) \myemph{right equivalent} to the germ of $f_2$ at $a_2$.
If the germ of $f_1$ at $a_1$ is (smoothly or topologically) equivalent to a germ of a map $g:(\bR^n,0) \to (\bR^m,0)$, then $g$ will be called a \myemph{local representation} of $f_1$ at $a_1$.
} to a germ at $0\in\bR^2$ of some non-zero homogeneous polynomial $\func_z:\bR^2\to\bR$ without multiple factors.
\end{axiom}
\xadded{0.1}{0}{1.4}{Let $\func\in\FF(\Mman,\Pman)$.
Then axiom~\AxCrPt\ implies that every critical point of $\func$ is isolated.
Hence, since $\Mman$ is compact, the set of critical points of $\func$ is finite.}
Figure~\ref{fig:isol_crit_pt} below describes possible singularities satisfying Axiom~\AxCrPt.
In particular, since the polynomial $\pm x^2 \pm y^2$ (non-degenerate singularity) is homogeneous and has no multiple factors, we see that $\FF(\Mman,\Pman)$ contains an open and everywhere dense subset $\Morse(\Mman,\Pman)$ consisting of maps with non-degenerate singularities only (Morse maps).

Denote by $\DiffM$ the group of all $C^{\infty}$ diffeomorphisms of $\Mman$ endowed with $C^{\infty}$-topology, and let $\DiffIdM$ be its identity path component, i.e. it consists of $\dif\in\DiffM$ isotopic to $\id_{\Mman}$.
Let also $\Stabilizer{\func}$ be the group of all $\func$-preserving diffeomorphisms, $\StabilizerIsotId{\func} = \Stabilizer{\func} \cap \DiffIdM$, and $\StabilizerId{\func}$ be the identity path component of $\Stabilizer{\func}$ being also the identity path component of $\StabilizerIsotId{\func}$.

The present paper is devoted to study of the group of $\func$-preserving diffeomorphisms up to $\func$-preserving isotopies for all $\func\in\FF(\Mman,\Pman)$.
More precisely, we will describe the algebraic structure of homotopy groups $\pi_0\fStabIsotId \cong \fStabIsotId / \fStabId$ being analogues of mapping class group $\pi_0\DiffM\cong \DiffM/\DiffIdM$ for isotopic to identity $\func$-preserving diffeomorphisms, see \S\ref{sect:main_results}.

This kind of problems is motivated by significant progress in classifications of Hamiltonian dynamical systems of small degrees of freedom initiated by A.~Fo\-menko~\cite{Fomenko:RANDAN:ENG:1986}, \cite{Fomenko:RMS:1989}.

Classification of different kinds of Morse functions on surfaces was given in many papers, see e.g.\!\! A.~Bolsinov, S.~Matveev, A.~Fomenko~\cite{BolsinovMatveevFomenko:RMS:1990}, A.~Bolsinov, A.~Fomenko~\cite{BolsinovFomenko:ENG:2004}, A.~Oshemkov~\cite{Oshemkov:PSIM:ENG:1995}, E.~Kulinich~\cite{Kulinich:MFAT:1998}, V.~Manturov~\cite{Manturov:Atoms:1998}, V.~Sharko~\cite{Sharko:Chekyabinsk:1996, Sharko:UMZ:2003}, S.~Maksymenko~\cite{Maksymenko:UMJ:1999}, A.~Prishlyak~\cite{Prishlyak:UMZ:2000, Prishlyak:MFAT:2002}, V.~Arnold~\cite{Arnold:MPAG:2007, Arnold:TMIS:2010}.
Connected components of the space of Morse maps $\Morse(\Mman,\Pman)$ were independently computed for $\Pman=\bR$ by H.~Zieschang (unpublished), S.~Matveev published in the paper by E.~Kudryavtseva~\cite{Kudryavtseva:MatSb:ENG:1999}, V.~Shar\-ko~\cite{Sharko:Zb:1998}, and for $\Pman = \Circle$ by the author in~\cite{Maksymenko:Zb:1998}, \cite{Maksymenko:CMH:2005}.
Cobordism groups of Morse functions on surfaces were calculated by K.~Ikegami and O.~Saeki~\cite{IkegamiSaeki:MPCPS:2009}, and B.~Kalmar~\cite{Kalmar:KJM:2005}.

The group $\pi_0\fStabIsotId$ was studied by the author~\cite{Maksymenko:AGAG:2006, Maksymenko:ProcIM:ENG:2010, Maksymenko:MFAT:2010}, and by E.~Kudryavtseva~\cite{Kudryavtseva:MathNotes:2012, Kudryavtseva:MatSb:ENG:2013} in connection with another interpretation which we will now describe.

\xadded{0.4}{0}{3.2}{Given a set $\Xman$, denote by $\ptnum{\Xman}$ the number of its points if $\Xman$ finite and put $\ptnum{\Xman}:=\infty$ otherwise.}
For a closed subset $\Xman\subset\Mman$ let $\DiffMX$ be the subgroup of $\DiffM$ consisting of diffeomorphisms fixed on $\Xman$ and $\DiffIdMX$ be its identity path component.
Notice that $\DiffMX$ acts from the right on $C^{\infty}(\Mman,\Pman)$ as follows: if $\func\in C^{\infty}(\Mman,\Pman)$ and $\dif\in\DiffMX$, then the result of the action of $\dif$ on $\func$ is the composition map $\func\circ\dif:\Mman\to\Pman$.
Let
\begin{align}\label{equ:right_act_DiffM}
\Stabilizer{\func,\Xman} &= \{\dif\in\DiffMX \mid \func\circ\dif=\func\}, &
\Orbit{\func,\Xman} = \{\func\circ\dif \mid \dif\in\DiffMX \}
\end{align}
be respectively the stabilizer and orbit of $\func$.
Endow $\DiffM$ and $C^{\infty}(\Mman,\Pman)$ with the $C^{\infty}$ Whitney topologies.
Then they induce certain topologies on the stabilizers and orbits of maps $\func\in C^{\infty}(\Mman,\Pman)$.
Let also $\StabilizerId{\func,\Xman}$ be the identity path component of $\Stabilizer{\func,\Xman}$ and $\OrbitPathComp{\func,\Xman}{\func}$ be the path component of $\Orbit{\func,\Xman}$ containing $\func$.
If $\Xman=\varnothing$, then we will omit $\Xman$ from notation.
From this point of view the group $\Stabilizer{\func}$ of $\func$-preserving diffeomorphisms is the stabilizer of $\func$.
Denote also
\begin{align*}
\StabilizerIsotId{\func,\Xman} &= \Stabilizer{\func}\cap\DiffIdMX.
\end{align*}

\begin{theorem}\label{th:stab_orb:full_info}
{\rm(\cite{Maksymenko:AGAG:2006, Maksymenko:TrMath:2008, Maksymenko:ProcIM:ENG:2010, Maksymenko:UMZ:ENG:2012}).}
Let $\Mman$ be a connected compact surface, $\func\in\FF(\Mman,\Pman)$, and $\Xman$ a \jadded{3.1}{possibly empty} union of finitely many connected components of some level-sets of $\func$ and some critical points of $\func$.
Then the following statements hold.

\begin{enumerate}[leftmargin=*, label={\rm(\arabic*)}, wide]
\item\label{enum:th:stab_orb_full_info:Serre}
The map $p:\DiffMX\to\Orbit{\func,\Xman}$ defined by $p(\dif) = \func\circ\dif$ is a locally trivial principal $\Stabilizer{\func,\Xman}$-fibration.
In particular, the restriction $p:\DiffIdMX\to\OrbitPathComp{\func,\Xman}{\func}$ is a locally trivial principal $\StabilizerIsotId{\func,\Xman}$-fibration.
\added{Also $\OrbitPathComp{\func}{\func}$ is a Fr\`echet manifold, whence \cite{Palais:Top:1966} it has a homotopy type of a CW-complex.}

\item\label{enum:th:stab_orb_full_info:Sidf}
\removed{If $\Mman$ is orientable, ...}
\added{
The group $\StabilizerId{\func,\Xman}$ is \myemph{contractible} except for the cases described in the following table when it is \myemph{homotopy equivalent to the circle}:
\begin{center}\rm
\begin{tabular}{|c||p{8cm}|c|c|}\hline
Case & \centering $\func$ & $\fSing$ & $\Xman$  \\ \hline\hline
(A) & $\func:S^2\to\Pman$ has precisely two critical points and they are non-degenerate local extremes
& $\{ a, b \}$ & $\varnothing$, $\{a\}$, $\{b\}$, or $\{a, b\}$ \\ \hline
(B) & $\func:\Disk\to \Pman$ has a unique critical point $a$, and this point is a non-degenerate local extreme & $\{a\}$ & $\varnothing$ or $\{a\}$  \\ \hline
(C) & $\func:\Circle\times\UInt\to \Pman$ has no critical points & $\varnothing$ & $\varnothing$  \\ \hline
(D) & $\func: T^2\to \Circle$ has no critical points, whence it is a locally trivial fibration with fibre $\Circle$  & $\varnothing$ & $\varnothing$  \\ \hline
\end{tabular}
\end{center}
where $\UInt=\UInt$ and $T^2$ is a $2$-torus.
}

\item\label{enum:th:stab_orb_full_info:pikOf}
\changed{If}{Suppose} $\StabilizerId{\func,\Xman}$ is contractible.
Then $\pi_k\Diff(\Mman,\Xman) \cong \pi_k\Orbit{\func,\Xman}$ \removed{$\cong \pi_k\Mman$} for $k\geq3$, $\pi_2\Orbit{\func,\Xman}=0$, and we have the following short exact sequence%
\footnote{ Throughout the paper $\monoArrow$ means a ``\myemph{monomorphism}'' and $\epiArrow$ an ``\myemph{epimorphism}''.}:
\begin{equation}\label{equ:exact_seq_for_pi1OfX}
\pi_1\DiffIdMX \xmonoArrow{~p_1} \pi_1\OrbitPathComp{\func,\Xman}{\added{\func}} \xepiArrow{~~~} \pi_0\StabilizerIsotId{\func,\Xman}.
\end{equation}%
\xadded{-3}{0}{3.3}{}
If $\chi(\Mman) < |\Xman|$, then $\DiffIdMX$ is contractible as well, and \eqref{equ:exact_seq_for_pi1OfX} yields an isomorphism 
\begin{equation}\label{equ:pi1OfX_pi0SfX}
\pi_1\OrbitPathComp{\func,\Xman}{\func} \cong \pi_0\StabilizerIsotId{\func,\Xman}.
\end{equation}

\item\label{enum:th:stab_orb_full_info:DZk_Of_G}
\added{If $\Xman=\varnothing$, then there is another short exact sequence 
\begin{equation}\label{equ:exact_seq_for_pi1OfX:2}
\pi_1\DiffIdM \times \bZ^l \longmonoArrow \pi_1\OrbitPathComp{\func}{\func} \longepiArrow G
\end{equation}
for some $l\geq0$ and some finite group $G$.}

\item\label{enum:th:stab_orb_full_info:Off_V}
$\OrbitPathComp{\func,\Xman}{\func} = \OrbitPathComp{\func,\Xman \cup \Vman}{\func}$ for any collection $\Vman$ of boundary components of $\Mman$.

\item\label{enum:th:stab_orb_full_info:Stab_dM__Stab_X}
If $\Vman$ is a union of some boundary components of $\Mman$ and either $\chi(\Mman)<0$ or both $\partial\Mman$ and $\Vman$ are non-empty, then the inclusion $i: \StabilizerIsotId{\func, \partial\Mman}\subset \StabilizerIsotId{\func,\Vman}$ is a homotopy equivalence.

\item\label{enum:th:stab_orb_full_info:Of_generic}
Suppose $\func$ is \myemph{generic}, i.e. it is Morse and each critical level set of $\func$ contains exactly one critical point.
If $\Mman=S^2$ and $\func$ has exactly two critical points being local extremes, then $\OrbitPathComp{\func}{\func}$ is homotopy equivalent to $S^2$.
Otherwise, if $\Mman=S^2$ or $\PrjPlane$, then $\OrbitPathComp{\func}{\func}$ is homotopy equivalent to $\SO(3) \times (\Circle)^k$ for some $k\geq0$.
In all other cases $\OrbitPathComp{\func}{\func}$ is homotopy equivalent to $(\Circle)^k$ for some $k\geq0$.

\item\label{enum:th:stab_orb_full_info:Of_braid}
If $\func$ is Morse and has exactly $n$ critical points, then $\OrbitPathComp{\func}{\func}$ is homotopy equivalent to a certain covering space of the $n$-th configuration space of $\Mman$, which in turn is homotopy equivalent to some (possibly non-compact) $(2n-1)$-dimensional CW-complex.
In particular, $\pi_1\OrbitPathComp{\func}{\func}$ is a subgroup of the $n$-th braid group $\mathcal{B}_{n}(\Mman)$ of $\Mman$.
\end{enumerate}
\end{theorem}
\begin{proof}[Discussion of proofs]
\xadded{1}{0}{3.4}{}
\xadded{0}{0}{6.1}{}
\ref{enum:th:stab_orb_full_info:Serre}
Let $C_{0}(\bR^2)$ be the algebra of germs of $C^{\infty}$ functions $(\bR^2,0)\to(\bR,0)$ and $\JacIdeal = \bigl\{ \alpha\,\tfrac{\partial\func_z}{\partial x} + \beta\,\tfrac{\partial\func_z}{\partial y} \mid \alpha,\beta\in C_{0}(\bR^2)\bigr\}$ the ideal in $C_{0}(\bR^2)$ generated by partial derivatives of $\func_z$.
Then the codimension $\mathrm{codim} \func_z := \dim_{\bR}\bigl[C_{0}(\bR^2)/\JacIdeal\bigr]$ over $\bR$ of $\JacIdeal$ in $C_{0}(\bR^2)$ is called the \myemph{codimension} of the singularity of $\func_z$.
Furthermore, say that a smooth \jchanged{6.2}{function $\func:\Mman\to\bR$}{map $\func:\Mman\to P$}  is of \myemph{finite codimension}, if it has only finitely many critical points and all of them have finite codimensions.

In particular, by~\cite[Lemma~12]{Maksymenko:ProcIM:ENG:2010} all singularities satisfying Axiom~\AxCrPt\ have finite codimension, whence every function $\func\in\FF(\Mman,\Pman)$ has finite codimension. 

By F.~Sergeraert~\cite{Sergeraert:ASENS:1972}, for each $\func$ of finite codimension, its orbit $\mathcal{O}_{M\bR}(\func)$ with respect to the left-right action of $\Diff(\Mman)\times\Diff(\bR)$ is a Fr\'echet submanifold of $C^{\infty}(\Mman,\bR)$, and the natural map $p:\Diff(\Mman)\times\Diff(\bR)\to\mathcal{O}_{M\bR}(\func)$ defined by $p(\dif,\phi)=\phi\circ\func\circ\dif$ is a locally trivial fibration.
Therefore\jadded{7.1}{,} $\mathcal{O}_{M\bR}(\func)$ has the homotopy type of a certain CW-complex and its  topological type is completely determined by the corresponding homotopy type, \cite{Palais:Top:1966},

These results were extended in~\cite[\S11]{Maksymenko:AGAG:2006} to a larger class of actions of \myemph{tame Fr\'echet groups} on \myemph{tame Fr\'echet manifolds} which includes the above action~\eqref{equ:right_act_DiffM} of the group $\Diff(\Mman)$ and its orbits $\Orbit{\func}$ and thus proves~\ref{enum:th:stab_orb_full_info:Serre} for the cases when $\partial\Mman$ and $\Xman$ are empty. 
See also~\cite{Maksymenko:BSM:2006}, for precise relations between the homotopy types of $\mathcal{O}_{MP}(\func)$ and $\Orbit{\func}$.
Further, in~\cite[Theorem~2.2]{Maksymenko:UMZ:ENG:2012} statement~\ref{enum:th:stab_orb_full_info:Serre} was deduced for all $\func\in\FF(\Mman,\Pman)$ and $\Xman$ satisfying assumption of the theorem.

Statement~\ref{enum:th:stab_orb_full_info:Sidf} is a consequence of a series of papers~\cite{Maksymenko:CEJM:2009, Maksymenko:hamv2, Maksymenko:IUMJ:2010, Maksymenko:OsakaJM:2011} on diffeomorphisms preserving orbits of flows.
It is stated in~\cite[Theorem~1.3]{Maksymenko:AGAG:2006}, \cite[Theorem~3.7]{Maksymenko:OsakaJM:2011}, and~\cite[Theorem~3]{Maksymenko:ProcIM:ENG:2010} for $\Xman=\varnothing$ or $\fSing$, and in~\cite[Theorem~2.1]{Maksymenko:UMZ:ENG:2012} for all other cases of $\Xman$.
See \S\ref{sect:shifts_along_orbits} and Theorem~\ref{th:charact_Stabf} below describing the orientable case.

Statement~\ref{enum:th:stab_orb_full_info:pikOf} is a direct consequence of the exact sequence of homotopy groups of the fibration $p$, \added{contractibility of $\StabilizerId{\func,\Xman}$}, and a description of homotopy types of groups $\DiffId(\Mman,\Xman)$, see e.g.~\cite[\jadded{4.1}{Theorem~A}]{Smale:ProcAMS:1959}, \cite[\added{Theorems~1,2}]{EarleEells:BAMS:1967}, \cite[Theorems 1B,1C,1D]{EarleSchatz:DG:1970}, \cite[\added{Th\'eor\`eme~1}]{Gramain:ASENS:1973}.
\added{For the convenience of the reader we collect that information in the following table:
\begin{center}
\begin{tabular}{|c|c|c|c|c|c|}\hline
$\Mman$  & $\ptnum{\Xman}$ &  \multicolumn{4}{|c|}{$\DiffIdM$} \\ \cline{3-6} 
         &                 & homotopy  type  & $\pi_1$ & $\pi_2$ & $\pi_k$, $k\geq3$  \\ \hline \hline 
$S^2$, $\bR{P}^2$ & $0$ & $SO(3) \equiv \bR{P}^3$ & $\bZ_2$ & $0$ & $\pi_k S^3 = \pi_k S^2 = \pi_k \bR{P}^2$ \\ \hline
$T^2$ & $0$ & $T^2$ & $\bZ^2$ & $0$ & $0$ \\ \hline
$\Disk$, $\Circle\times\UInt$, $Mo$, $K$ & $0$ & & & & \\  \cline{1-2}
$\Disk$, $S^2$, $\bR{P}^2$           & $1$ & $\Circle$ & $\bZ$ & $0$ & $0$ \\ \cline{1-2}
$S^2$                              & $2$ & & & &\\ \hline
\multicolumn{2}{|c|}{other cases} & point  & $0$ & $0$ & $0$ \\ \hline 
\end{tabular}
\end{center}
where $Mo$ is a M\"obius band and $K$ is a Klein bottle.}
\jadded{4.2}{}In particular, we see that $\pi_2\DiffIdM=0$ and $\pi_k\Mman = \pi_k\DiffIdM$ for $k\geq3$.
Moreover, isomorphisms $\pi_k\OrbitPathComp{\func,\Xman}{\func} = \pi_k\DiffIdMX$ for $k\geq3$ follow from aspherity of $\StabilizerId{\func,\Xman}$.

Statement~\ref{enum:th:stab_orb_full_info:DZk_Of_G} is proved in~\cite[Eq.~(1.6)]{Maksymenko:AGAG:2006} for Morse functions, and in~\cite[Eq.~(5.2)]{Maksymenko:ProcIM:ENG:2010} for $\func\in\FF(\Mman,\Pman)$.
In fact, for the case when $\chi(\Mman)<0$, the sequence~\eqref{equ:exact_seq_for_pi1OfX:2} is isomorphic with the ``Bieberbach'' sequence~\eqref{equ:exact_seq_DSG}, see below, being the main object of study in this paper.

Statement~\ref{enum:th:stab_orb_full_info:Off_V} is a simple consequence of~\ref{enum:th:stab_orb_full_info:Serre} and is proved in~\cite[Corollary~2.1]{Maksymenko:UMZ:ENG:2012}.
Statement~\ref{enum:th:stab_orb_full_info:Of_generic} is established in~\cite[Theorems~1.5(3) \& 1.9]{Maksymenko:AGAG:2006}, and statement~\ref{enum:th:stab_orb_full_info:Of_braid} is showed in~\cite{Maksymenko:TrMath:2008}.

\ref{enum:th:stab_orb_full_info:Stab_dM__Stab_X}
\jadded{4.4}
By assumption, the identity path components of $\StabilizerIsotId{\func, \partial\Mman}$ and $\StabilizerIsotId{\func,\Vman}$ are contractible.
Therefore it suffices to show that $i$ induces an isomorphism between the corresponding $\pi_0$-groups.
As $\OrbitPathComp{\func,\Vman}{\func} = \OrbitPathComp{\func,\partial\Mman}{\func}$ due to~\ref{enum:th:stab_orb_full_info:Off_V}, we get from~\ref{enum:th:stab_orb_full_info:pikOf} the following isomorphisms:
$\pi_0 \StabilizerIsotId{\func, \partial\Mman} \cong \pi_1 \OrbitPathComp{\func,\partial\Mman}{\func}=
\pi_1 \OrbitPathComp{\func,\Vman}{\func} \cong
\pi_0 \StabilizerIsotId{\func, \Vman}$.
\end{proof}

E.~Kudryavtseva~\cite{Kudryavtseva:SpecMF:VMU:2012, Kudryavtseva:MathNotes:2012, Kudryavtseva:MatSb:ENG:2013}, computed homotopy types of spaces $\Morse(\Mman,\bR)$ for orientable surfaces $\Mman$, and, rediscovering many arguments from~\cite{Maksymenko:AGAG:2006}, obtained, in particular, the following result generalizing  Theorem~\ref{th:stab_orb:full_info}\ref{enum:th:stab_orb_full_info:Of_generic}:
\begin{theorem}\label{th:Kudryavtseva}{\rm\cite{Kudryavtseva:SpecMF:VMU:2012, Kudryavtseva:MathNotes:2012, Kudryavtseva:MatSb:ENG:2013}.}
Let $\Mman$ be a connected orientable surface, $\func\in\Morse(\Mman,\bR)$, \xadded{0.5}{0}{4.5}{ $\Fix{\StabilizerIsotId{\func}}$ the set of all critical points of $\func$ being fixed under each diffeomorphism $\dif\in\StabilizerIsotId{\func}$, and $|\Fix{\StabilizerIsotId{\func}}|$ the number of such points.}
Suppose $\chi(\Mman) < |\Fix{\StabilizerIsotId{\func}}|$ (which holds e.g.\! if $\chi(\Mman)<0$ or if $\func$ is generic and has at least one saddle critical point).
\added{Then there exists a free action of the group $G$ from~\eqref{equ:exact_seq_for_pi1OfX:2} on an $k$-dimensional torus $T^k := (\Circle)^k$ such that $\OrbitPathComp{\func}{\func}$ has the homotopy type of the quotient $T^k/G$ if $\Mman\not=S^2$, and the homotopy type of $(T^k/G) \times \SO(3)$ if $\Mman=S^2$.}
\end{theorem}

Moreover, in~\cite{Kudryavtseva:ENG:DAN2016} she also described the homotopy types of the spaces of smooth functions $f$ with prescribed local singularities of $A_{\mu}$-types, $\mu\in\bN$, without any restrictions, as well as the homotopy types of orbits $\mathcal{O}(f)$ of such $f$.

The group ${G}$ from Theorem~\ref{th:Kudryavtseva} can be regarded as the group of automorphisms of the Kronrod-Reeb graph of $\func$ generated by diffeomorphisms from $\StabilizerIsotId{\func}$\added{, see Section~\ref{sect:FolStab_Grp}.}
It is trivial in generic case, so that result extends Theorem~\ref{th:stab_orb:full_info}\ref{enum:th:stab_orb_full_info:Of_generic} to non-generic case%
\footnote{
The papers~\cite{Kudryavtseva:MathNotes:2012} and \cite{Kudryavtseva:MatSb:ENG:2013} use terminology and notation different from~\cite{Maksymenko:AGAG:2006} and this is probably the reason of such an overlapping of results. 

\quad
To be precise, E.~Kudryavtseva, \cite{Kudryavtseva:MathNotes:2012}, \cite{Kudryavtseva:MatSb:ENG:2013}, studied the orbit, denoted by $[\func]$, of a function $\func \in \Morse(\Mman,\bR)$ with respect to the \textit{left-right} action of the product group $\Diff^{+}(\bR)\times \DiffM$ on $C^{\infty}(\Mman,\bR)$ defined by $(\phi,\dif) \func = \phi\circ\func\circ\dif$. 
The orbit of $\func$ with respect to the action of the identity path component $\Diff^{+}(\bR)\times \DiffIdM$ of $\Diff^{+}(\bR)\times \DiffM$ is denoted by $[\func]_{\isot}$ and is called the \myemph{isotopy class} of $\func$.
Then~\cite[Theorem~2.6 C)]{Kudryavtseva:MatSb:ENG:2013}, in particular, claims that in the case $\chi(\Mman) < |\Fix{\StabilizerIsotId{\func}}|$ the orbit $[\func]_{\isot}$ is homotopy equivalent to
$\bigl( (\Circle)^k/ \Gamma_{[f]}\bigr) \times R$, where $\Gamma_{[f]}$ is a certain finite group freely acting on torus $(\Circle)^k$,
and $R$ is defined as follows: $R:=\SO(3)$ if $\Mman = S^2$, $R := T^2$ if $\Mman = T^2$, $R := \Circle$ if $M = \Circle \times\UInt$ or $M = \Disk$, and $R$ is a point if $\chi(\Mman)<0$.
In fact the isotopy class $[\func]_{\isot}$ is a path component of $[\func]$, since by F.~Sergeraert \cite{Sergeraert:ASENS:1972}, the mapping $p:\Diff^{+}(\bR)\times \DiffM \to [\func]$ defined by $p(\phi,\dif) = \phi\circ\func\circ\dif$ is a Serre fibration.

\quad
Evidently, $\fOrbComp \subset [\func]_{\isot}$.
Moreover, by~\cite[Theorem~1.5]{Maksymenko:BSM:2006} for all $\func\in\FF(\Mman,\bR)$ this inclusion of the \myemph{right} orbit $\fOrbComp$ into the corresponding \myemph{left-right} orbit $[\func]_{\isot}$ is a homotopy equivalence.
Thus $\fOrbComp$ is homotopy equivalent to $\bigl((\Circle)^k/ \Gamma_{[f]}\bigr)\times R$ as well.

\quad
However, as one can show (see below), the group $\Theta^{*}$ defined by~\cite[\S3.2, Eq~(22)]{Kudryavtseva:MatSb:ENG:2013} and playing the main role in the appearance of torus $(\Circle)^k$ is exactly the same as the group $\pi_0\mathcal{D}^{+}(\Delta_f)$ of~\cite{Maksymenko:AGAG:2006} and denoted by $\pi_0\FolStabilizer{\func}$ in the present paper.
The generators of this group are isotopy classes of Dehn twists in $\Stabilizer{\func}$ along special regular components of level-sets of $\func$ corresponding to ``internal'' edges of Kronrod-Reeb graph of $\func$, see \added{Section~\ref{sect:FolStab_Grp}}, \cite[\S6, Theorem 6.2]{Maksymenko:AGAG:2006} and \cite[\S3.2]{Kudryavtseva:MatSb:ENG:2013}.

\quad
Let $\Stabilizer{\func}^H$ be the subgroup of $\Stabilizer{\func}$ consisting of diffeomorphisms isotopic to $\id_{\Mman}$ in the class of \textit{$\func$-preserving homeomorphisms}.
In other words, $\Stabilizer{\func}^H$ is the identity path component of $\Stabilizer{\func}$ with respect to $C^{0}$ (i.e. compact open) topology, while $\StabilizerId{\func}$ is the identity path component of $\Stabilizer{\func}$ with respect to $C^{\infty}$ topology.
Then $\StabilizerId{\func} \subset \StabilizerId{\func}^H$.
Moreover, by~\cite[Remark~1.4]{Maksymenko:AGAG:2006} $\StabilizerId{\func} = \StabilizerId{\func}^H$.

\quad
It then follows that the group $\widetilde{\Gamma}_{[f]_{\isot}}$ defined in~\cite[\S3.2, Eq~(37)]{Kudryavtseva:MatSb:ENG:2013} as
$\StabilizerIsotId{\func}/\StabilizerId{\func}^H$ is the same as the group $\pi_0\StabilizerIsotId{\func} = \StabilizerIsotId{\func}/\StabilizerId{\func}$ being also the kernel of the homomorphism $i_0:\pi_0\fStab\to\pi_0\DiffM$ induced by inclusion, see~\cite[\S8-9]{Maksymenko:AGAG:2006}, \cite{Maksymenko:MFAT:2010}, \cite{Maksymenko:UMZ:ENG:2012}.

\quad
Finally, the group $\Gamma_{[f]}$ acting on the torus (denoted by ${G}$ in Theorem~\ref{th:Kudryavtseva}) is defined by~\cite[\S3.2, Eq~(32)]{Kudryavtseva:MatSb:ENG:2013} and coincides with the group $G$ defined in~\cite[\S9]{Maksymenko:AGAG:2006}.
This group will be denoted by $\GrpKRIsotId{\func}$ in the present paper.
}.
\added{
On the other hand, the statement of Theorem~\ref{th:Kudryavtseva} for the case $\Mman\not=S^2$ also holds for non-orientable surfaces distinct from $\bR{P}^2$ and can be easily deduced from Theorem~\ref{th:stab_orb:full_info} as follows.
\begin{corollary}\label{cor:Of_hom_type}
Suppose $\Mman\not=S^2, \bR{P}^2$.
Then for every $\func\in\Morse(\Mman,\Pman)$ there exists a free action of the group $G$ from~\eqref{equ:exact_seq_for_pi1OfX:2} on an $k$-dimensional torus $T^k$ such that $\OrbitPathComp{\func}{\func}$ is homotopy equivalent to $T^k/G$.
\end{corollary}
\begin{proof}
If $\StabilizerId{\func}$ is homotopy equivalent to $\Circle$ and $\Mman\not=S^2$, then by~\cite[Theorem~1.9]{Maksymenko:AGAG:2006} $\OrbitPathComp{\func}{\func}$ is either contractible or homotopy equivalent to $1$-torus $\Circle$.
Hence we may restrict ourselves to the case when $\StabilizerId{\func}$ is contractible and thus by Theorem~\ref{th:stab_orb:full_info}\ref{enum:th:stab_orb_full_info:pikOf}, $\OrbitPathComp{\func}{\func}$ is aspherical.

Moreover, under assumptions on $\Mman$, the group $\pi_1\DiffIdM$ is free abelian, and therefore $\pi_1\DiffIdM\times\bZ^l$ in~\eqref{equ:exact_seq_for_pi1OfX:2} is also free abelian of some rank $k$.
Then~\eqref{equ:exact_seq_for_pi1OfX:2} means that $\pi_1\OrbitPathComp{\func}{\func}$ is a crystallographic group, see Definition~\ref{def:solvable_bieberbach}.

Also due to Theorem~\ref{th:stab_orb:full_info}\ref{enum:th:stab_orb_full_info:Of_braid}, $\pi_1\OrbitPathComp{\func}{\func}$ is a subgroup of the braid group of $\Mman$, and therefore it is torsion free, \cite{FoxNeuwirth:MS:1962, Dyer:MZ:1980}.
Hence $\pi_1\OrbitPathComp{\func}{\func}$ is a Bieberbach group and~\eqref{equ:exact_seq_for_pi1OfX:2} is its ``Bieberbach'' sequence.
So by Bieberbach theorem (see below Theorem~\ref{th:bieberbach}) there exists a free action of $G$ on $T^k$ such that $\OrbitPathComp{\func}{\func}$ is homotopy equivalent to $T^k/G$.
\end{proof}
}


We will describe the algebraic structure of the groups $\pi_0\StabilizerIsotId{\func,\Xman}$, ${G}$, and the quotient map $\pi_0\StabilizerIsotId{\func,\Xman} \to {G}$ for all $\func\in\FF(\Mman,\Pman)$ and all $\func$-adapted submanifolds $\Xman$ in the case when $\Mman$ is orientable and distinct from $S^2$ and $T^2$, see Theorems~\ref{th:stab:chi_neg}, \ref{th:stab:annulus}, \ref{th:stab:disk:one_crpt}, \ref{th:stab:disk_ann:gen_case}, and~\ref{th:solvable_groups}.
\added{In particular, we will show that those groups are \myemph{solvable} and $\pi_0\StabilizerIsotId{\func,\Xman}$ is also \myemph{Bieberbach} without referring to the fact that braid groups for such surfaces are torsion free.}

\medskip
{\bf Class $\FF(\Mman,\Pman)$.}
Notice that every real homogeneous polynomial $\func_z:\bR^2\to\bR$ splits into a product $\func_z = L_1\cdots L_p \cdot Q_1 \cdots Q_q$, where each $L_i(x,y) = a_i x + b_i y$ is a linear function and each $Q_j(x,y) = c_j x^2 + d_j xy + e_j y^2$ is an irreducible over $\bR$ quadratic form.
The assumption that $\func_z$ has no multiple factors means that $L_i/L_{i'}\not=\mathrm{const}$ for $i\not=i'$, and $Q_j/Q_{j'}\not=\mathrm{const}$ for $j\not=j'$.
In particular, it also follows that the origin $0\in\bR^2$ is an \myemph{isolated} critical point of $\func_z$.

It is known that every germ of a smooth isolated singularity $\func:(\bC,0)\to(\bR,0)$ is \myemph{topologically} equivalent to a germ $g:(\bC,0) \to (\bR,0)$ of the form:
\begin{equation}\label{equ:isol_sing}
g(z) =
\begin{cases}
|z|^2, & \text{if $0\in\bC$ is a \myemph{local extreme} of $\func$, \cite{Dancer:2:JRAM:1987}}, \\
Re(z^m), & \text{otherwise, \cite{Prishlyak:TA:2002}}.
\end{cases}
\end{equation}

In the second case the critical point $0\in\bC$ is called a \myemph{saddle}.
For $m=1$ we will also call this point a \myemph{quasi-saddle}.
Topological structure of level-set of $\func\in\FF(\Mman,\Pman)$ near its critical point and comparison with level sets of simple singularities $A-D-E$ is presented in Figure~\ref{fig:isol_crit_pt}.

Let $\StabilizerPlus{\func_z}$ be the group of germs of orientation preserving diffeomorphisms $\dif$ of $(\bR^2,0)$ such that $\func_z\circ \dif = \func_z$.
In particular $\dif(0)=0$, whence we get a natural \myemph{Jacobi} homomorphism
$J:\StabilizerPlus{\func_z}\longrightarrow \GL^{+}(\bR,2)$
associating to each $\dif\in\StabilizerPlus{\func_z}$ its Jacobi matrix $J(\dif)$ at $0$.

Let also $\LStab(\func_z)$ be the subgroup of $\GL^{+}(2,\bR)$ preserving $\func_z$, that is the group of orientation preserving \textit{linear} isomorphisms $A:\bR^2\to\bR^2$ such that $\func_z(Au) = \func_z(u)$ for all $u\in\bR^2$.
Since every such $A$ can be regarded as a diffeomorphism of $\bR^2$ with $J(A)=A$, we see that $\LStab(\func_z) \subset J(\StabilizerPlus{\func_z})$.
\begin{lemma}\label{lm:LStabf}{\rm\cite[Lemma~6.2]{Maksymenko:MFAT:2009}.}
$\LStab(\func_z) = J(\StabilizerPlus{\func_z})$.
Moreover, after a proper linear change of coordinates in $\bR^2$ and replacing if necessary $\func_z$ with $-\func_z$, one can assume that $\func_z$ has the following properties.
\begin{enumerate}[label=$\mathrm{(\roman*)}$, leftmargin=*]
\item\label{enum:sing:nondeg_loc_extr}
If $z$ is a \myemph{non-degenerate local extreme}, then $\func_z(x,y) =x^2+y^2$ and $\LStab(\func_z) = \SO(2)$.
\item\label{enum:sing:nondeg_saddle}
If $z$ is a \myemph{non-degenerate saddle}, then $\func_z(x,y)=xy$ and
$\LStab(\func_z) = \{ \left(\begin{smallmatrix} t^{-1} & 0 \\ 0 & t\end{smallmatrix} \right) \mid t\not=0\}$.
\item\label{enum:sing:deg_saddle}
In all other cases there exists $m\geq1$ such that $\LStab(\func_z)$ is a finite cyclic subgroup of $\SO(2)$ of order $m$ generated by the rotation by $2\pi/m$.
\added{If $\deg\func_z$ is even, e.g. if $z$ is a \myemph{degenerate local extreme}, then $-\id_{\bR^2} \in \LStab(\func_z)$ is an element of order $2$, whence $m$ is even as well.}
That number $m$ will be called the \myemph{symmetry index} of $z$.
\end{enumerate}
\end{lemma}

\begin{remark}\label{rem:deformation_in_fStab}\rm
Suppose $z$ is a local extreme of $\func$ and $\gamma:\UInt \to \StabilizerPlus{\func_z}$ is a continuous map with respect to $C^{k}$-topology, $k\geq1$, on $\StabilizerPlus{\func_z}$.
This assumption implies, in particular, that the Jacobi matrix of $J(\gamma(t))$ of $\gamma(t)$ at $0$ is continuous in $t$.
Hence if $z$ is a non-degenerate local extreme of $\func$, then by Lemma~\ref{lm:LStabf}\ref{enum:sing:nondeg_loc_extr}, $J(\gamma(t))$ can be arbitrary matrix from $\SO(2)$.
However, if $z$ is a degenerate local extreme, then by Lemma~\ref{lm:LStabf}\ref{enum:sing:deg_saddle}, $J(\gamma(t))$ can take only finitely many values, whence $J(\gamma(t))$ must be the same for all $t$.
A consequence of this difference between non-degenerate and degenerate local extremes is reflected in Theorem~\ref{th:stab:disk:one_crpt}.
\end{remark}

\section{Algebraic preliminaries}\label{sect:algebraic}

Let $\uSeq: \kA_1 \xrightarrow{\alpha_1} \cdots \xrightarrow{\alpha_{k-1}} \kA_{k}$ and $\vSeq: \kB_1 \xrightarrow{\beta_1}  \cdots \xrightarrow{\beta_{k-1}} \kB_{k}$ be two sequences of homomorphisms of groups.
By a \myemph{homomorphism} $\gamma=(\gamma_1,\ldots,\gamma_{k}): \uSeq \to \vSeq$ we will mean a collection of homomorphisms $\gamma_i:\kA_i\to \kB_i$, $i=1,\ldots,k$, making commutative the following diagram:
\[
\aligned
\xymatrix@R=1.1em{
\kA_1 \ar[r]^-{\alpha_1}  \ar[d]^-{\gamma_1} & 
\kA_2 \ar[r]^-{\alpha_2}  \ar[d]^-{\gamma_2} & 
\cdots \ar[r]^-{\alpha_{k-1}} & \kA_{k} \ar[d]^-{\gamma_k} \\
\kB_1 \ar[r]^-{\beta_1}   & 
\kB_2 \ar[r]^-{\beta_2}   & 
\cdots \ar[r]^-{\beta_{k-1}} & \kB_{k}
}
\endaligned
\]
We will say that $\gamma$ is an \myemph{epimorphism} (resp. \myemph{monomorphism}, \myemph{isomorphism}) if each $\gamma_i$ is so.
Given $n$ sequences of homomorphisms $\aSeq^i: \kA^{i}_1 \xrightarrow{\alpha^{i}_1} \kA^{i}_2 \xrightarrow{\alpha^{i}_2} \cdots \xrightarrow{\alpha^{i}_{k-1}} \kA^{i}_{k}$, $i=1,\ldots,n$, their \myemph{product} is the sequence
\[
\prod_{i=1}^{n}\aSeq^i: \  \prod_{i=1}^{n} \kA^{i}_1 \xrightarrow{\alpha^{1}_1\times \cdots \times\alpha^{n}_1}
\prod_{i=1}^{n} \kA^{i}_2 \xrightarrow{\alpha^{1}_2\times \cdots \times\alpha^{n}_2}
\cdots \xrightarrow{\alpha^{1}_{k-1}\times \cdots \times\alpha^{n}_{k-1}}
\prod_{i=1}^{n}\kA^{i}_{k}.
\]

More generally, one can define in an obvious way similar notions for commutative diagrams of arbitrary fixed type not only for chains of homomorphisms.

{\bf Wreath products.}
Recall that a group $G$ \myemph{splits into a direct product of its subgroups $G_1,\ldots,G_k$} if $G_1,\ldots,G_k$ generate $G$, $G_i\cap G_j = \{e\}$ for $i\not=j$, and $g_i g_j = g_j g_i$ for all $g_i\in G_i$, $g_j \in G_j$, $i\not=j$.
In this case the map $\psi:\prod_{i=1}^{k} G_i \to G$, $\psi(g_1,\ldots,g_k) = g_1\cdots g_k$, is an isomorphism.

Also, a group $G$ \myemph{splits into a semidirect product $A\rtimes Z$ of its subgroups $A$ and $Z$} if $A$ and $Z$ generate $G$, $A$ is normal in $G$, and $A \cap Z = \{e\}$.
In this case the we have the following commutative diagram with exact rows:
\begin{equation}\label{equ:semidirect_prod}
\aligned
\xymatrix@R=1.2em{
\ A            \           \ar@{=}[d] \ar@{^{(}->}[rr]^-{a \mapsto (a,e)} &&
\ A  \rtimes Z \ \ar[d]^-{\psi}_-{\cong} \ar@{->>}[rr]^-{(a,z) \mapsto z \ \ } &&
\ Z            \ \ar[d]_-{\cong} \\
\ A            \ \ar@{^{(}->}[rr] && \ G  \  \ar@{->>}[rr]^-{q} && \ G/A \
}
\endaligned
\end{equation}
where $A \rtimes Z$ is a semidirect product corresponding to the action of $Z$ on $A$ by conjugations.
The latter group $A \rtimes Z$ is by definition the Cartesian product $A\times Z$ of sets with the operation defined by $(a,y)(b,z) = (ayby^{-1}, yz)$.
Then in~\eqref{equ:semidirect_prod}, $\psi:A\rtimes Z \to G$ is given by $\psi(a,z)=az$, and $q:G\to G/A$ is a quotient homomorphism yielding an isomorphism of $Z$ onto $G/A$.

In particular, 
\jadded{8.1}{}
if $\eta: G \to \added{Z}$ is an \myemph{epimorphism} admitting a section, i.e.\! a homomorphism $s:Z \to G$ such that $\eta\circ s = \id_{Z}$, then $G$ splits into a semidirect product $\ker(\eta) \rtimes s(Z)$.

Notice that each epimorphism $\eta:G \to \bZ$ onto the\jchanged{8.2}{ groups}{ group} $\bZ$ of integers, always has a section: just take any $g\in G$ with $\eta(g)=1$, and put $s(n) = g^n$.
Then $G \cong \ker(\eta) \rtimes \bZ$ is a Cartesian product of sets $\ker(\eta) \times \bZ$ with the following multiplication:
\begin{equation}\label{equ:prod_in_GrtimesZ}
(a,m)(b,n) = \bigl(a g^{m} b g^{-m}, \ m+n\bigr).
\end{equation}

\begin{definition}\label{defn:wreath_products}
For a group ${G}$ and $m \geq 1$ let
\[ {G}\wrm{m}\bZ \ := \ \underbrace{{G}\times\cdots\times {G}}_{m} \ \rtimes \ \bZ\]
be the semidirect product ${G}^m \rtimes \bZ$ corresponding to the \myemph{non-effective} action of $\bZ$ on ${G}^m$ by cyclic shifts of coordinates to the left.
Thus ${G}\wrm{m}\bZ$ is a cartesian product of sets $G^m \times \bZ$ with the following operation:
\begin{equation}\label{equ:mult_in_GmZ}
(a_0,\ldots,a_{m-1}, k)   (b_0,\ldots,b_{m-1}, l)  =  (a_0 b_k, a_{1} b_{k+1}, \cdots, a_{m-1} b_{k-1}, k+l),
\end{equation}
where all indices are taken modulo $m$.
Let also
\[ {G}\wr\bZ_m \ := \ \underbrace{{G}\times\cdots\times {G}}_{m} \ \rtimes \ \bZ_m\]
be the standard \myemph{wreath product} of ${G}$ and $\bZ_m$, i.e.\! the semidirect product ${G}^m \rtimes \bZ_m$ corresponding to the \myemph{effective} action of $\bZ_m$ on ${G}^m$ by cyclic shifts of coordinates to the left.
Evidently, ${G}\wrm{1}\bZ \cong {G} \times \bZ$ and ${G}\wr\bZ_1 = {G}$.
\end{definition}

\begin{lemma}\label{lm:no_fin_el}
If $G$ is \myemph{torsion free}, then $G\wrm{m}\bZ$, $m\geq1$, is torsion free as well.
\end{lemma}
\begin{proof}
Let $g = (g_1,\ldots,g_m; n) \in G\wrm{m}\bZ$ be an element of finite order $k\geq2$, where $g_i\in G$, $(i=1,\ldots,m)$, and $n\in\bZ$, so
$g^k = (e,\ldots,e;0)$ is the unit of $G\wrm{m}\bZ$.
Let also $q:G\wrm{m}\bZ \to \bZ$ be the canonical projection.
Then $q(g)=n$, whence $0 = q(e,\ldots,e;0) = q (g^k) = kn$ and $n=0$.
But then $g^k = (g_1,\ldots,g_m; 0)^k = (g_1^k,\ldots,g_m^k; 0) = (e,\ldots,e;0)$, that is each $g_1,\ldots,g_m$ is of order $k$ in $G$.
As $G$ is torsion free, $g_1=\ldots=g_m=e$, whence $g=(e,\ldots,e;0)$ is the unit of $G\wrm{m}\bZ$.
\end{proof}

{\bf Homomorphisms of short exact sequences.}
By an \myemph{exact $(3\times3)$-diagram} we will mean a commutative diagram shown on the left:
\begin{equation}\label{equ:iso_3x3_diagrams}
\aligned
\xymatrix@C=0.8em@R=1em{
\uSeq_0: &
K \ar@{^(->}[d] \ar@{^(->}[r]  & 
L \ar@{^(->}[d]  \ar@{->>}[r]  & 
M \ar@{^(->}[d] &&&&
A\cap L \ar@{^(->}[d] \ar@{^(->}[r]  & 
L \ar@{^(->}[d]  \ar@{->>}[r]  & 
L/(A\cap L) \ar@{^(->}[d] \\ 
\uSeq_1: &
A \ar@{^(->}[r] \ar@{->>}[d] &
B \ar@{->>}[r] \ar@{->>}[d] &
C \ar@{->>}[d] &
\quad\ar[rr]^-{\cong} &&&
A \ar@{^(->}[r] \ar@{->>}[d] &
B \ar@{->>}[r] \ar@{->>}[d] &
B/A \ar@{->>}[d] \\
\uSeq_2: &
P \ar@{^(->}[r] &
Q \ar@{->>}[r]  &
R &&&&
A/(A\cap L) \ar@{^(->}[r] &
B/L \ar@{->>}[r]  &
\frac{B/L}{A/(A\cap L)} \cong \frac{B/A}{M/(A\cap L)}
}
\endaligned
\end{equation}
in which each row and column is a short exact sequence.
In our terminology, such a diagram can be viewed as a \myemph{short exact sequence of its rows (or columns) being in turn short exact sequences of groups homomorphisms}.
If all monomorphisms are inclusions of subgroups, then $L$ and $A$ are normal subgroups of $B$ such that $K=L\cap A$.
Conversely, one easily checks that for every pair of normal subgroups $L$ and $A$ of $B$ there is a unique isomorphism of the above exact $(3\times3)$-diagrams~\eqref{equ:iso_3x3_diagrams} being identity on $B$, $A$, and $L$.

In what follows we will also use the following short exact sequences:
\begin{align}\label{equ:ex_seq_TZ}
\seqTriv& : \{1\}\monoArrow \{1\} \epiArrow \{1\}, &
\seqZ{1}& : \bZ\xmonoArrow{\id} \bZ \epiArrow \{1\}, &
\seqZ{m}& : m\bZ\xmonoArrow{~~} \bZ \xepiArrow{~\mathrm{mod}\, m~} \bZ_m,
\end{align}
for $m\geq2$. 
Then for a short exact sequence $\added{\aSeq}: \kA\xmonoArrow{\alpha} \kB \xepiArrow{~\beta~} \kC$ and $m\geq1$ we have the following exact $(3\times3)$-diagram:
\begin{equation}\label{equ:wr_ex_seq}
\aligned
\xymatrix@C=1.8em@R=1.2em{
\ \aSeq^{m}: \ar@{^(->}[d] &
\kA^m\times0 \ar@{^(->}[d] \ar@{^(->}[r]  & 
\kB^m\times0 \ar@{^(->}[d]  \ar@{->>}[r]  & 
\kC^m\times0 \ar@{^(->}[d] \\ 
\ \seqWrm{\aSeq}{m}: \ar@{->>}[d] &
\kA^m\times m\bZ \ \ar@{^(->}[r]^-{\alpha'} \ar@{->>}[d] &
\kB\wrm{m}\bZ \ar@{->>}[r]^-{\beta'} \ar@{->>}[d]^-{\,p} &
\kC\wr\bZ_m \ar@{->>}[d] \\
\ \seqZ{m}: &
m\bZ \ar@{^(->}[r] &
\bZ \ar@{->>}[r]  &
\bZ_m
}
\endaligned
\end{equation}
where $p$ is the projection to the last coordinate, and 
\begin{gather*}
\alpha'(a_1,\ldots,a_m, mk) = (\alpha(a_1),\ldots,\alpha(a_m), mk), \\
\beta'(b_1,\ldots,b_m, n) = (\beta(b_1),\ldots,\beta(b_m), \MOD{n}{m}),
\end{gather*}
$\alpha_i\in \kA$, $\beta_i\in \kB$ for all $i=1,\ldots,m$, and $k,n\in\bZ$.
The middle horizontal sequence will be called the \myemph{wreath product of $\aSeq$ with $\seqZ{m}$} and denoted by $\seqWrm{\aSeq}{m}$.
Thus~\eqref{equ:wr_ex_seq} is a short exact sequence of its rows: $\aSeq^{m}\monoArrow\seqWrm{\aSeq}{m}\epiArrow\seqZ{m}$.
Evidently,
\begin{gather}\label{equ:W_Z1}
\seqWrm{\aSeq}{1} = \aSeq \times \seqZ{1} : \ \
\kA \times \bZ \ \xmonoArrow{~~~~} \ \kB \wrm{1} \bZ\equiv \kB\times\bZ \ \xepiArrow{~~~~~} \ \kC, \\
\label{equ:Zm_Zn}
\seqWrm{\seqZ{k}}{m}: \ \ (k\bZ)^m \times m\bZ \ \xmonoArrow{~~~} \ \bZ\wrm{m}\bZ \ \xepiArrow{~~~~} \  \bZ_k\wr\bZ_m, \  (k,m\geq1).
\end{gather}

Suppose we have a
\myemph{short exact sequence} of short exact sequences $\kSeq \monoArrow \lSeq \epiArrow \seqZ{m}$:
\begin{equation}\label{equ:3x3_general}
\aligned
\xymatrix@C=0.8em@R=1.1em{
**[l] \kSeq: &
\ \kA \ \ar@{^(->}[d] \ar@{^(->}[rrr]         &&& 
\ \kB \ \ar@{^(->}[d]  \ar@{->>}[rrr]  &&& 
\ \kC \ \ar@{^(->}[d] \\ 
**[l] \lSeq: &
\ \xA \  \ar@{^(->}[rrr] \ar@{->>}[d]         &&&
\ \xB \   \ar@{->>}[rrr]^-{\fc}  \ar@{->>}[d]^-{\hb} &&&
\ \xC \  \ar@{->>}[d] \\
**[l] \seqZ{m}: &
 m\bZ \ar@{^(->}[rrr]                    &&&
\bZ    \ar@{->>}[rrr]^-{\mathrm{mod}\,m} &&&
\bZ_m
}
\endaligned
\end{equation}
in which $\kA, \xA, \kB$ are normal subgroups of $\xB$ and
\begin{align}
\xA &= \ker(\fc) = \hb^{-1}(m\bZ), &
\kB &= \ker(\hb), &
\kA &= \xA \cap \kB.
\end{align}

The following lemma extends to short exact sequences a partial result from the forthcoming paper~\cite{KuznietsovaMaksymenko:Mob:2020} characterizing several kinds of wreath products.

\begin{lemma}\label{lm:charact_seq_wr}{\rm(cf.\cite{KuznietsovaMaksymenko:Mob:2020}).}
In the notation of~\eqref{equ:3x3_general} let $g\in \xB$ and $\zB{0} \subset \kB$ be a subgroup.
Denote $\zA{0} = \kA \cap \zB{0}$ and  $\zC{0} = \zB{0}/\zA{0}$, so we get a short exact sequence $\uSeq: \zA{0} \monoArrow \zB{0} \epiArrow \zC{0}$. 
Also let $\zB{i} := g^{-i} \zB{0} g^{i}$ and $\zA{i} := g^{-i} \zA{0} g^{i}$ for $i=0,\ldots,m-1$.
Then $\zB{i}\subset\kB$ and $\zA{i}\subset\kA$ since $\kB$ and $\kA$ are normal.
Moreover, suppose that
\begin{enumerate}[label={\rm(\alph*)}]
\item\label{enum:3x3:g} 
$\hb(g)=1$ and $g^m$ commutes with $\kB$;

\item\label{enum:3x3:prod} 
$\zB{0},\ldots,\zB{m-1}$ generate $\kB$, pairwise commute, and $\zB{i}\cap \zB{j}=\{e\}$ for all $i\not=j$;

\item\label{enum:3x3:A_i}
$\zA{0},\ldots,\zA{m-1}$ generate $\kA$.
\end{enumerate}
Then the map $\beta:\zB{0}\wrm{m}\bZ \to \xB$ defined by
\begin{align}\label{equ:beta_iso}
\beta(\xv{0},\xv{1},\ldots,\xv{m-1}, k) 
&= \xv{0} \, (g^{-1} \xv{1} g^{1}) \, (g^{-2} \xv{2} g^{2}) \, \cdots \, (g^{-m+1} \xv{m-1} g^{m-1}) g^{k}\\ 
&= \xv{0} \, g^{-1} \, \xv{1}\,  \cdots \, g^{-1} \, \xv{{}m-1} \, g^{-1+m+k}, \nonumber
\end{align}
for $\xv{i}\in \zB{0}$, $i=0,\ldots,m-1$, and $k\in\bZ$, is an isomorphism of groups inducing an isomorphism of exact $(3\times3)$-diagrams:
\begin{equation}\label{equ:iso_3x3_3x3}
\gathered
\xymatrix@C=2em@R=1.3em{
\zA{0}^m\times0 \ar@{^(->}[d] \ar@{^(->}[r]  & 
\zB{0}^m\times0 \ar@{^(->}[d]  \ar@{->>}[r]  & 
\zC{0}^m\times0 \ar@{^(->}[d]                &
&
&
\ \kA \ \ar@{^(->}[d] \ar@{^(->}[r] & 
\ \kB \ \ar@{^(->}[d]  \ar@{->>}[r] & 
\ \kC \ \ar@{^(->}[d]               \\
\zA{0}^m\times m\bZ \ \ar@{^(->}[r]^-{\alpha'} \ar@{->>}[d] &
\zB{0}\wrm{m}\bZ \ar@{->>}[r]^-{\beta'} \ar@{->>}[d]^{p}    &
\zC{0}\wr\bZ_m \ar@{->>}[d]                                 &
\ar@{=>}[r]^-{~~\beta~~}_-{\cong} &
&
\ \xA \  \ar@{^(->}[r] \ar@{->>}[d]                &
\ \xB \   \ar@{->>}[r]^-{\fc}  \ar@{->>}[d]^-{\hb} &
\ \xC \  \ar@{->>}[d]                              \\
 m\bZ \ar@{^(->}[r] &
\bZ    \ar@{->>}[r] &
\bZ_m               &
\ar@{=}[r]^-{~~\id~~} &
&
m\bZ \ar@{^(->}[r] &
\bZ \ar@{->>}[r]   &
\bZ_m
}
\endgathered
\end{equation}
being identity on the lower sequence.
In other words, we get an \myemph{``isomorphism over $\seqZ{m}$''} of short exact sequences
\ $\uSeq^m \monoArrow \seqWrm{\uSeq}{m} \epiArrow \seqZ{m}$ \ and \ $\kSeq \monoArrow \lSeq \epiArrow \seqZ{m}$.
\end{lemma}
\begin{proof}
1) First we will check that $\beta$ induces an isomorphism of upper rows. 
Indeed, let $\delta_i:\zB{0} \to \zB{i}\equiv g^{-i} \zB{0} g^{i}$, $i=0,\ldots,m-1$, be the isomorphism defined by $\delta_i(q)= g^{-i}qg^{i}$.
Then $\delta_i(\zA{0}) = \zA{i}$.
Also notice, that condition~\ref{enum:3x3:prod} just says that $\kB$ splits into a direct product of its subgroups $\zB{i}$, $i=0,\ldots m-1$, i.e. the map $\ia:\zB{0}\times\cdots\times \zB{m-1} \to \kB$ defined by $\ia(\xv{0},\ldots,\xv{m-1}) = \xv{0} \cdots \xv{m-1}$ is an isomorphism.

Now formulas~\eqref{equ:beta_iso} for $\beta$ show that the restriction of $\beta$ to $\zB{0}^m\times0 \equiv \zB{0}^m$ coincides with the composition $\delta:=\alpha \circ (\delta_0\times\cdots\times\delta_{m-1}):\zB{0}^m\to\kB$, and so it is an isomorphism as well.
Moreover, by condition~\ref{enum:3x3:A_i}, the groups $\zA{i}$ generate $\kA$ which implies that 
\[ \beta(\zA{0}^m\times0) = \delta(\zA{0}^m) = \alpha(\zA{0}\times\cdots\times\zA{m-1}) = \kA,\]
whence $\beta$ induces an isomorphism of the upper rows.

2) Next we will show that $\beta$ is an \myemph{isomorphism} such that $\beta(\zA{0}^m\times m\bZ)=\xA$.
First notice that~\ref{enum:3x3:g} and~\ref{enum:3x3:prod} imply that for all $\xv{},\xw{}\in\zB{0}$ and $i,j\in\bZ$
\begin{align}\label{equ:g_act_on_Li}
g^{-i} \, \xv{} \, g^{i} &= g^{\MOD{-i}{m}} \, \xv{} \, g^{\MOD{i}{m}} \in \zB{i}, &
(g^{-i} \, \xv{} \, g^{i})(g^{-j} \, \xw{} \, g^{j})&=
(g^{-j} \, \xw{} \, g^{j})(g^{-i} \, \xv{} \, g^{i}).
\end{align}
Therefore, if $v = (\xv{0},\ldots,\xv{m-1}; k)$ and $w=(\xw{0},\ldots,\xw{m-1}; l) \in \zB{0}\wrm{m}\bZ$, then  
\begin{align*}
\beta(v)\beta(w)&= 
\xv{0} \, (g^{-1} \xv{1} g^{1}) \, (g^{-2} \xv{2} g^{2}) \, \cdots \, (g^{-m+1} \xv{m-1} g^{m-1}) \ g^{k} \ \times   \\
&\qquad \times 
\xw{0}\, (g^{-1} \xw{1} g^{1}) \, (g^{-2} \xw{2} g^{2}) \, \cdots \, (g^{-m+1} \xw{m-1} g^{m-1}) \ g^{l} \ =  \\
&= 
\xv{0} \, (g^{-1} \xv{1} g^{1}) \, (g^{-2} \xv{2} g^{2}) \, \cdots \, (g^{-m+1} \xv{m-1} g^{m-1}) \ \times   \\
&\qquad \times 
(g^{k} \xw{0} g^{-k}) \, (g^{k-1} \xw{1} g^{1-k}) \, (g^{k-2} \xw{2} g^{2-k}) \, \cdots \, (g^{k-m+1} \xw{m-1} g^{-k+m-1}) \, g^{k+l} =\\
&\stackrel{\eqref{equ:g_act_on_Li}}{=}(\xv{0} \xw{k}) \, (g^{-1} \xv{1} \xw{k+1} g^{1}) \, (g^{-2} \xv{2} \xw{k+2} g^{2}) \, \cdots \, (g^{-m+1} \xv{{}m-1+k} g^{m-1}) \ g^{k+l} = \\
&= \beta(\xv{0} \xw{k},\,\xv{1} \xw{k+1},\,\ldots,\, \xv{m-1} \xw{k-1}; k+l) 
\stackrel{\eqref{equ:mult_in_GmZ}}{=} 
\beta(v\,w),
\end{align*}
where all the subscript indices are taken modulo $m$.
Hence $\beta$ is a \myemph{homomorphism}.

Since $\hb(g)=1$, each $x\in \xB$ can be written in the form $x=(x g^{-\eta(x)}) g^{\eta(x)}$, where $x g^{-\eta(x)}\in\kB = \delta(\zB{0}^m)$.
This easily implies that the map $\gamma:\xB \to \zB{0}\wrm{m}\bZ \equiv \zB{0}^m\rtimes\bZ$, $\gamma(x) = \bigl( \delta^{-1}(x g^{-\eta(x)}), g^{\eta(x)}\bigr)$, is the inverse of $\beta$, whence $\beta$ is an \myemph{isomorphism}.

Moreover, it also follows that $\xA = \fc^{-1}(m\bZ)$ is generated by $\kA=\delta(\zA{0}^m)$ and $g^{m}$, whence $\beta(\zA{0}^m\times m\bZ)=\xA$.

Thus $\beta$ induces an isomorphism between left upper $(2\times2)$-squares of diagrams~\eqref{equ:iso_3x3_3x3} whose arrows are monomorphisms.
Hence it yields isomorphisms of all the corresponding quotients, and thus gives an isomorphism of entire diagrams.
Also notice that $p=\hb\circ \beta$:
\[
\hb\circ \beta(\xv{0},\ldots,\xv{m-1}; k)=\hb(g^k)=k = p(\xv{0},\ldots,\xv{m-1}; k),
\]
whence that isomorphism of diagrams is identity on $\bZ$, and therefore on all the lower row.
\end{proof}

\begin{definition}\label{def:solvable_bieberbach}
Recall that a group $G$ is \myemph{solvable} if it has a finite increasing sequence of subgroups $1 = G_0 \vartriangleleft G_1 \vartriangleleft\cdots \vartriangleleft G_n = G$ such that $G_i$ is normal in $G_{i+1}$ and each quotient $G_{i+1}/G_i$ is abelian.
It is well known and is easy to check that in a short exact sequence of groups $A \monoArrow B\epiArrow C$ if any two groups are solvable, then so is the third one.

A short exact sequence $\aSeq:\ A\monoArrow B\epiArrow C$ will be called \myemph{crystallographic} if $A$ is free abelian and $C$ is finite.
In this case the middle group $B$ is also called \myemph{crystallographic}.
If in addition $B$ is torsion free, then $\aSeq$ and $B$ are called \myemph{Bieberbach}.
\end{definition}

The following statement is a ``homotopy'' formulation of one of principal results about Bieberbach groups.
\begin{theorem}[L.~Bieberbach]\label{th:bieberbach}{\rm(e.g.~\cite[Corollary~5.1]{Charlap:BiebGr:1986}).}
Let $\aSeq: \bZ^k \monoArrow B \epiArrow C$ be a Bieberbach sequence.
Then there exists a free action of $C$ on a $k$-torus $T^{k}$ such that the non-zero part $\pi_1 T^k \monoArrow \pi_1(T^{k}/C) \epiArrow C$ of the exact sequence of homotopy groups of the quotient map $p:T^{k} \to T^{k}/C$ (being in this case a covering map) is isomorphic with $\aSeq$.
In particular, $B \cong \pi_1 (T^{k}/C)$.

Moreover, $T^{k}/C$ is, evidently, aspherical, whence every \myemph{aspherical} path-connected topological space $X$ having the homotopy type of a CW-complex and satisfying $\pi_1 X \cong B$ is homotopy equivalent to $T^{k}/C$.
\qed
\end{theorem}

\begin{itemize}[leftmargin=*, wide]
\item 
Let $\classStab$ be the \myemph{minimal} set of isomorphism classes of groups such that:
\begin{align*}
\mathrm{(i)}&~\{1\} \in \classStab, &
\mathrm{(ii)}&~ \text{if $A,B\in\classStab$ and $m\geq 1$ then $A\times B, \, A \wrm{m} \bZ \, \in \, \classStab$;}
\end{align*}

\item
$\classGrp$ the \myemph{minimal} set of isomorphism classes of groups such that:
\begin{align*}
\mathrm{(i)}&~\{1\} \in \classGrp, &
\mathrm{(ii)}&~ \text{if $A,B\in\classGrp$ and $m\geq 1$ then $A\times B, \, A \wr\bZ_m \, \in \, \classGrp$;}
\end{align*}

\item 
$\classZBP$ the set of isomorphisms of \myemph{crystallographic sequences $A\monoArrow B\epiArrow C$} in which $B\in\classStab$, and $C\in\classGrp$;

\item
$\ZBP$ the minimal set of isomorphism classes of \myemph{short exact sequences} satisfying the following conditions:
\begin{align*}
\mathrm{(i)}&~\seqTriv \in \ZBP, &
\mathrm{(ii)}&~ \text{if $\aSeq,\aSeq'\in\ZBP$ and $m\in\bN$ then $\aSeq\times\aSeq', \, \seqWrm{\aSeq}{m} \, \in \, \ZBP$.}
\end{align*}
\end{itemize}

For example, $\bigl((\bZ^4\wrm{3}\bZ) \times (\bZ\wrm{2}\bZ)\bigr) \wrm{7}\bZ\in\classStab$ and  $((\bZ_{12}\times\bZ_{22})\wr\bZ_7) \times \bZ_4\in\classGrp$.

\begin{lemma}\label{lm:classes_SG}
\begin{enumerate}[wide, label={\rm(\arabic*)}]
\item\label{enum:lm:classes_SG:seq_op}
$\seqZ{m} \in \classZBP$ for all $m\geq0$.
Also if $\aSeq_i:  \bZ^{k_i} \monoArrow B_i \epiArrow C_i$, $i=1,2$, are crystallographic (resp. Bieberbach, belonging to $\classZBP$) and $m\geq1$, then so are
\begin{align*}
\aSeq_1\times\aSeq_2 &: \bZ^{k_1+k_2} \monoArrow B_1\!\times\!B_2 \epiArrow C_1\!\times\!C_2, & \ \ \
\seqWrm{\aSeq_1}{m} &: \bZ^{k_1 m}\!\times\! m\bZ \monoArrow B_1\wrm{m}\bZ \epiArrow C_1\wr\bZ_m.
\end{align*}

\item\label{enum:lm:classes_SG:solv_bib}
Every $B\in\classStab$ is solvable and Bieberbach, and every $C\in\classGrp$ is solvable and finite.

\item\label{enum:lm:classes_SG:ZBP}
$\ZBP \subset \classZBP$.
\end{enumerate}
\end{lemma}
\begin{proof}
Statement~\ref{enum:lm:classes_SG:seq_op} is a direct consequence of definitions and we leave it for the reader.
A unique nontrivial point here is a verification that if $B_1$ is torsion free, then $B_1\wrm{m}\bZ$ is so.
But this is proved in Lemma~\ref{lm:no_fin_el}.

\ref{enum:lm:classes_SG:solv_bib}
Minimality of $\classStab$ means that a group $B\in\classStab$ iff $B$ is obtained from the unit group $\{1\}$ by a finite number of operations of two types: $(U,V)\mapsto U\times V$ and $U\mapsto U\wrm{m}\bZ$ for some $m\geq1$.
Since $\{1\}$ is solvable Bieberbach and by 1) those operations preserve solvability and Bieberbach property, $B$ is solvable Bieberbach as well.

Similarly, a group $C\in\classGrp$ iff $C$ is obtained from the unit group $\{1\}$ by a finite number of operations of types $(U,V)\mapsto U\times V$ and $U\mapsto U\wr\bZ_m$ for some $m\geq1$.
Since $\{1\}$ is solvable and finite and those operations preserve solvability and finiteness, it follows that $C$ is solvable and finite as well.

\ref{enum:lm:classes_SG:ZBP}
Again, minimality of $\ZBP$ means that $\aSeq\in\ZBP$ iff $\aSeq$ is obtained from the trivial sequence $\seqTriv$ by a finite number of operations of two types:
$(\uSeq,\vSeq)\mapsto \uSeq\times\vSeq$ and $\uSeq\mapsto \seqWrm{\vSeq}{m}$ for some $m\geq1$.
Since $\seqTriv$ is $\classZBP$ and these operations preserve $\classZBP$-property, it follows that $\aSeq \in \classZBP$ as well.
\end{proof}

\section{$\func$-adapted subsets and their $\func$-regular neighborhoods}
Let $\Xman\subset\Mman$ be a compact submanifold.
If $\Xman$ is connected and $\dim\Xman < \dim\Mman$, then one can define its \myemph{closed tubular} or \myemph{regular} neighborhood, e.g. \cite[Chapter~4]{Hirsch:DiffTop}.
If $\dim\Xman=\dim\Mman$, so $\Xman$ is a subsurface, then by a \myemph{regular} neighborhood of $\Xman$ in $\Mman$ we will mean a union of $\Xman$ with mutually disjoint regular neighborhoods of its boundary components.
A \myemph{regular} neighborhood of an arbitrary submanifold is a union of mutually disjoint regular neighborhoods of its connected components.
\begin{lemma}\label{lm:DiffMX_DiffMV}
Let $\regU{\Xman}$ be a regular neighborhood of $\Xman$, and $Z$\jremoved{9.1}{be} another submanifold disjoint from $\regU{\Xman}$.
Then the inclusion $i:\Diff(\Mman,Z\cup \regU{\Xman}) \subset \Diff(\Mman,Z\cup \Xman)$ is a homotopy equivalence.

In particular, the induced homomophism $i_0:\pi_0\Diff(\Mman,Z\cup \regU{\Xman}) \to \pi_0\Diff(\Mman,Z\cup \Xman)$ is \myemph{injective}\jadded{9.2}{}\footnote{\added{~The homomorphism $i_0$ is of course bijective, however we will use only its injectivity, and therefore explicitly formulate that property.}}, which means that if $\dif \in\Diff(\Mman,Z\cup \regU{\Xman})$ is isotopic to $\id_{\Mman}$ rel.~$Z\cup \Xman$, then $\dif$ is also isotopic to $\id_{\Mman}$ rel.~$Z\cup \regU{\Xman}$.
\hfill\qed
\end{lemma}
This lemma seems to be well-known and follows from ambient isotopy extension theorem for smooth submanifolds and contractibility of the space of tubular neighborhoods of a submanifold $\Xman$ in a manifold $\Mman$, see~\cite[Chapter~A, Proposition~31]{Godin:arxiv:2008}.
The latter fact is an adaptation of the proof of uniqueness theorem for tubular neighborhoods of submanifolds presented in \cite[Chapter~4, Theorem~5.3]{Hirsch:DiffTop}.
In fact, we will use only that last statement of Lemma~\ref{lm:DiffMX_DiffMV} which can be established easier.

Let $\func\in\FF(\Mman,\Pman)$, $c\in \Pman$, and $\Xman$ be a connected component of the corresponding level set $\func^{-1}(c)$.
We will call $\Xman$ \myemph{regular} if it contains no critical points of $\func$.
Otherwise, $\Xman$ will be called \myemph{critical}.
In this case $\Xman$ will also be called \myemph{extremal} whenever it consists of a unique point (being therefore a local extreme of $\func$), and \myemph{non-extremal} otherwise.

Evidently, a regular component of a level set of $\func$ is a submanifold of $\Mman$ diffeomorphic to the circle.
On the other hand, due to Axiom~\AxCrPt, see also Figure~\ref{fig:isol_crit_pt}, a critical component $\Xman$ has a structure of a $1$-dimensional CW-complex whose $0$-cells are critical points of $\func$ belonging to $\Xman$.
Notice that if $\Xman$ contains only quasi-saddles of $\func$, see Figure~\ref{fig:isol_crit_pt}\,b), then it is a smooth submanifold of $\Mman$ diffeomorphic to the circle, however it is nevertheless \myemph{critical} as a connected component of a level set of $\func$.

\begin{definition}\label{def:func_adapted}{\rm(cf.~\cite{Maksymenko:UMZ:ENG:2012})}
Let $\func\in\FF(\Mman,\Pman)$.
A compact submanifold $\Xman \subset\Mman$ will be called \myemph{$\func$-adapted} if each of its connected components is either
\begin{enumerate}
\item[$(1)$]
a regular component of a level-set of $\func$, or
\item[$(2)$]
a compact subsurface whose boundary consists of regular components of some level-sets of $\func$.
\end{enumerate}
In particular, when all connected components of $\Xman$ have dimension $2$, then $\Xman$ will be called an \myemph{$\func$-adapted subsurface}.
In this case the restriction $\func|_{\Xman}:\Xman\to\Pman$ satisfies Axioms~\AxBd\ and \AxCrPt, that is $\func|_{\Xman} \in \FF(\Xman,\Pman)$.
\end{definition}
\begin{example}\rm
A (not necessarily disjoint) union of finitely many $\func$-adapted submanifolds, a connected component of an $\func$-adapted submanifold, and the submanifolds $\varnothing$, $\Mman$, $\partial\Mman$ are $\func$-adapted.
If $[a,b]\subset\Pman$ is a closed segment such that $a$ and $b$ are regular values of $\func$, then $\func^{-1}[a,b]$ is $\func$-adapted.
\end{example}

Let $\Xman\subset\Mman$ be a connected $\func$-adapted submanifold.
Then by an \myemph{$\func$-regular neighborhood} of $\Xman$ we will mean an arbitrary connected $\func$-adapted subsurface $\regU{}$ such that $\regU{}\setminus\Xman$ contains no critical points of $\func$.

More generally, let $\Xman=\mathop{\cup}\limits_{i=1}^{k}\Xman_i$ be a disjoint union of connected $\func$-adapted submanifolds $\Xman_i$.
For each $i=1,\ldots,k$ take any $\func$-regular neighborhood $\regU{i}$ of $\Xman_i$ such that $\regU{i}\cap\regU{j}=\varnothing$ for $i\not=j$.
Then their union $\regU{} = \mathop{\cup}\limits_{i=1}^{k}\regU{i}$, will be called an \myemph{$\func$-regular neighborhood} of $\Xman$.
\jremoved{10.1}{Evidently, each $\func$-regular neighborhood is regular.}

\newcommand\invhp[1]{$(#1,+)$-invariant}
\newcommand\invhm[1]{$(#1,-)$-invariant}

\jadded{38.4}{Let $\dif:\Mman\to\Mman$ be a homeomorphism and $\Vman\subset\Mman$ be an \myemph{orientable} submanifold.
We will say that $\Vman$ is \myemph{\invhp{\dif}} (resp. \myemph{\invhm{\dif}}, if $\dif(\Vman)=\Vman$ and the restriction map $\dif|_{\Vman}:\Vman\to\Vman$ preserves (resp. reverses) orientation of $\Vman$.}

\section{Groups $\FolStabilizer{\func}$ and $\GrpKR{\func}$}\label{sect:FolStab_Grp}
{\bf Graph of $\func\in\FF(\Mman,\Pman)$.}
Consider the partition $\KRGraphf$ of $\Mman$ into connected components of level sets of $\func$, and let $p:\Mman\to\KRGraphf$ be the natural map associating to each $x\in\Mman$ the corresponding element of $\KRGraphf$ containing $x$.
Endow $\KRGraphf$ with the quotient topology, so a subset (a collection of leaves) $A\subset \KRGraphf$ is open iff $p^{-1}(A)$ (i.e. their union) is open in $\Mman$.
It follows from axioms \AxBd\ and \AxCrPt\ that $\KRGraphf$ has a natural structure of $1$-dimensional CW-complex, whose $0$-cells correspond to boundary components of $\Mman$ and critical components of level sets of $\func$.
We will call $\KRGraphf$ the \myemph{graph}%
\footnote{The space $\KRGraphf$ is often called \myemph{Reeb} or \myemph{Kronrod-Reeb} or \myemph{Lyapunov} graph of $\func$, though it was independently introduced by G.~M.~Adelson-Velsky and A.~S.~Kronrod~\cite{Adelson-Welsky-Kronrode:DAN:1945}, and G.~Reeb~\cite{Reeb:CR:1946}, see also \cite{Kronrod:UMN:1950}, \cite{Franks:Top:1985}.}
of $\func\in\FF(\Mman,\Pman)$.
Since $\func$ takes constant values on elements of $\KRGraphf$, it induces a function $\PF{\func}:\KRGraphf\to\Pman$ such that $\func = \PF{\func}\circ p$.

Let $\Homeo(\KRGraphf)$ be the group of homeomorphisms of $\KRGraphf$.
Then each $\dif\in\Stabilizer{\func}$ leaves invariant every level set of $\func$, i.e. $\dif(\func^{-1}(c))=\func^{-1}(c)$ for all $c\in\Pman$, and therefore induces a homeomorphism $\rho(\dif)$ of $\KRGraphf$ making commutative the following diagram:
\begin{equation}\label{equ:2x2_M_Graph}
\aligned
\xymatrix@R=3ex{
\Mman \ar[rr]^-{p} \ar[d]_-{\dif} &&
\KRGraphf \ar[rr]^-{\widehat{\func}} \ar[d]^-{\rho(\dif)} &&
\Pman \ar@{=}[d]  \\
\Mman \ar[rr]^-{p} &&
\KRGraphf \ar[rr]^-{\widehat{\func}} &&
\Pman
}
\endaligned
\end{equation}
so that the correspondence $h\mapsto g(h)$ is a homomorphism of groups $\rho :\Stabilizer{\func} \to \Homeo(\KRGraphf)$.
It can be easily shown, e.g.~\cite[Lemma 2.2]{KravchenkoMaksymenko:EJM:2020}, that its image $\rho(\Stabilizer{\func})$ is a finite subgroup of $\Homeo(\KRGraphf)$.

{\bf Enhanced graph of $\func\in\FF(\Mman,\Pman)$.}
We will now add to $\KRGraphf$ new edges in order to encode certain information coming from \myemph{degenerate local extremes} of $\func$, \cite[Section 6]{Maksymenko:ProcIM:ENG:2010}.

Suppose $\func\in\FF(\Mman,\Pman)$ has a \myemph{degenerate} local extreme $z\in\Mman$ of symmetry index $m$.
Denote by $\StabilizerPlus{\func,z}$ the subgroup of $\Stabilizer{\func}$ consisting of diffeomorphism $\dif$ such that $\dif(z)=z$ and $\dif$ preserves local orientation at $z$.
Then by Lemma~\ref{lm:LStabf} there are local coordinates at $z$ in which for every $\dif\in\StabilizerPlus{\func,z}$ its differential $T_{z}\dif:T_{z}\Mman\to T_{z}\Mman$, regarded as a map $\bR^2\to\bR^2$, is a rotation by some angle $2\pi k/ m$.

\begin{figure}[ht]
\centering
\begin{tabular}{ccc}
\includegraphics[height=1.8cm]{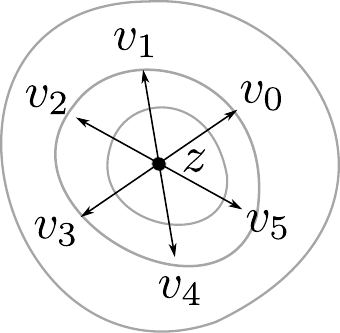}  &  \qquad  &
\includegraphics[height=1.8cm]{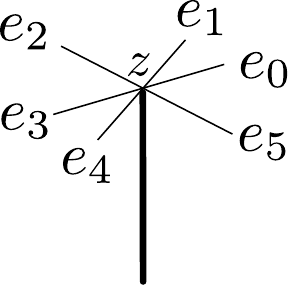}   \\
a) Framing at a degenerate local & & b) Enhanced graph  \\
extreme with symmetry index $m=6$ 
\end{tabular}
\caption{}\label{fig:framings}
\end{figure}

Choose any non-zero tangent vector $v_0\in T_{z}\Mman$ and let $\zfrm{z, v_0} = \{ T_{z}\dif(v_0) \mid \dif\in\StabilizerPlus{\func,z} \}$ be the images of $v_0$ under differentials of all elements from $\StabilizerPlus{\func,z}$.
It follows that $\zfrm{z, v_0}$ consists of $m$ tangent vectors at $z$ which can be enumerated as $v_0,v_1,\ldots,v_{m-1}$, so that for each $\dif\in\Stabilizer{\func,z}$ there exists $k\in\{0,\ldots,m-1\}$ such that $T_{z}\dif(v_{i}) = v_{\MOD{i + k}{m} }$ for all $i=0,\ldots,m-1$.
The set $\zfrm{z, v_0}$ will be called a \myemph{framing} at $z$, see Figure~\ref{fig:framings}a).

Suppose in addition that there exists $\dif\in\Stabilizer{\func,z} \setminus \StabilizerPlus{\func,z}$, i.e. $\dif$ reverses orientation at $z$.
Then $T_z\dif$ is a reflection, and therefore it has a fixed tangent vector $v_0$.
It easily follows that the corresponding framing $\zfrm{z, v_0}$ is then invariant with respect to the differentials of all $\gdif\in\Stabilizer{\func,z}$.

Now let $z_i$, $i=1,\ldots,l$, be all the \myemph{degenerate} local extremes of $\func$, $m_i$ the symmetry index of $z_i$, and $\lfrm_i=\zfrm{z_i,v^i_0}$ some framing at $z_i$.
We will say that such a collection of framings $\lfrm = \{\lfrm_i\}_{i=1,\ldots,l}$ is \myemph{$\func$-adapted}, if it is invariant with respect to $\Stabilizer{\func}$, that is if $\dif\in \Stabilizer{\func}$ and $\dif(z_i) = z_j$ for some $i,j$, then $T_{z_i}\dif(\lfrm_i)=\lfrm_j$.
One easily checks, that $\func$-adapted framings always exists, \cite[Corolary 1]{Maksymenko:ProcIM:ENG:2010}.

Thus $\Stabilizer{\func}$ naturally acts on $\KRGraphf$ as well as on each $\func$-adapted framing $\lfrm$ via the corresponding differentials of maps.
We want to ``join'' these actions.

Regard every point $z_i$, $i=1,\ldots,l$, as a vertex of $\KRGraphf$ of degree $1$, and glue to this vertex $m_i$ edges $e^i_{0},\ldots,e^i_{m_i-1}$, see Figure~\ref{fig:framings}b).
The obtained graph will be denoted by $\EKRGraphf$ and called the \myemph{enhanced} graph of $\func$.

\medskip 

{\bf Action of $\Stabilizer{\func}$ on $\EKRGraphf$.}
Fix a metric on $\EKRGraphf$ in which each edge $e^i_{j}$ has length $1$ and let $\lfrm = \{\lfrm_i\}_{i=1,\ldots,l}$ be any $\func$-adapted framing.
Now for $\dif\in\Stabilizer{\func}$ define a homeomorphism $\widehat{\rho}(\dif):\EKRGraphf\to\EKRGraphf$ by the following rule:
\begin{itemize}
\item 
$\widehat{\rho}(\dif) = \rho(\dif)$ on $\KRGraphf$ regarded a subcomplex of $\EKRGraphf$;
\item
if $T_{z_i}\dif(v^i_{\alpha}) = v^j_{\beta}$, then $\widehat{\rho}(\dif)|_{e^i_{\alpha}}: e^i_{\alpha} \to e^j_{\beta}$ is a unique isometry which maps $z_i$ to $z_j$, $i,j\in\{ 1,\ldots,l\}$, $\alpha\in \{0,\ldots,m_i-1\}$, $\beta\in \{0,\ldots,m_j-1\}$.
\end{itemize}
Then the correspondence $\dif\mapsto\widehat{\rho}(\dif)$ is a well-defined homomorphism $\widehat{\rho}:\Stabilizer{\func}\to\Homeo(\EKRGraphf)$.
Denote
\begin{align*}
\FolStabilizer{\func}&:=\ker(\widehat{\rho}), & \GrpKR{\func}&:= \widehat{\rho}(\Stabilizer{\func}) \cong 
\Stabilizer{\func} \bigl/ \FolStabilizer{\func}.
\end{align*}
One easily check that the following conditions are equivalent:
\begin{itemize}
\item
$\dif\in\ker(\widehat{\rho})=\FolStabilizer{\func}$;
\item
$\dif\in\ker(\rho)$ and the \jadded{10.2}{differential} $T\dif$ fixes every vector from the framing $\lfrm$;
\item
$\dif$ preserves every connected component of every level set of $\func$ and for every \myemph{degenerate local extreme} $z$ the differential $T_{z}\dif: T_z\Mman\to T_z\Mman$ of $\dif$ at $z$ is the identity map.
\end{itemize}

For a closed subset $\Xman$ of $\Mman$ define the following groups:
\begin{align*}
\Diff(\Mman,\Xman) &:= \{\dif\in\Diff(\Mman) \mid \dif \ \text{is fixed on} \ \Xman \}, \\
\DiffNbh(\Mman,\Xman) &:= \{\dif\in\Diff(\Mman) \mid \dif \ \text{is fixed on some neighborhood of} \ \Xman \}, \\
\DiffId(\Mman,\Xman) &:= \{\dif\in\Diff(\Mman,\Xman) \mid \dif \ \text{is isotopic to $\id_{\Mman}$ rel. $\Xman$} \},
\end{align*}
and their intersections with $\FolStabilizer{\func}$ and $\Stabilizer{\func}$:
\begin{equation}\label{equ:subgroups_of_fStabX}
\begin{aligned}
\FolStabilizer{\func,\Xman}&:= \FolStabilizer{\func} \cap \Diff(\Mman,\Xman), & \
\Stabilizer{\func,\Xman} &:= \Stabilizer{\func}\cap \Diff(\Mman,\Xman), \\
\FolStabilizerNbh{\func,\Xman}&:= \FolStabilizer{\func} \cap \DiffNbh(\Mman,\Xman), & \
\StabilizerNbh{\func,\Xman} &:= \Stabilizer{\func}\cap \DiffNbh(\Mman,\Xman), \\
\FolStabilizerIsotId{\func,\Xman}&:= \FolStabilizer{\func} \cap \DiffId(\Mman,\Xman), & \
\StabilizerIsotId{\func,\Xman} &:= \Stabilizer{\func}\cap \DiffId(\Mman,\Xman), \\
\FolStabilizerNbhIsotId{\func,\Xman}&:= \FolStabilizerIsotId{\func} \cap \DiffNbh(\Mman,\Xman), & \
\StabilizerNbhIsotId{\func,\Xman} &:= \StabilizerIsotId{\func}\cap \DiffNbh(\Mman,\Xman),
\end{aligned}
\end{equation}
where we follow the notation convention that $\Xman$ is omitted if it is empty.
In particular, $\FolStabilizerIsotId{\func}=\FolStabilizer{\func} \cap \DiffIdM$.

\begin{lemma}\label{lm:Lf_Sf}{\rm(cf.~\cite[Lemma~3.4]{Maksymenko:AGAG:2006})}
{\rm1)}~Every group in~\eqref{equ:subgroups_of_fStabX} is normal in $\Stabilizer{\func,\Xman}$.

\xadded{1}{0}{33.5}{}
{\rm2)}~The groups $\FolStabilizerIsotId{\func,\Xman}$, $\FolStabilizer{\func,\Xman}$, $\StabilizerIsotId{\func,\Xman}$ are unions of path components of $\Stabilizer{\func,\Xman}$.
In particular, $\StabilizerId{\func,\Xman}$ is the identity path component of each of these groups.
Similarly, the groups $\FolStabilizerNbhIsotId{\func,\Xman}$, $\FolStabilizerNbh{\func,\Xman}$, and $\StabilizerNbhIsotId{\func,\Xman}$ are also normal in $\StabilizerNbh{\func,\Xman}$ and are unions of path connected components of $\StabilizerNbh{\func,\Xman}$.
\end{lemma}
\begin{proof}
1) Normality of the corresponding groups is easy and we leave it for the reader.

2) To prove that $\FolStabilizer{\func,\Xman}$ is a union of path components of $\Stabilizer{\func,\Xman}$, we should check that if $\gamma:I \to \Stabilizer{\func,\Xman}$ is a continuous map with  $\gamma(0)\in\FolStabilizer{\func,\Xman}$, then $\gamma(t)\in\FolStabilizer{\func,\Xman}$ for all $t\in\UInt$.
Moreover, since by assumption $\gamma(t)$ is already in $\Stabilizer{\func,\Xman}$, it suffices to verify that $\gamma(t)\in\FolStabilizer{\func}$.
Then we will get that $\gamma(t)\in\FolStabilizer{\func} \cap \Stabilizer{\func,\Xman} =: \FolStabilizer{\func,\Xman}$.

a) We claim that \myemph{$\gamma(t)$ leaves invariant every connected component $\crLev$ of each level set of $\func^{-1}(c)$, $c\in\Pman$}.
Indeed, since $\func\circ\gamma(t)=\func$, it follows that $\gamma(t)(\func^{-1}(c))=\func^{-1}(c)$.
But $\crLev$ is an open-closed subset of $\func^{-1}(c)$ and $\gamma(0)(\crLev)=\crLev$, therefore $\gamma(t)(\crLev)=\crLev$ for all $t$ due to the continuity of $\gamma$.

b) Suppose $z$ is a degenerate local extreme of symmetry index $m$.
\jadded{11.2}{}%
Then for each $t\in\UInt$ the \xchanged{-0.8}{0}{11.1}{tangent maps}{differential} $T_z\gamma(t)$ is a rotation by the angle $2\pi k/m$ for some $k\in\{ 0,\ldots,m-1 \}$, and at the same time, cf. Remark~\ref{rem:deformation_in_fStab}, it must be continuous in $t$.
Hence $T_z\gamma(t) = T_z\gamma(0)=\id_{T_z\Mman}$.

It follows from a) and b) that $\gamma(t)\in\FolStabilizer{\func}$ for all $t\in\UInt$.
As noted above this implies that $\FolStabilizer{\func,\Xman}$ is a union of path components of $\Stabilizer{\func,\Xman}$.

Since $\DiffIdMX$ is a path component of $\DiffMX$, it follows that the intersections $\StabilizerIsotId{\func,\Xman}$ and $\FolStabilizerIsotId{\func,\Xman}$ with $\DiffIdMX$ are unions of path components of $\Stabilizer{\func,\Xman}$ as well.
The proofs for $\FolStabilizerNbhIsotId{\func,\Xman}$, $\FolStabilizerIsotId{\func,\Xman}$, and $\StabilizerNbhIsotId{\func,\Xman}$ are similar.
\end{proof}
\added{By} Lemma~\ref{lm:Lf_Sf}\added{,}\removed{implies that} $\pi_0\FolStabilizer{\func,\Xman}$ can be regarded as a normal subgroup of $\pi_0\Stabilizer{\func,\Xman}$.
\changed{Consider the corresponding factor group}{Put}:
\begin{align}\label{equ:GfX}
\GrpKR{\func,\Xman}&:=\dfrac{\Stabilizer{\func,\Xman}}{\FolStabilizer{\func,\Xman}} =
\Bigl. \frac{\Stabilizer{\func,\Xman}}{\StabilizerId{\func,\Xman}} \Bigr/
\frac{\FolStabilizer{\func,\Xman}}{\StabilizerId{\func,\Xman}}
= \dfrac{\pi_0\Stabilizer{\func,\Xman}}{\pi_0\FolStabilizer{\func,\Xman}}.
\end{align}
\added{
In particular, we get the following short exact sequence:
\[
\seqStab{\func,\Xman}: \
\pi_0\FolStabilizer{\func,\Xman} \monoArrow
\pi_0\Stabilizer{\func,\Xman} \epiArrow
\GrpKR{\func,\Xman}.
\]
One can define is a similar way the groups
\begin{align*}
\GrpKRIsotId{\func,\Xman}&:=\dfrac{\StabilizerIsotId{\func,\Xman}}{\FolStabilizerIsotId{\func,\Xman}}, & \ \
\GrpKRNbh{\func,\Xman}&:=\dfrac{\StabilizerNbh{\func,\Xman}}{\FolStabilizerNbh{\func,\Xman}}, \ \ &
\GrpKRNbhIsotId{\func,\Xman}&:=\dfrac{\StabilizerNbhIsotId{\func,\Xman}}{\FolStabilizerNbhIsotId{\func,\Xman}},
\end{align*}
and short exact sequences $\seqStabIsotId{\func,\Xman}$, $\seqStabNbh{\func,\Xman}$, $\seqStabNbhIsotId{\func,\Xman}$.
}

Notice that for any two closed subsets $\Xman\subset\Yman$ of $\Mman$ \added{we have} the following commutative diagram of homomorphisms:
\begin{equation}\label{equ:bX_bY_seq}
\begin{gathered}
\xymatrix@R=3ex{
\FolStabilizer{\func,\Yman} \ar@{^(->}[r] \ar@{^(->}[d]  & \Stabilizer{\func,\Yman} \ar@{^(->}[d] \ar@{->>}[r]
& \GrpKR{\func,\Yman} \ar@{^(->}[d]^{\GHom{\Yman,\Xman}} \\
\FolStabilizer{\func,\Xman} \ar@{^(->}[r] & \Stabilizer{\func,\Xman} \ar@{->>}[r]
& \GrpKR{\func,\Xman}
}
\end{gathered}
\end{equation}
whose rows are exact, the left square consists of natural inclusions, and $\GHom{\Yman,\Xman}$ is the induced homomorphism of the corresponding quotients.
Since $\FolStabilizer{\func,\Yman} = \FolStabilizer{\func,\Xman} \cap \Stabilizer{\func,\Yman}$, it follows from isomorphism~\eqref{equ:iso_3x3_diagrams} that $\GHom{\Yman,\Xman}$ is a \myemph{monomorphism} as well.
In particular, diagram~\eqref{equ:bX_bY_seq} induces a homomorphism $\seqStab{\func,\Yman} \monoArrow \seqStab{\func,\Xman}$ of short exact sequences.

\medskip 

\xadded{0}{0}{13.1}{
{\bf$\func$-adapted Dehn twists.}
Let $\gamma$ be a regular component of some level set of $\func$ and $\Xman$ its $\func$-regular neighborhood.
Put $J = [-1,1]$ if $\gamma\subset\Int{\Mman}$ and $J = \UInt$ if $\gamma\subset\partial\Mman$.
Let also $\Cylinder = \Circle\times J$ and $\Bman = \Circle\times[0.2,0.8]$.
Then there exists a diffeomorphism $\phi:\Cylinder \to \Xman$ such that $\phi(\Circle\times 0) = \gamma$ and $\phi(\Circle\times t)$, $t\in J$, is also a regular component of some level set of $\func$.

Fix any $\Cinfty$-function $\mu:J\to\UInt$ such that $\mu=0$ on $J \cap [-1,0.2]$ and $\mu=1$ on $[0.8,1]$, and define the following diffeomorphism $\bar{\tau}:\Cylinder\to\Cylinder$ by $\tau(z,t) = (z e^{2\pi i \mu(t)}, t)$.
Since $\tau$ is fixed near $\partial\Cylinder$, it induces a diffeomorphism $\tau_{\gamma}:\Mman\to\Mman$ given by 
\begin{equation}\label{equ:Dehn_twist}
\tau(x) = \begin{cases}
\phi\circ\bar{\tau} \circ \phi^{-1}(x), & x\in \Xman,\\
x, & x\in \Mman\setminus \Xman.
\end{cases}
\end{equation} 
It follows from the construction that $\tau \in \FolStabilizer{\func,\overline{\Mman\setminus\phi(\Bman)} } \subset \FolStabilizer{\func} \subset \Stabilizer{\func}$.
In particular, $\tau$ is fixed near $\gamma$.

In spite that $\tau$ depends on the choice of an $\func$-regular neighborhood $\Xman$ of $\gamma$ as well as on the diffeomorphism $\phi$, we will call any such $\tau$ an \myemph{$\func$-adapted Dehn twist along $\gamma$}.
}

\section{Main results}\label{sect:main_results}
Suppose $\Mman$ is orientable and distinct from $S^2$ and torus $T^2$.
Let $\func\in\FF(\Mman,\Pman)$ and $\Vman\subset\Mman$ be an $\func$-adapted submanifold.
Our aim is to describe up to isomorphism the \changed{following short exact sequence}{corresponding $\DSG$-sequence} 
\begin{equation}\label{equ:exact_seq_DSG}
\seqStabIsotId{\func,\Vman}: \ \ \
   \pi_0\FolStabilizerIsotId{\func,\Vman} \ \monoArrow \
    \pi_0\StabilizerIsotId{\func,\Vman} \ \epiArrow \
    \GrpKRIsotId{\func,\Vman},
\end{equation}
which will also be called \myemph{Bieberbach sequence} for $(\func,\Vman)$.

We will prove that~\eqref{equ:exact_seq_DSG} can be expressed in an \myemph{explicit} way via the corresponding \changed{groups}{$\DSG$-sequences} of the restriction of $\func$ onto certain subsurfaces $\{\Xman_i\}$ being only $2$-disks (Lemmas~\ref{lm:reduction_Mconn_V_dM}, \ref{lm:f_null_hom_repl_S1_R1}, and Theorems~\ref{th:stab:chi_neg}, \ref{th:stab:annulus}, \ref{th:stab:disk:one_crpt}, \ref{th:stab:disk_ann:gen_case}).
Moreover, those disks are \added{also} defined explicitly and the above results can therefore be applied to each of the restrictions $\func|_{\Xman_i}$ as well.
Notice that each of those restriction has smaller number of critical points than $\func$.
Hence for computation of  \added{$\DSG$-sequence} \eqref{equ:exact_seq_DSG} one can use an induction on the total number $|\fSing|$ of critical points of $\func$, see Example~\ref{exmp:computaiton} below.
The initial inductive step is provided by Theorems~\ref{th:stab:annulus} and \ref{th:stab:disk:one_crpt} describing the case of functions with minimal number of critical points on $\Disk$ and $\Circle\times\UInt$.

In a series of papers~\cite{MaksymenkoFeshchenko:UMZ:ENG:2014, MaksymenkoFeshchenko:MFAT:2015, MaksymenkoFeshchenko:MS:2015, Feshchenko:Zb:2015} the case of torus was reduced to the cases $\Mman=\Disk$ or $\Circle\times\UInt$.
\xadded{-0.8}{0}{12.1}{Also, due to Theorem~\ref{th:stab:chi_neg}, computation of \changed{(5.1)}{$\DSG$-sequence} for all (orientable and non-orientable) compact surfaces $\Mman$ with $\chi(\Mman)<0$ reduces to the case of $2$-disk, annulus and M\"obius band.}
Thus the remained open problems are to compute \changed{(5.1)}{$\DSG$-sequences} for functions on $2$-sphere, projective plane, M\"obius band, and Klein bottle.

\begin{lemma}[Reduction to the case of non-empty boundary]\label{lm:incl_SprfX_Sprf}
Let $\Mman$ be a connected compact surface, $\func\in\FF(\Mman,\Pman)$, $\Vman$\jremoved{12.2}{be} a non-empty union of several boundary components of $\Mman$, and $\jInclZ:\pi_0\StabilizerIsotId{\func,\Vman} \to \pi_0\StabilizerIsotId{\func}$\jremoved{12.3}{be} the homomorphism induced by the inclusion $\jIncl:\StabilizerIsotId{\func,\Vman} \subset \StabilizerIsotId{\func}$.
\xadded{-0.5}{0}{12.4}{Suppose also that $\StabilizerId{\func}$ is contractible}.
\begin{enumerate}[wide, label={\rm(\roman*)}]
\item\label{enum:xx:DOS__jSdS}
Then we have the following commutative diagram:
\xadded{5}{0}{28.3}{}
\begin{equation}\label{equ:DOS__jSdS}
\begin{aligned}
\xymatrix{
   \pi_1\DiffId(\Mman) \ar@{^(->}[r]^-{p_1} \ar@/_2.8pc/[ddrrr]_-{\alpha}  &
   \pi_1\OrbitPathComp{\func}{\func}  \ar@{->>}[rr]^-{\partial} \ar[drr]^-{\beta} &&
   \pi_0\StabilizerIsotId{\func}   \\
   & 
   \pi_1\OrbitPathComp{\func,\Vman}{\func}
   \ar[u]_-{\cong}^{\text{Theorem~\ref{th:stab_orb:full_info}\ref{enum:th:stab_orb_full_info:Off_V}}} 
   \ar@{->>}[rr]^-{\cong}_-{\eqref{equ:pi1OfX_pi0SfX}} &&
   \pi_0\StabilizerIsotId{\func,\Vman} \ar@{->>}[u]_-{\jInclZ} \\
    & && \ker(\jInclZ) \ar@{^(->}[u] 
}
\end{aligned}
\end{equation}
In particular, it induces an isomorphism between the right column and the upper row which coincides with the sequence~\eqref{equ:exact_seq_for_pi1OfX}.

\item\label{enum:xx:3x3_diagram}
We also have the following \jadded{13.2}{exact $(3\times3)$-diagram:} 
\begin{equation}\label{equ:3x3_diagram}
\begin{gathered}
\xymatrix@R=3ex{
& \ker(\jInclZ) \ \ar@{^{(}->}[d]^(.25){} \ar@{=}[r]         &
\ \ker(\jInclZ) \ \ar@{^{(}->}[d] \ar[r]                   &
\ \{1\} \ \ar@{^{(}->}[d] \\
**[l] \seqStabIsotId{\func,\Vman}: &
\pi_0\FolStabilizerIsotId{\func,\Vman} \ \ar@{->>}[d] \ar@{^{(}->}[r] &
\ \pi_0\StabilizerIsotId{\func,\Vman} \ \ar@{->>}[d]^{\jInclZ} \ar@{->>}[r]^-{\fc_\Vman} &
\ \GrpKRIsotId{\func,\Vman} \ar[d]_-{\cong}^-{\GHomIsotId{\Vman,\varnothing}} \\
**[l] \seqStabIsotId{\func}: &
\pi_0\FolStabilizerIsotId{\func} \ \ar@{^{(}->}[r]              &
\ \pi_0\StabilizerIsotId{\func} \ \ar@{->>}[r]^-{\fc} &
\ \GrpKRIsotId{\func}
}
\end{gathered}
\end{equation}

\item\label{enum:xx:chiM}
If $\chi(\Mman)<0$, and thus $\DiffId(\Mman)$ is contractible, then by~\eqref{equ:DOS__jSdS}, $\jInclZ$ is an isomorphism, and~\eqref{equ:3x3_diagram} yields an isomorphism  $\seqStabIsotId{\func,\Vman}\cong \seqStabIsotId{\func}$.
On the other hand, if $\chi(\Mman)\geq0$, then $\Mman$ is either a $2$-disk, or an annulus, or a M\"obius band, $\ker(\jInclZ) \cong \pi_1\DiffId(\Mman)\cong\bZ$ due to~\eqref{equ:DOS__jSdS}, and the diagram \eqref{equ:3x3_diagram} is a short exact sequence $\seqZ{1} \monoArrow \seqStabIsotId{\func,\Vman} \epiArrow \seqStabIsotId{\func}$.
\end{enumerate}
\end{lemma}
\begin{proof}
\jadded{29.1-29.7}{}
\ref{enum:xx:DOS__jSdS}
First notice that the right upper square of~\eqref{equ:DOS__jSdS} is commutative due to functorial properties of short exact sequences of homotopy groups of fibrations. 
Moreover, by assumption $\Vman\not=\varnothing$, whence the group $\DiffId(\Mman,\Vman)$ is contractible, and we get isomorphism~\eqref{equ:pi1OfX_pi0SfX}. 
Then $\jInclZ$ is surjective due to surjectivity of $\partial$.
Finally, $\beta$ is just the compositions of the corresponding isomorphisms, and $\alpha$ is induced by $\beta$.

\ref{enum:xx:3x3_diagram}
The homomorphism from the middle row to the lower one is induced by the inclusion of pairs
$\jIncl: \bigl(\StabilizerIsotId{\func,\Vman}, \FolStabilizerIsotId{\func,\Vman}\bigr) \subset
 \bigl(\StabilizerIsotId{\func}, \FolStabilizerIsotId{\func}\bigr)$, see~\eqref{equ:bX_bY_seq}.
Since $\jInclZ$, $\fc_\Vman$, and $\fc$ are surjective it follows that $\GHomIsotId{\Vman,\varnothing}$ is surjective as well.
Moreover, as $\FolStabilizerIsotId{\func,\Vman} = \FolStabilizerIsotId{\func} \cap \StabilizerIsotId{\func,\Vman}$, it also follows that $\jInclZ$ maps $\pi_0\FolStabilizerIsotId{\func,\Vman}$ onto $\pi_0\FolStabilizerIsotId{\func}$ and that $\GHomIsotId{\Vman,\varnothing}$ is injective as well, i.e.\! it is an isomorphism.
This implies all the diagram~\eqref{equ:3x3_diagram}.

Statement \ref{enum:xx:chiM} is a direct consequence of~\ref{enum:xx:DOS__jSdS} and~\ref{enum:xx:3x3_diagram}.
\end{proof}

\removed{
Lemma 5.1. will be proved in~\S*. It suggests to describe ...}
The following lemma reduces the problem of computing $\DSG$-sequences to the case when $\Mman$ is connected.
It will be proved in~\S\ref{sect:proof:lm:reduction_Mconn_V_dM}.
\begin{lemma}[Reduction to connected surfaces]\label{lm:reduction_Mconn_V_dM}
\xadded{0.2}{0}{13.3}{}
Let $\Mman$ be a compact possibly non-connected surface, $\func\in\FF(\Mman,\Pman)$, $\Vman$\jremoved{13.4}{be} an $\func$-adapted submanifold, $\regN{}$\jremoved{13.5}{be} an $\func$-regular neighborhood of $\Vman$, $\Xman_1,\ldots,\Xman_{\cnt}$ all the connected components of $\overline{\Mman \setminus \regN{}}$, and $\hXman_i = \Xman_i \cap \regN{}$, so $\hXman_i$ consists of some boundary components of $\Xman_i$.
Then
\added{there exists an isomorphism $\seqStab{\func,\Vman} \cong \prod\limits_{i=1}^{\cnt} \seqStab{\func|_{\Xman_i},\hXman_i}$ inducing an isomorphism $\seqStabIsotId{\func,\Vman} \cong \prod\limits_{i=1}^{\cnt} \seqStabIsotId{\func|_{\Xman_i},\hXman_i}$.
Moreover, denote by $\zB{i}$, $i=1,\ldots,\cnt$, the subgroup of $\pi_0\Stabilizer{\func,\Vman}$ corresponding to the group
\[ [\id_{\Xman_1}]\times\cdots\times\pi_0\Stabilizer{\func|_{\Xman_i},\hXman_i} \times\cdots\times[\id_{\Xman_\cnt}]\] under above isomorphism, so $\pi_0\Stabilizer{\func,\Vman}$ splits into the direct product $\zB{1}\times\cdots\times\zB{\cnt}$.
Let also $\dif\in\Stabilizer{\func,\Vman}$ be such that $\dif(\Xman_i)=\Xman_j$ for some $i,j$.
Then $[\dif] \zB{i} [\dif]^{-1} = \zB{j}$ in $\pi_0\Stabilizer{\func,\Vman}$.
}
\end{lemma}

\begin{lemma}[Reduction to null-homotopic maps]\label{lm:f_null_hom_repl_S1_R1}
\xadded{1}{0}{13.7-13.9}{}
Let $\Mman$ be a connected compact surface, $p:\tPman\to \Pman$ a smooth covering map between connected $1$-dimensional manifolds each being either $\bR$ or $\Circle$, $\tfunc:\Mman\to\tPman$ a smooth map, and $\func=p\circ\tfunc:\Mman\to\Pman$.
Then $\func\in\FF(\Mman, \Pman)$ if and only if $\tfunc\in\FF(\Mman, \tPman)$.
Moreover, in this case 
\begin{enumerate}[label=$(\arabic*)$]
\item\label{enum:lm:nullhom:f-adapted}
a submanifold $\Xman \subset \Mman$ is $\func$-adapted if and only if $\Xman$ is $\tfunc$-adapted;

\item\label{enum:lm:nullhom:lift_in_F}
$\FolStabilizer{\func} = \FolStabilizer{\tfunc}$ and $\Stabilizer{\func} \supset \Stabilizer{\tfunc}$;

\item\label{enum:lm:nullhom:R}
if $\tPman=\bR$, then $\Stabilizer{\func} = \Stabilizer{\tfunc}$.
\end{enumerate}
In particular, if $\Mman$ is one of the surfaces: $2$-disk, annulus, $2$-sphere, projective plane or M\"obius band\added{, then for computations of $\DSG$-sequences} one can assume that $\Pman=\bR$.
\end{lemma}
\begin{proof}
\xadded{1}{0}{13.6}{}
Since $p:\tPman\to\Pman$ is a covering map, the relation $\func = p \circ\tfunc$ implies that 
\begin{enumerate}[label={\rm(\alph*)}]
\item\label{enum:nulhom:bd} $\func$ is constant on each boundary component $\Vman$ of $\Mman$ and has no critical points on $\partial\Mman$ iff $\tfunc$ does so;
\item\label{enum:nulhom:crpt} $\func$ and $\tfunc$ have the same critical points and they are smoothly equivalent at each of those points;
\item\label{enum:nulhom:reg_comp} $\func$ and $\tfunc$ have the same partitions into connected components of level sets.
\end{enumerate}
Then~\ref{enum:nulhom:bd} and~\ref{enum:nulhom:crpt} mean that $\func$ satisfies axioms \AxBd\ and \AxCrPt\ if and only if $\tfunc$ does so.
In other words, $\func\in\FF(\Mman, \Pman)$ if and only if $\tfunc\in\FF(\Mman, \tPman)$.

Assume now that $\func\in\FF(\Mman, \Pman)$.
Then~\ref{enum:lm:nullhom:f-adapted} and the relation $\FolStabilizer{\func} = \FolStabilizer{\tfunc}$ follows from~\ref{enum:nulhom:reg_comp}.
Moreover, if $\dif\in\Stabilizer{\tfunc}$, that is $\tfunc\circ\dif=\tfunc$, then $\func\circ\dif=p \circ \tfunc\circ \dif =  p \circ \tfunc = \func$, i.e.\! $\dif\in\Stabilizer{\func}$.
In other words, $\Stabilizer{\tfunc} \subset \Stabilizer{\func}$, which completes~\ref{enum:lm:nullhom:lift_in_F}.

\ref{enum:lm:nullhom:R} 
Let $\tPman = \bR$ and $\dif\in\Stabilizer{\func}$.
We should prove that $\dif\in\Stabilizer{\tfunc}$ as well, i.e. $\tfunc\circ\dif = \tfunc$.

If $\Pman=\bR$, then $p$ is a diffeomorphism and our lemma is evident.
Suppose $\Pman=\Circle$.
Then one can assume that $p:\bR\to \Circle$ is the universal covering map given by $p(t)=e^{2\pi i t}$.
Notice that the assumption $\func\circ\dif=\func$ means that
\[ 
p\circ \tfunc \circ \dif(x) = 
e^{2\pi i \tfunc (\dif(x))} = 
e^{2\pi i \tfunc (x)} = 
p \circ \tfunc (x),
\]
for each $x\in\Mman$, whence $\tfunc\circ\dif(x) = \tfunc(x) + k(x)$ for some $k(x)\in\bZ$.
In other words, $k(x) = \tfunc\circ\dif(x) - \tfunc(x)$ is a continuous function taking only integer values, and therefore it is constant.
Hence $\tfunc\circ\dif = \tfunc + k$ for some $k\in\bZ$ which does not depend on $x\in\Mman$.

As $\Mman$ is connected, $\tfunc(\Mman) = [a,b] \subset \bR$ for some $a,b\in\bR$.
Then $\tfunc\circ\dif(\Mman) = \tfunc(\Mman) = [a,b]$ as well.
On the other hand, $\tfunc\circ\dif = \tfunc + k$ implies that $\tfunc\circ\dif(\Mman) = [a+k, b+k]$, whence $k=0$, and so $\tfunc\circ\dif=\tfunc$.
\end{proof}

{\bf Surfaces of negative Euler characteristics.}
Theorem~\ref{th:stab:chi_neg} below \changed{furhter}{further} reduces description of \changed{(5.1)}{$\DSG$-sequences} to the case when $\chi(\Mman)\geq0$.
It extends~\cite[Theorem~1.7]{Maksymenko:MFAT:2010} and we will sketch its proof in~\S\ref{sect:proof:th:pi0SfX:decomp}.

Let $\Xman$ be an $\func$-adapted submanifold, $\regN{\Xman}$ be its $\func$-regular neighborhood and $D_1,\ldots,D_q$ be all the connected components of $\overline{\Mman\setminus\regN{\Xman}}$ diffeomorphic with a $2$-disk.
Then the union
\begin{equation}\label{equ:canon_nbh}
\canN{\Xman}: = \regN{\Xman} \cup D_1\cup\ldots\cup D_q
\end{equation}
will be called a \myemph{canonical neighborhood} $\canN{\Xman}$ of $\Xman$ (corresponding to $\regN{\Xman}$), see~\cite{JacoShalen:Topology:1977}.

\begin{theorem}\label{th:stab:chi_neg}{\rm\jadded{14.1}{(cf.\,\cite{Maksymenko:MFAT:2010}).}}
Let $\Mman$ be a compact connected surface (possibly non-orientable) with $\chi(\Mman)<0$, $\func\in\FF(\Mman,\Pman)$, and $\Vman\subset\partial\Mman$\jremoved{14.2}{be} an $\func$-adapted submanifold.
Let also $\crLev$ be the union of all non-extremal critical components of level sets of $\func$ whose canonical neighborhoods have negative Euler characteristic, see~\eqref{equ:canon_nbh}, $\regN{}$ \jremoved{14.3}{be} an $\func$-regular neighborhood of $\crLev$, $\Xman_1,\ldots,\Xman_{\cnt}$ \jremoved{14.4}{be} all the connected components of $\overline{\Mman\setminus\regN{}}$, and $\hXman_i = \Xman_i \cap (\Vman\cup\regN{})$, $i=1,\ldots,\cnt$.
Then
\begin{enumerate}[leftmargin=5ex, topsep=0pt, label=$(\arabic*)$]
\item\label{enum:pi0SfX:decomp:1}
each $\Xman_i$ is diffeomorphic either with a $2$-disk or an annulus or a M\"obius band;
\item\label{enum:pi0SfX:decomp:2}
the inclusion $\bigl( \StabilizerIsotId{\func,\Vman\cup\regN{}}, \FolStabilizerIsotId{\func,\Vman\cup\regN{}} \bigr) \subset \bigr(\StabilizerIsotId{\func, \Vman}, \FolStabilizerIsotId{\func, \Vman}\bigr)$ is a homotopy equivalence, implying therefore an isomorphism $\seqStabIsotId{\func,\Vman} \cong \seqStabIsotId{\func,\Vman\cup\regN{}}
\cong \myprod\limits_{i=1}^{\cnt}\seqStabIsotId{\func|_{\Xman_i},\hXman_i}$.
\end{enumerate}
\end{theorem}

{\bf Functions on $2$-disk and annulus.}
Suppose $(\Mman,\Vman)$ is one of the pairs $(\Disk,\partial\Disk)$ or $(\Circle\times\UInt, \Circle\times 0)$. 
Then the group $\Diff(\Mman, \Vman)$ is connected (in fact even contractible, see~\cite{Smale:ProcAMS:1959, Gramain:ASENS:1973} and the table in the proof of Theorem~\ref{th:stab_orb:full_info}\ref{enum:th:stab_orb_full_info:pikOf}.
Hence $\Diff(\Mman, \Vman) = \DiffId(\Mman, \Vman)$ and 
\begin{equation}\label{equ:StabIsotIdD2}
 \StabilizerIsotId{\func,\Vman} = \Stabilizer{\func,\Vman} \cap\DiffId(\Mman, \Vman)= \Stabilizer{\func,\Vman} \cap\Diff(\Mman, \Vman) = \Stabilizer{\func,\Vman}.
\end{equation}
Similarly, $\FolStabilizerIsotId{\func,\Vman} = \FolStabilizer{\func,\Vman}$.
This implies that $\seqStab{\func,\Vman} = \seqStabIsotId{\func,\Vman}$.

\begin{theorem}[Functions on the annulus]\label{th:stab:annulus}
Let $(\Cylinder,\Vman) = (\Circle\times\UInt, \Circle\times0)$.
\begin{enumerate}[wide,label={\rm(\arabic*)}]
\item\label{eqnu:th:cyl:func:no_cr_pt}
If $\func\in\FF(\Cylinder,\Pman)$ has \myemph{no critical points}, then
\begin{enumerate}[leftmargin=11ex, label={\rm(\alph*)}]
\item\label{enum:th:cyl:func:no_cr_pt:Spr}
$\pi_0\StabilizerIsotId{\func,\Vman}$ is trivial, whence $\sDSG{\func,\Vman} \cong \seqTriv: \{1\} \monoArrow \{1\} \epiArrow \{1\}$;
\item\label{enum:th:cyl:func:no_cr_pt:S}
$\pi_0\Stabilizer{\func,\partial\Cylinder}=\bZ$.
\end{enumerate}

\item\label{eqnu:th:cyl:func:incl}
For each $\func\in\FF(\Cylinder,\Pman)$ the inclusion\changed{s}{ of pairs}
\begin{align}\label{equ:func_on_cylinder_hom_eq1}
\bigl(
\StabilizerIsotId{\func, \partial\Cylinder},\,  \FolStabilizerIsotId{\func, \partial\Cylinder}
\bigr)
\ \subset \
\bigl( \StabilizerIsotId{\func, \Vman},\, \FolStabilizerIsotId{\func, \Vman} \bigr) 
\ \stackrel{\eqref{equ:StabIsotIdD2}}{\equiv} \ 
\bigl( \Stabilizer{\func, \Vman},\, \FolStabilizer{\func, \Vman} \bigr)
\end{align}
\changed{are}{is a} homotopy equivalence\removed{s}, whence \added{$\seqStabIsotId{\func,\partial\Cylinder} \cong \seqStabIsotId{\func,\Vman} \equiv \seqStab{\func,\Vman}$.}
\removed{exact sequences ...}

\item\label{eqnu:th:cyl:func:subset}
Let $\Mman$ be a not necessarily orientable connected compact surface, $\Vman$ a connected component of $\partial\Mman$, and $\func\in\FF(\Mman,\Pman)$.
Suppose there exists a regular component $\Wman$ of some level set of $\func$ separating $\Mman$ and let $\Bman$ and $\Cylinder$ be the connected components of $\overline{\Mman\setminus\Wman}$.
Assume that $\Vman\subset\Cylinder$ and let \jchanged{14.5}{$\Xman\subset \Cylinder\setminus \Wman$}{$\Xman\subset \Bman\setminus \Wman$} be an $\func$-adapted submanifold, see Figure~\ref{fig:cylinder_split}.
Suppose $\Cylinder$ is an \myemph{annulus} 
\jadded{35.10}{and $\dif(\Cylinder)=\Cylinder$ for all $\dif\in\Stabilizer{\func,\Xman\cup\Vman}$.}
Then we have the following homotopy equivalence:
\[
\Stabilizer{\func, \Xman \cup \Vman} \cong \Stabilizer{\func|_{\Bman}, \Xman \cup \Wman} \times \StabilizerIsotId{\func|_{\Cylinder}, \partial\Cylinder}
\]
which induces the homotopy equivalence
$\StabilizerIsotId{\func, \Xman \cup \Vman} \cong \StabilizerIsotId{\func|_{\Bman}, \Xman \cup \Wman} \times \StabilizerIsotId{\func|_{\Cylinder}, \partial\Cylinder}$ and homotopy equivalences between the corresponding $\Delta$- and $\Delta'$-groups.
In particular, we get isomorphisms of short exact sequences:
\begin{equation}\label{equ:cyl:split}
\begin{aligned}
\seqStab{\func,\Xman\cup\Vman}&\cong \seqStab{\func|_{\Bman}, \Xman \cup \Wman} \times \seqStabIsotId{\func|_{\Cylinder}, \partial\Cylinder}, \\
\seqStabIsotId{\func,\Xman\cup\Vman} &\cong \seqStabIsotId{\func|_{\Bman}, \Xman \cup \Wman} \times \seqStabIsotId{\func|_{\Cylinder}, \partial\Cylinder}.
\end{aligned}
\end{equation}
In particular, if $\Bman$ is either a $2$-disk or an annulus and $\Xman=\varnothing$, then
\begin{equation}\label{equ:bseq_two_cyl}
\seqStabIsotId{\func,\partial\Mman} \cong 
\seqStabIsotId{\func|_{\Bman}, \partial\Bman} \times 
\seqStabIsotId{\func|_{\Cylinder}, \partial\Cylinder}.
\end{equation}
\begin{figure}[ht]
\centering
\includegraphics[height=1.8cm]{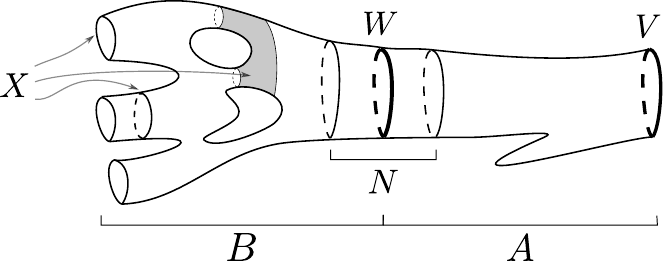}
\caption{}\label{fig:cylinder_split}
\end{figure}
\end{enumerate}
\end{theorem}

\begin{theorem}[Simplest functions on $2$-disk]\label{th:stab:disk:one_crpt}
Suppose $\func\in\FF(\Disk,\Pman)$ has a \myemph{unique critical point} $z$ being therefore a \myemph{local extreme}.\removed{, and let $\Vman=\partial\Disk$.}
\begin{enumerate}[label={\rm(\arabic*)}]
\item\label{enum:2disk:1crpt:nondeg}
If $z$ is \myemph{non-degenerate}, then
\added{$\sDSG{\func,\partial\Disk} \cong \seqTriv: \{1\} \monoArrow \{1\} \epiArrow \{1\}$.}
\item\label{enum:2disk:1crpt:deg}
If $z$ is \myemph{degenerate} of symmetry index $m$, then \added{$\sDSG{\func,\partial\Disk} \cong \seqZ{m}: m\bZ\monoArrow \bZ\epiArrow\bZ_m$.}
\end{enumerate}
\end{theorem}

Theorem~\ref{th:stab:annulus}\ref{eqnu:th:cyl:func:no_cr_pt} and Theorem~\ref{th:stab:disk:one_crpt} are reformulations of results established in different papers by the author, \cite{Maksymenko:AGAG:2006, Maksymenko:NO:ENG:2009, Maksymenko:NO:ENG:2010}.
We will present in~\S\ref{sect:proof:th:cyl:func}-\ref{sect:proof:th:disk:func_one_crpt} a unified approach to their proofs.

{\bf Functions on $2$-disk and annulus. General case.}
Again we assume that $(\Mman,\Vman)$ is one of the pairs $(\Disk,\partial\Disk)$ or $(\Circle\times\UInt, \Circle\times 0)$.
Then in both cases we have isomorphisms:
\[
\seqStab{\func,\Vman} 
\stackrel{\eqref{equ:StabIsotIdD2}}{\cong}
\seqStabIsotId{\func,\Vman}
\stackrel{\eqref{equ:func_on_cylinder_hom_eq1}}{\cong}
\seqStabIsotId{\func,\partial\Mman}.
\]
Our aim is to compute this sequence.

\begin{figure}[!htbp]
\centering
\begin{tabular}{ccc}
\includegraphics[height = 2.2cm]{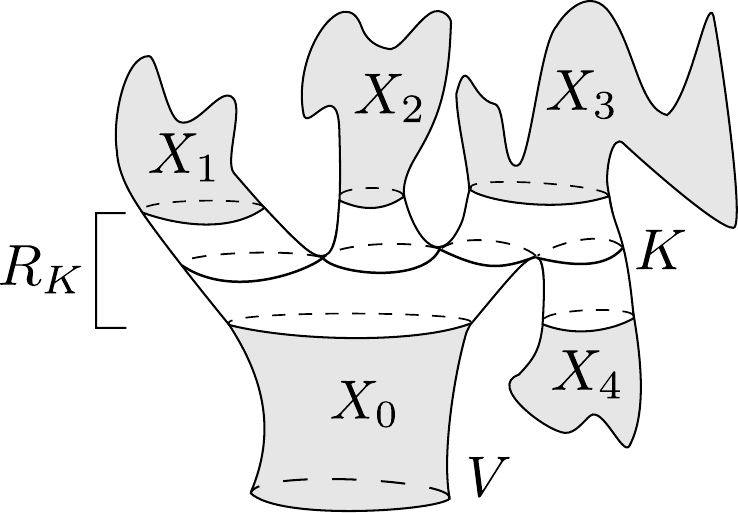} & \qquad\qquad &
\includegraphics[height = 2.2cm]{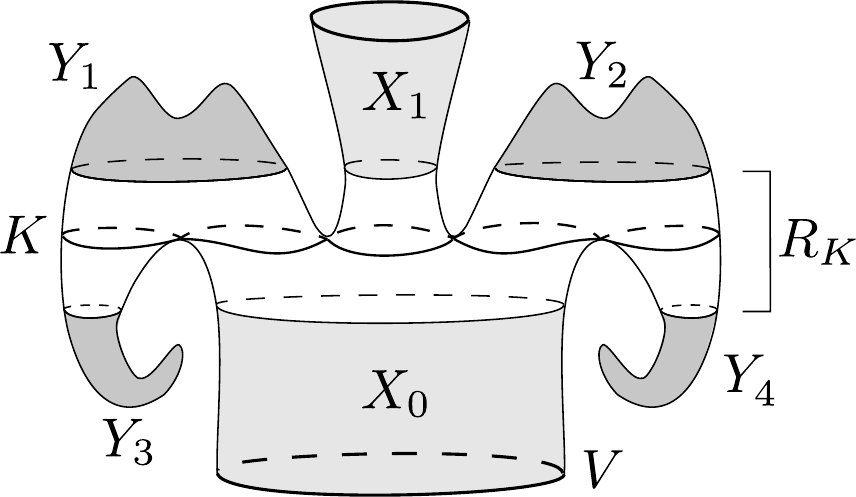} \\
(a) & & (b)
\end{tabular}
\caption{}\label{fig:critical_comp}
\end{figure}

Let $\crLev$ be the ``closest'' to $\Vman$ non-extremal critical component of $\func$.
More precisely, $\crLev$ is a unique critical component of a level set of $\func$ such that the connected component of $\Mman\setminus\crLev$ containing $\Vman$ includes no critical points of $\func$, see~Figure~\ref{fig:critical_comp}.
Then 
\begin{equation}\label{equ:hK_K}
\dif(\crLev) =\crLev, \qquad \forall\ \dif\in\Stabilizer{\func,\Vman}.
\end{equation}

Let $c=\func(\crLev) \in \Pman$ and $\eps>0$.
Denote by $\regNK$ the connected component of $\func^{-1}[c-\eps,c+\eps]$ containing $\crLev$.
Decreasing $\eps$ we can assume that $\regNK \setminus \crLev$ contains no critical points of $\func$ and $\regNK \cap \partial\Mman \subset\partial\regNK$ (though this intersection can be empty).
Then $\regNK$ is an $\func$-regular neighborhood of $\crLev$.
\removed{being invariant with respect to $\Stabilizer{\func,\Vman}$.}

\added{
Denote by $\bZman$ the family of all connected components of $\overline{\Mman\setminus\regNK}$.
It follows from~\eqref{equ:hK_K} that $\dif(\regNK)=\regNK$ for all $\dif\in\Stabilizer{\func,\Vman}$, and therefore $\dif$ also interchanges elements of $\bZman$.
Thus we get a natural action of $\Stabilizer{\func,\Vman}$ on $\bZman$.
Let 
\begin{equation}\label{equ:StabZ}
 \Stabilizer{\bZman} = \{ \dif\in\Stabilizer{\func,\Vman} \mid \dif(\Zman)=\Zman \ \text{for each} \ \Zman\in\bZman\}
\end{equation}
be its kernel of non-effectiveness.
Then the quotient $\Stabilizer{\func,\Vman}/\Stabilizer{\bZman}$ \myemph{effectively} acts on $\bZman$.
}

\begin{itemize}[label=$\bullet$, wide]
\item
Let $\bZmanX=\{ \XFixA, \Xman_1,\ldots,\Xman_a\}$ be all the elements of $\bZman$ invariant under $\Stabilizer{\func,\Vman}$, i.e. fixed points of $\Stabilizer{\func,\Vman}/\Stabilizer{\bZman}$, see Figure~\ref{fig:critical_comp}a).
We will always enumerate them so that $\Vman\subset \XFixA$, and in particular, $\XFixA$ is always an annulus.
Notice that if $\Mman=\Disk$, them all $\Xman_i$, $i\geq1$, are $2$-disks.
On the other hand, if $\Mman=\Circle\times\UInt$, then the element of $\bZman$ containing another boundary component $\Circle\times1$ is  invariant under $\Stabilizer{\func,\Vman}$ and we will always denote it by $\XFixB$.
In this case $\XFixB$ is an annulus, and all others $\Xman_i$, $i\geq2$, must be $2$-disks.

\item
Let also $\bZmanY:=\bZman\setminus\bZmanX = \{ \Yman_1,\ldots,\Yman_b \}$ be all other $2$-disks of $\bZman$ being not invariant under some elements of $\Stabilizer{\func,\Vman}$, see Figure~\ref{fig:critical_comp}b).
\end{itemize}

\begin{remark}\label{rem:get_out_inv_comp}\rm
Suppose $\Xman_i\in\bZmanX$, $i\geq1$, is an $\Stabilizer{\func,\Vman}$-invariant element of $\bZman$ such that the surface $\Cyli{i}=\overline{\Mman\setminus\Xman_i}$ is an \myemph{annulus}.
Then $\regNK \cap \partial\Mman = \varnothing$ and we have an isomorphism
\[\seqStabIsotId{\func,\partial\Mman}
\stackrel{\eqref{equ:bseq_two_cyl}}{\cong}
\seqStabIsotId{\func|_{\Xman_i},\partial\Xman_i} \times
\seqStabIsotId{\func|_{\Cyli{i}},\partial\Cyli{i}}.\]
Consider the restriction $\func|_{\Cyli{i}} \in \FF(\Cyli{i},\Pman)$.
Then $\regNK$ is still an $\func|_{\Cyli{i}}$-regular neighborhood of $\crLev$ in $\Cyli{i}$, however now $\regNK \cap \partial\Cyli{i} = \Xman_i\cap\Cyli{i} \not=\varnothing$.
Moreover, let $\bZman_{\Cyli{i}}:=\bZmanX_{\Cyli{i}}\sqcup\bZmanY_{\Cyli{i}}$
be the collection of connected components of $\overline{\Cyli{i}\setminus\regNK}$, where $\bZmanX_{\Cyli{i}}$ is the set of $\Stabilizer{\func|_{\Cyli{i}},\Vman}$-invariant elements of $\bZman_{\Cyli{i}}$ and $\bZmanY_{\Cyli{i}} := \bZman_{\Cyli{i}}\setminus\bZmanX_{\Cyli{i}}$.
Then
\begin{align}\label{equ:Z_A}
\bZmanX_{\Cyli{i}} &= \bZmanX \setminus\{\Xman_i\}, &
\bZmanY_{\Cyli{i}} &= \bZmanY, &
\bZman_{\Cyli{i}}  &= \bZman \setminus\{\Xman_i\}.
\end{align}
\end{remark}

\begin{theorem}[General case]\label{th:stab:disk_ann:gen_case}
\begin{enumerate}[wide, label={\rm(\Alph*)}, itemsep=1ex]
\item\label{enum:SS:A}
Suppose $\bZmanY=\varnothing$, i.e. $\bZman=\bZmanX=\{\XFixA,\Xman_1,\ldots,\Xman_a\}$ so all its elements are invariant under $\Stabilizer{\func,\Vman}$, which also means that $\Stabilizer{\func,\Vman}=\Stabilizer{\bZman}$.
Then 
\begin{equation}\label{equ:DOS__jSdS_case_A}
\seqStabIsotId{\func,\partial\Mman} \cong  
\Bigl( \,
 \mprod\limits_{i=1}^{a}\seqStabIsotId{\func|_{\Xman_i},\partial\Xman_i}
\, \Bigr)
 \ \times \ \seqZ{1}
\ \stackrel{\eqref{equ:W_Z1}}{\equiv} \
\seqWrm{ 
\Bigl( \,
 \mprod\limits_{i=1}^{a}\seqStabIsotId{\func|_{\Xman_i},\partial\Xman_i}
\, \Bigr)}{1}.
\end{equation}

\item\label{enum:SS:B}
Suppose $\bZmanY=\{\Yman_i\}_{i=1}^{b}\not=\varnothing$.
Then $\Stabilizer{\func,\Vman}/\Stabilizer{\bZman} \cong \bZ_m$ for some $m\geq2$ and its action on $\bZman$ is \myemph{semifree}, i.e.\! free on the set $\bZmanY$ of non-fixed elements.
Hence $m$ divides $b$ and that action has exactly $c:=b/m$ orbits.
Moreover, it also has exactly either \myemph{one} or \myemph{two} fixed elements, and in each of the corresponding cases $\sDSG{\func,\Vman}$ can be described as follows.

Fix any ($2$-disks) $\Yman_1,\ldots,\Yman_c\in\bZmanY$ belonging to mutually distinct non-fixed $\bZ_m$-orbits and define the following short exact sequence:
\begin{equation}\label{equ:DOS__jSdS_case_B}
\aligned
\aSeq &:= \seqWrm{ 
   \Bigl( \myprod\limits_{j=1}^{c}
          \seqStabIsotId{\func|_{\Yman_i}, \partial\Yman_i} \Bigr) 
}{m}
\endaligned
\end{equation}

\begin{enumerate}[wide, label={\rm(\alph*)}, itemsep=1ex]
\item\label{enum:SS:B:B0}
If $\bZmanX=\{\XFixA\}$, then $\seqStabIsotId{\func,\partial\Mman} \cong \aSeq$.

\item\label{enum:SS:B:B0_X1}
Otherwise, $\bZmanX=\{\XFixA, \XFixB\}$.
Let $\Cyli{1}=\overline{\Mman\setminus\XFixB}$.
Then 
$\bZmanX_{\Cyli{1}} \stackrel{\eqref{equ:Z_A}}{=} \bZmanX \setminus\{\XFixB\} = \{\XFixA\}$, whence the restriction $\func|_{\Cyli{1}}$ satisfies condition~\ref{enum:SS:B}\ref{enum:SS:B:B0}.
Therefore $\seqStabIsotId{\func|_{\Cyli{1}},\partial\Cyli{1}} \cong \aSeq$ and 
$\seqStabIsotId{\func,\partial\Mman} \stackrel{\eqref{equ:bseq_two_cyl}}{\cong}\seqStabIsotId{\func|_{\Xman_1},\partial\Xman_1} \times \aSeq$.
\end{enumerate}
\end{enumerate}
\end{theorem}

\begin{corollary}[Decomposition into annuli and possibly one disk]\label{cor:stab:disk_ann:ann_decomp}
Let $\Mman$ be either a $2$-disk or an annulus and $\func\in\FF(\Mman,\Pman)$.
Then there are $\func$-adapted subsurfaces $\Cyli{1},\ldots,\Cyli{n} \subset \Mman$ having the following properties.
\begin{enumerate}[label={\rm(\alph*)}, wide]
\item\label{enum:cor:split_into_cyls:1} 
For $i=1,\ldots,n-1$ the surface $\Cyli{i}$ is an \myemph{annulus}, while $\Cyli{n}$ is either a $2$-disk or an annulus.
Moreover, the intersection $\Cyli{i}\cap \Cyli{j}$ for $i<j$ is non-empty only for $j=i+1$ and in this case it is a common boundary component of these subsurfaces.
Also $\Cyli{1}$ contains some boundary component of $\Mman$.

\item\label{enum:cor:split_into_cyls:2}
For each $i=1,\ldots,n-1$ there exist $m_i,c_i\geq1$ and certain $\func$-adapted mutually disjoint \myemph{$2$-disks} $\Yman_{i,1},\ldots,\Yman_{i,c_i} \subset \Cyli{i}$, such that we have an isomorphism
\begin{equation}\label{equ:seq:SDG:Cyl_i}
\seqStabIsotId{\func|_{\Cyli{i}}, \partial\Cyli{i}} \cong 
\seqWrm{\biggl(\myprod\limits_{j=1}^{c_i} \seqStabIsotId{\func|_{\Yman_{i,j}}, \partial\Yman_{i,j}}\biggr)}{m_i},
\end{equation}
while $\seqStabIsotId{\func|_{\Cyli{n}}, \partial\Cyli{n}}$ is isomorphic either with a sequence of type~\eqref{equ:seq:SDG:Cyl_i} or with $\seqTriv$ or with $\seqZ{m_n}$ for some $m_n\geq1$.

\item\label{enum:cor:split_into_cyls:3}
$\seqStabIsotId{\func, \partial\Mman}\cong \prod\limits_{i=1}^{n}\seqStabIsotId{\func|_{\Cyli{i}}, \partial\Cyli{i}} $.

\item\label{enum:cor:split_into_cyls:4} 
$\seqStabIsotId{\func, \partial\Mman} \in \ZBP \subset \classZBP$, see Lemma~\ref{lm:classes_SG}\ref{enum:lm:classes_SG:ZBP}.
\end{enumerate} 
\end{corollary}
\begin{proof}
1) First suppose $\func$ satisfies assumptions either of Theorem~\ref{th:stab:annulus}\ref{eqnu:th:cyl:func:no_cr_pt} or Theorem~\ref{th:stab:disk:one_crpt} or Theorem~\ref{th:stab:disk_ann:gen_case}\ref{enum:SS:B}\ref{enum:SS:B:B0}.
Then $\seqStabIsotId{\func, \partial\Mman}$ is either $\seqTriv$ or $\seqZ{m_1}$ for some $m_1\geq1$ or has the form~\eqref{equ:DOS__jSdS_case_B}.
In all these cases condition~\ref{enum:cor:split_into_cyls:2} holds and we can put $\Cyli{1}=\Mman$.
Thus Corollary~\ref{cor:stab:disk_ann:ann_decomp} holds with $n=1$.

\jadded{19.2-20.3}{}
2) Assume now that $\func$ satisfies either of the cases \ref{enum:SS:A} or \ref{enum:SS:B}\ref{enum:SS:B:B0_X1} of Theorem~\ref{th:stab:disk_ann:gen_case}.
Then we get a direct product splitting 
$\seqStabIsotId{\func, \partial\Mman} \stackrel{\eqref{equ:bseq_two_cyl}}{\cong}
\seqStabIsotId{\func|_{\Qman_1}, \partial\Qman_1}
\times 
\seqStabIsotId{\func|_{\Cyli{1}}, \partial\Cyli{1}}$, where $\Cyli{1}$ is an annulus containing $\Vman$, $\Qman_1 = \overline{\Mman\setminus\Cyli{1}}$ is either a $2$-disk or an annulus, and $\Qman_1\cap\Cyli{1}$ is a common boundary component of these subsurfaces.

Then, due to Remark~\ref{rem:get_out_inv_comp}, the restriction $\func|_{\Cyli{1}}$ satisfies 
either of the cases \ref{enum:SS:A} or \ref{enum:SS:B}\ref{enum:SS:B:B0} of Theorem~\ref{th:stab:disk_ann:gen_case}.
In particular, $\seqStabIsotId{\func|_{\Cyli{1}}, \partial\Cyli{1}}$ has the form~\eqref{equ:seq:SDG:Cyl_i}.
Moreover, the restriction $\func|_{\Qman_1}$ again satisfies exactly one of the cases 1) or 2).
In the case 1) we just stop and put $\Cyli{2}:=\Qman_1$, which will complete the proof.
On the other hand, in the case 2) we again get a direct product splitting 
$\seqStabIsotId{\func|_{\Qman_1}, \partial\Qman_1}\cong\seqStabIsotId{\func|_{\Qman_2}, \partial\Qman_2} \times \seqStabIsotId{\func|_{\Cyli{2}}, \partial\Cyli{2}}$, which yields an isomorphism $\seqStabIsotId{\func, \partial\Mman} \cong 
\seqStabIsotId{\func|_{\Qman_2}, \partial\Qman_2}
\times 
\seqStabIsotId{\func|_{\Cyli{2}}, \partial\Cyli{2}}
\times 
\seqStabIsotId{\func|_{\Cyli{1}}, \partial\Cyli{1}}$.

Applying the same arguments to $\func|_{\Qman_2}$ and so on we will definitely stop at case 1) in a finite number of steps, and thus get the desired direct product splitting~\ref{enum:cor:split_into_cyls:3}, since at each step the number of critical points of $\func$ in the remained subsurface $\Qman_i$ decreases.

Thus $\seqStabIsotId{\func, \partial\Mman}$ is obtained from $\seqTriv$ by finite number of operations of direct product and ``wreath product'' with $\seqZ{m}$.
Hence by Lemma~\ref{lm:classes_SG}\ref{enum:lm:classes_SG:ZBP}, $\dDSG{\func}{\Mman}\in\ZBP$.
\end{proof}

\xadded{2}{0}{19.2-20.3}{}
\begin{theorem}\label{th:solvable_groups}{\rm(cf. Theorem~\ref{th:Kudryavtseva} and Corollary~\ref{cor:Of_hom_type}).}
Let $\Mman$ be a connected orientable compact surface distinct from $S^2$ and $T^2$, $\func\in\FF(\Mman,\Pman)$, and $\Xman \subset\Mman$ an $\func$-adapted submanifold such that $\Xman\not=\varnothing$ whenever $\Mman=\Disk$ or $\Circle\times\UInt$, so $\chi(\Mman)<\ptnum{\Xman}$, and thus $\DiffId(\Mman,\Xman)$ is always contractible.
Then the following statements hold.
\begin{enumerate}[wide,label={\rm(\alph*)}]
\item\label{enum:th:solvable_groups:a}
The sequence $\sDSG{\func, \Xman}:\pi_0\FolStabilizerIsotId{\func,\Xman}\monoArrow\pi_0\StabilizerIsotId{\func,\Xman}\epiArrow\GrpKRIsotId{\func,\Xman}$ belongs to the class $\ZBP \subset \classZBP$.
In particular, $\pi_0\FolStabilizerIsotId{\func,\Xman}$ is a free abelian group of some rank $k$, $\pi_0\StabilizerIsotId{\func,\Xman}$ is a solvable Bieberbach group, and $G := \GrpKRIsotId{\func,\Xman}$ is solvable and finite.
Hence, by Theorem~\ref{th:bieberbach}, there is a free action of $G$ on $k$-torus $T^{k}$ such that $\pi_1 (T^{k}/G) \cong \pi_0\StabilizerIsotId{\func,\Xman}$.

\item\label{enum:th:solvable_groups:b}
$\OrbitPathComp{\func,\Xman}{\func}$ is homotopy equivalent to the quotient $T^{k}/G$.
\end{enumerate}
\end{theorem}
\begin{proof}
\ref{enum:th:solvable_groups:a}
If $\Mman=\Disk$ or $\Circle\times\UInt$ and $\Xman\not=\varnothing$, then $\sDSG{\func, \Xman}\in\ZBP$ by Corollary~\ref{cor:stab:disk_ann:ann_decomp}\ref{enum:cor:split_into_cyls:4}.
In all other cases $\chi(\Mman)<0$, whence by Theorem~\ref{lm:reduction_Mconn_V_dM}, $\sDSG{\func, \Xman}$ splits into a direct product of finitely many $\DSG$-sequences of the form $\seqStabIsotId{\func|_{\Zman},\partial\Zman}$, where $\Zman$ is an $\func$-adapted $2$-disk or an annulus.
As just noted every such sequence belongs to $\ZBP$.
Therefore by Lemma~\ref{lm:classes_SG}\ref{enum:lm:classes_SG:ZBP} so does $\sDSG{\func, \Xman}$.

\ref{enum:th:solvable_groups:b}
Since by~\ref{enum:th:stab_orb_full_info:pikOf} of Theorem~\ref{th:stab_orb:full_info}, $\OrbitPathComp{\func,\Xman}{\func}$ is aspherical, we get from Theorem~\ref{th:bieberbach} that $\OrbitPathComp{\func,\Xman}{\func}$ and $T^{k}/G$ are homotopy equivalent.
\end{proof}

\begin{remark}\rm
The statement that $\pi_0\FolStabilizerIsotId{\func}$ is free abelian is proved in~\cite[Theorem~6.2]{Maksymenko:AGAG:2006} for Morse maps and in \cite[Theorem~5]{Maksymenko:ProcIM:ENG:2010} for all $\func\in\FF(\Mman,\Pman)$.
Moreover, that group appears as $\bZ^l$ in~\eqref{equ:exact_seq_for_pi1OfX:2}.
\end{remark}
The above statements allow easily to compute the groups $\pi_0\FolStabilizerIsotId{\func, \Xman}$, $\pi_0\StabilizerIsotId{\func, \Xman}$, $\GrpKRIsotId{\func, \Xman}$ for arbitrary $\func\in\FF(\Mman,\Pman)$ and an $\func$-adapted submanifold $\Xman$ whenever we know ``discrete symmetries of $\func$'' and $\Mman$ is connected, orientable and distinct from $S^2$ and $\Torus$.
To simplify notations we will use the symbol ``='' for isomorphisms between of $\funcSeq$- and $\funcSeq'$-sequences. 

\begin{example}\label{exmp:computaiton}\rm
Let $\func\in\FF(\Disk,\bR)$ be a function having precisely three critical points: $\crsdl$, $\cra$, $\crb$, where $\crsdl$ is a non-degenerate saddle, $\cra$ and $\crb$ are local maximums.
We will assume that $\func(\partial\Disk) = 0$, $\func(\crsdl)=1$, and $\func(\cra), \func(\crb) > 1.5$.
Then $\crLev = \func^{-1}(1)$ is a \myemph{connected} critical level set containing $\crsdl$.
Let also $\regNK = \func^{-1}[0.5, 1.5]$, $\XFixA=\func^{-1}[0,0.5]$ be the connected component of $\overline{\Disk\setminus\regNK}$ containing $\partial\Disk$, and $\Xman_{1}$, $\Xman_{2}$ be the connected components of $\overline{\Disk\setminus\regNK}$ containing $\cra$ and $\crb$ respectively.
Consider three cases.

1) Suppose $\func(\cra) > \func(\crb)$.
Since $\func$ takes distinct values at $\cra$ and $\crb$, every $\dif\in\Stabilizer{\func,\partial\Disk}$ leaves invariant each $\Xman_i$ and we have the case~\ref{enum:SS:A} of Theorem~\ref{th:stab:disk_ann:gen_case}.
Therefore
\[ \dDSG{\func}{\Disk}  = \dDSG{\func}{\Xman_1} \times \dDSG{\func}{\Xman_2} \times \seqZ{1}. \]

a) Assume that $\cra$ and $\crb$ are non-degenerate local extremes, see Figure~\ref{fig:example_1_2}-1a). 
Then the restrictions $\func|_{\Xman_1}$ and $\func|_{\Xman_2}$ satisfy statement~\ref{enum:2disk:1crpt:nondeg} of Theorem~\ref{th:stab:disk:one_crpt}, whence their $\DSG$-sequences consist of trivial groups, i.e. $\dDSG{\func}{\Xman_i} = \seqTriv$, $i=1,2$.
Therefore
\begin{align*}
\dDSG{\func}{\Disk} &  =  \seqTriv \times \seqTriv \times \seqZ{1}  =  \seqZ{1}: \bZ \monoArrow \bZ \epiArrow 1.
\end{align*}

\begin{figure}[ht]
\centering
\begin{tabular}{ccc}
\includegraphics[height=1.9cm]{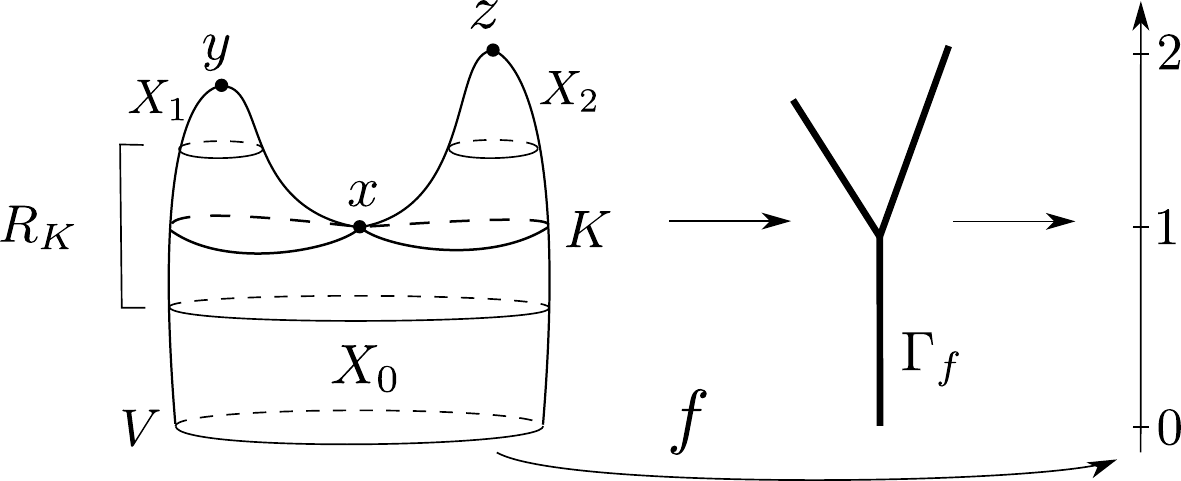} & \qquad &
\includegraphics[height=1.9cm]{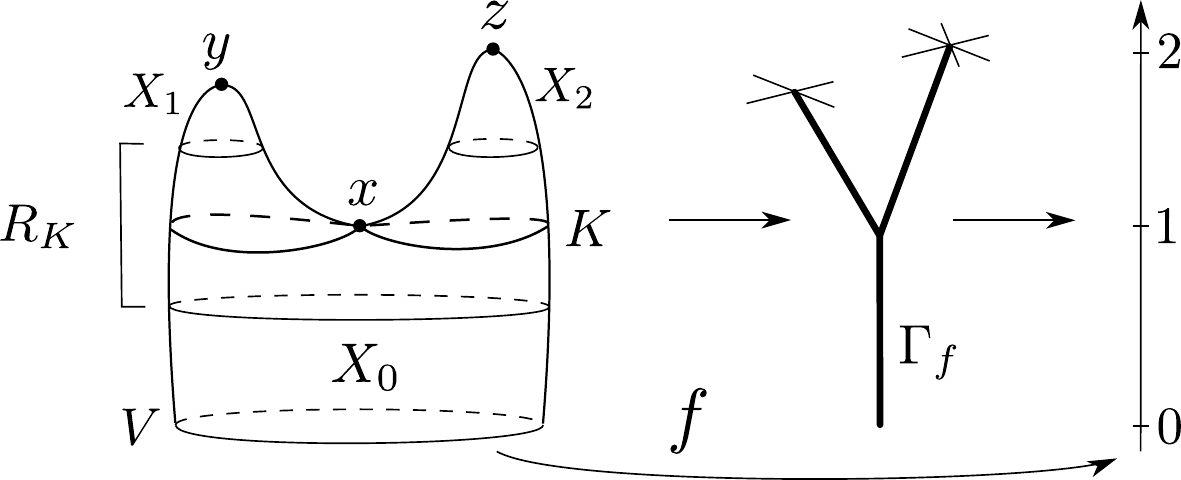} \\
1a) & & 1b) \\[1ex]
\includegraphics[height=1.9cm]{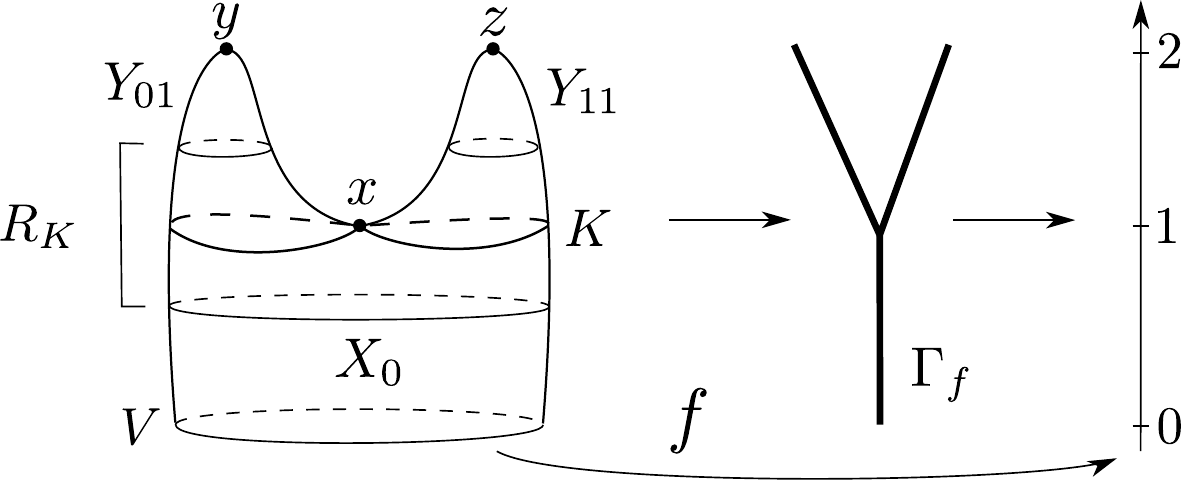}  & \qquad &
\includegraphics[height=1.9cm]{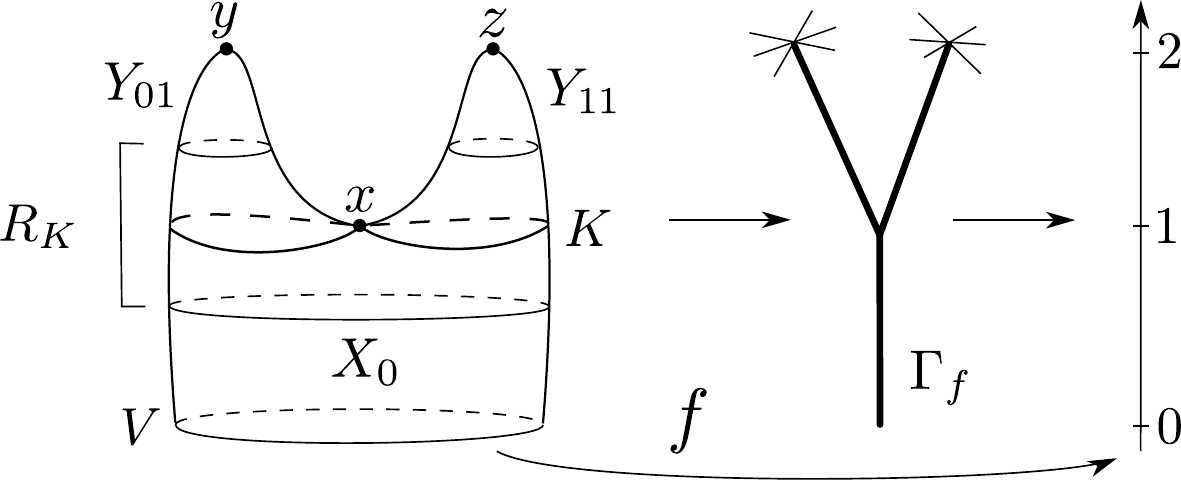} \\
2a) & & 2b)
\end{tabular}
\caption{}\label{fig:example_1_2}
\end{figure}

b) Suppose both $\cra$ and $\crb$ are \myemph{degenerate} and have symmetry indices $4$ and $6$ respectively, see Figure~\ref{fig:example_1_2}-1b).
Then by Theorem~\ref{th:stab:disk:one_crpt}\ref{enum:2disk:1crpt:deg}, $\dDSG{\func}{\Xman_1} = \seqZ{4}$, $\dDSG{\func}{\Xman_2} = \seqZ{6}$, whence
\[
\dDSG{\func}{\Disk}  = \seqZ{4} \times \seqZ{6} \times \seqZ{1}: \
\bZ^4 \times \bZ^6 \times \bZ 
 \monoArrow  (\bZ\wrm{4}\bZ)\times(\bZ\wrm{6}\bZ)\times\bZ \epiArrow \bZ_4 \times \bZ_6.
\]

2) Suppose $\func(\cra)=\func(\crb)=2$, and the germs of $\func$ at these points are smoothly equivalent.
Then there exists $\dif\in\Stabilizer{\func,\partial\Disk}$ such that $\dif(\cra)=\crb$ and $\dif(\crb) =\cra$.
Redenote $\Yman_{01}:=\Xman_1$ and $\Yman_{11}:=\Xman_2$.
Then $\dif$ must also interchange $\Yman_{01}$ and $\Yman_{11}$.
Therefore the action of $\Stabilizer{\func,\partial\Disk}$ on the set $\{\Yman_{01}, \Yman_{11}\}$ of non-invariant connected components of $\overline{\Disk\setminus\regNK}$ reduces to a free action of $\bZ_2$.
Thus we have a situation described in Theorem~\ref{th:stab:disk_ann:gen_case}\ref{enum:SS:B}\ref{enum:SS:B:B0} with $m=2$ and $c=1$ orbit.
Hence 
\begin{equation}\label{equ:DSG_pants}
\dDSG{\func}{\Disk} = \seqWrm{\dDSG{\func}{\Yman_{01}}}{2}.
\end{equation}

a) If $\cra$ and $\crb$ are non-degenerate, see Figure~\ref{fig:example_1_2}-2a), then $\dDSG{\func}{\Yman_{01}} = \seqTriv$, whence 
\[\dDSG{\func}{\Disk} = \seqWrm{\seqTriv}{2} = \seqZ{2}: \ 2\bZ \monoArrow \bZ\epiArrow \bZ_2.\]

b) Suppose $\cra$ (and therefore $z$) is degenerate of symmetry index $k$, see Figure~\ref{fig:example_1_2}-2b).
Then $\dDSG{\func}{\Yman_{01}}=\seqZ{k}$, whence, see~\eqref{equ:Zm_Zn}, 
\begin{equation}\label{equ:DSG_ZkZ2}
\dDSG{\func}{\Disk} = \seqWrm{\seqZ{k}}{2}: \
(k\bZ)^2 \times 2\bZ \monoArrow 
\bZ\wrm{2}\bZ \epiArrow
\bZ_k\wr\bZ_2.
\end{equation}

3) To illustrate once again such a recursive character of calculations consider another function $\func\in\FF(\Disk,\bR)$ having three critical values $1$, $2$ and $3$ such that $\crLev = \func^{-1}(1)$ contains a unique non-degenerate saddle, $\func^{-1}(2)$ contains two non-degenerate saddles, $\func^{-1}(3)$ consists of $4$ degenerate global extremes of the same symmetry index $k$ such that the germs of $\func$ at these extremes are mutually smoothly right equivalent, see Figure~\ref{fig:example_3}.
Let $\regNK = \func^{-1}[0.5, 1.5]$, and $X_0$, $\Yman_{01}$ and $\Yman_{11}$ be the connected components of $\overline{\Disk\setminus\regNK}$.
Then, similarly to the case 2a), the action of $\Stabilizer{\func,\Vman}$ on the set $\{\Yman_{01}, \Yman_{11}\}$ of non-invariant connected components of $\overline{\Disk\setminus\regNK}$ reduces to a free action of $\bZ_2$ with one $c=1$ orbit.
Moreover, by the construction the restriction of $\func$ to $\Yman_{01}$ is ``the same'' as the function from the case 2b).
Hence 
\begin{align*}
\dDSG{\func}{\Disk} & 
   \ \stackrel{\eqref{equ:DSG_pants}}{=} \
\seqWrm{\dDSG{\func}{\Yman_{01}}}{2} 
   \ \stackrel{\eqref{equ:DSG_ZkZ2}}{=}
\seqWrm{ \bigl( \seqWrm{\seqZ{k}}{2} \bigr) }{2}    \\
& \ \stackrel{\phantom{\eqref{equ:DSG_pants}}}{=} \
    \bigl((k\bZ)^2 \times 2\bZ\bigr)^2 \times 2\bZ \ \monoArrow \
    (\bZ \wrm{2} \bZ)  \wrm{2} \bZ \ \epiArrow \
    (\bZ_k \wr \bZ_2)  \wr \bZ_2.
\end{align*}

\begin{figure}[ht]
\centering
\includegraphics[height=2.4cm]{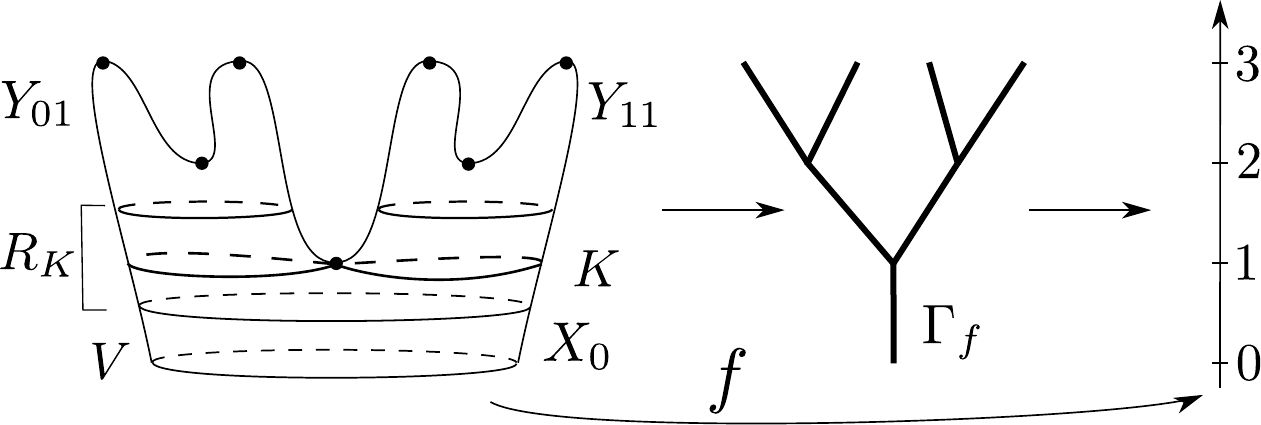}
\caption{Case 3)}\label{fig:example_3}
\end{figure}

4) Let $\Mman=\Circle\times\UInt$ be an annulus and $\func:\Mman\to\bR$ a Morse function having exactly two critical values $1$ and $3$.
Suppose also that $\func^{-1}(0)=\Circle\times 0$, $\crLev:=\func^{-1}(1)$ is connected and contains $k$ saddles, and $\func^{-1}(3)$ consists of $m$ non-degenerate local maximums, see Figure~\ref{fig:example_4}-4a) for $m=4$.
Let $\regNK$ be an $\func$-regular neighborhood of $\crLev$.
Then $\overline{\Mman\setminus\regNK}$ consists of $m+2$ connected components: $\XFixA$, $\XFixB$, $\Yman_{0i}$, $i=0,\ldots,m-1$, such that $\Circle\times s \subset \Xman_s$ for $s=0,1$, and each $\Yman_{0i}$ is a $2$-disk containing a unique local maximum of $\func$.
Evidently, each restriction $\func|_{\Yman_{0i}}$ satisfies assumption~\ref{enum:2disk:1crpt:nondeg} of Theorem~\ref{th:stab:disk:one_crpt}, whence $\dDSG{\func}{\Yman_{0i}}=\seqTriv$ consists of trivial groups.
Moreover, every $\dif\in\Stabilizer{\func,\partial\Mman}$ cyclically interchanges $\{\Yman_{0i}\}_{i=0,\ldots,m-1}$, whence we have case~\ref{enum:SS:B} of Theorem~\ref{th:stab:disk_ann:gen_case} for one orbit (i.e. $c=1$) of length $m$, and so
\[
\dDSG{\func}{\Mman} =
\seqWrm{ \dDSG{\func}{\Yman_{00}} }{m} =
\seqWrm{\seqTriv}{m} = \seqZ{m}.
\]
\begin{figure}[ht]
\centering
\begin{tabular}{ccc}
\includegraphics[height=3cm]{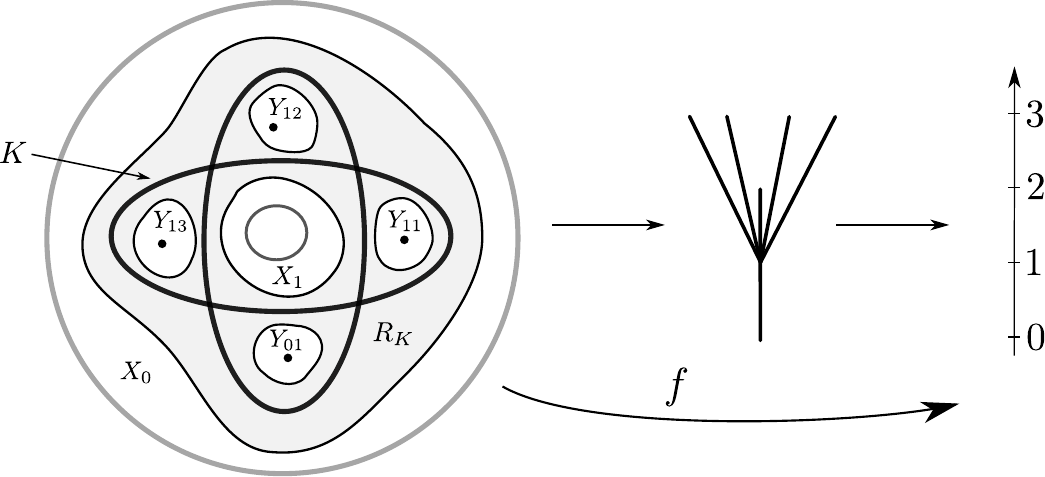} &  \qquad
\includegraphics[height=3cm]{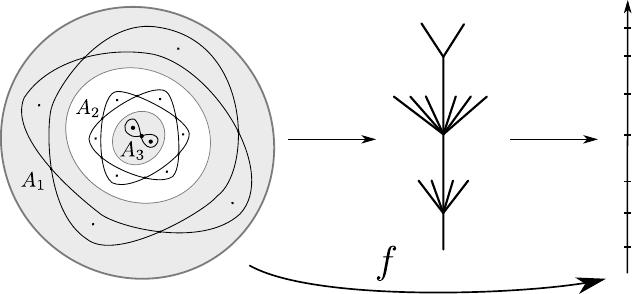} \\
4a) & 4b)
\end{tabular}
\caption{}\label{fig:example_4}
\end{figure}

5) Finally, let $\func:\Disk\to\bR$ be a Morse function whose critical components of level sets are shown in Figure~\ref{fig:example_4}-4b).
\jadded{19.1}{}%
Thus $\Disk$ is a union of two $\func$-adapted annuli $\Cyli{1}$, $\Cyli{2}$ and an $\func$-adapted disk $\Cyli{3}$.
Assume also that $\func$ takes the same value at all local extremes belonging to the same $\Cyli{i}$, $i=1,2,3$.
Then $\func|_{\Cyli{1}}$ and $\func|_{\Cyli{2}}$ are ``the same'' as in 3a) for $k=4$ and $6$ respectively, while $\func|_{\Cyli{3}}$ is ``the same'' as in 2a).
Hence by Corollary~\ref{cor:stab:disk_ann:ann_decomp}\ref{enum:cor:split_into_cyls:3},
\begin{align*}
\dDSG{\func}{\Mman} & =
\dDSG{\func}{\Cyli{1}} \times  \dDSG{\func}{\Cyli{2}} \times \dDSG{\func}{\Cyli{3}} 
= \seqZ{4} \times  \seqZ{6} \times  \seqZ{2}: \\
&  (\bZ^4 \times 4\bZ) \times  (\bZ^6 \times 6\bZ) \times 2\bZ \ \monoArrow \
    (\bZ\wrm{4}\bZ) \times (\bZ\wrm{6}\bZ) \times \bZ  \ \xepiArrow{} \
    \bZ_4 \times \bZ_6 \times \bZ_2.
\end{align*}
\end{example}

\section{Shifts along orbits of flows}\label{sect:shifts_along_orbits}
In order to make the paper as self-contained as possible, we will present short proofs of some results obtained in~\cite{Maksymenko:TA:2003}.
Let $\Mman$ be a smooth \added{possibly non-connected} manifold and $\fld$ a smooth vector field on $\Mman$ tangent to $\partial\Mman$ and generating a global flow $\flow:\Mman\times\bR\to\Mman$.
For each smooth function $\alpha:\Mman\to\bR$\jadded{20.4}{,} let $\fld(\alpha)$ be the Lie derivative of $\alpha$ with respect to $\fld$ and $\flow_{\alpha}:\Mman\to\Mman$ \jremoved{20.5}{be} the map defined by
\begin{equation}\label{equ:shift_via_alpha}
\flow_{\alpha}(x) = \flow(x,\alpha(x)), \qquad x\in\Mman.
\end{equation}
We will call $\flow_{\alpha}$ the \myemph{shift along orbits of $\flow$} via the function $\alpha$, and $\alpha$ will in turn be called a \myemph{shift function} for $\flow_{\alpha}$.
\begin{lemma}\label{lm:shift}{\rm\jadded{20.6}{(\cite{Maksymenko:TA:2003}).}}
The following statements hold.
\begin{enumerate}[label=$\mathrm{(\roman*)}$, wide]
\item\label{enum:shift:local_existence}
Let $\Uman$ be an open subset of $\Mman$, $\dif:\Uman\to\Mman$\jremoved{20.7}{be} a $C^{k}$ map, $k\geq0$, such that $\dif(\orb\cap\Uman)\subset\orb$ for each orbit $\orb$ of $\flow$, $z\in\Uman$\jremoved{20.8}{be} a non-singular point of $\fld$, i.e.\! $\fld(z)\not=0$, and $a\in\bR$ \jchanged{20.9}{be}{a value} such that $\dif(z)=\flow(z,a)$.
Then there exists an open neighborhood $\Vman$ of $z$ and a \myemph{unique} $C^{k}$ function $\alpha:\Vman\to\bR$ such that $\alpha(z)=a$ and $\dif(x) = \flow(x,\alpha(x))$ for all $x\in\Vman$.
Moreover, $\dif$ is a local diffeomorphism at $z$ if and only if $1+\fld(\alpha)(z)\not=0$.

\item\label{enum:shift:shift_func_uniqueness}
Suppose $\Uman\subset\Mman$ is an open connected subset containing no singular points of $\fld$, and $\alpha,\beta:\Uman\to\bR$ are continuous functions such that $\flow(x,\alpha(x))=\flow(x, \beta(x))$ for all $x\in\Uman$.
If $\alpha(z)=\beta(z)$ for some $z\in\Uman$, then $\alpha=\beta$ on all of $\Uman$.
In particular, if $\Uman$ contains a non-closed orbit $\omega$ of $\flow$, then $\alpha=\beta$ on $\omega$ and therefore on all of $\Uman$.

\item\label{enum:shift:composition}
\added{For any smooth functions $\alpha,\beta:\Mman\to\bR$ we have that $\flow_{\alpha}\circ\flow_{\beta} = \flow_{\beta + \alpha\circ\flow_{\beta}}$.
Moreover, if $\flow_{\beta}$ is a diffeomorphism, then $\flow_{\beta}^{-1} = \flow_{-\beta\circ\flow_{\beta}^{-1}}$, whence}
\begin{equation}\label{equ:shift_fa_fbinv}
\added{\flow_{\alpha}\circ\flow_{\beta}^{-1} = 
\flow_{\alpha}\circ\flow_{-\beta\circ\flow_{\beta}^{-1}} = 
\flow_{-\beta\circ\flow_{\beta}^{-1} + \alpha\circ\flow_{\beta}^{-1}} =
\flow_{(\alpha-\beta)\circ\flow_{\beta}^{-1}}.}
\end{equation}

\item\label{enum:shift:shift_func_existence}
Suppose all orbits of $\flow$ are non-singular and non-closed, and $\dif:\Mman\to\Mman$ is a $C^{k}$ map, $k\geq0$, such that $\dif(\orb)=\orb$ for each orbit of $\flow$.
Then there exists a unique $C^{k}$ function $\alpha:\Mman\to\bR$ such that $\dif = \flow_{\alpha}$.

\item\label{enum:shift:diffeomorphism}
Consider the following \removed{convex} set of smooth functions:
\begin{equation}\label{equ:shift_func_for_diff}
 \DFunc =\{ \alpha\in C^{\infty}(\Mman,\bR) \mid \fld(\alpha)(x) > -1 \ \text{for all}  \ x\in\Mman \},
\end{equation}
\jadded{21.1}{being \myemph{convex} as it is defined by a linear inequality.}
Suppose $\Mman$ is compact and let $\alpha\in C^{\infty}(\Mman,\bR)$.
Then the map $\flow_{\alpha}$ is a diffeomorphism of $\Mman$ preserving orientation of orbits of $\fld$ if and only if $\alpha\in\DFunc$.

\item\label{enum:shift:isotopy_between_shifts}
Let $\alpha\in\DFunc$ and $\mu:\Mman\to\UInt$ be a $C^{\infty}$ function such that $\fld(\mu)=0$, i.e. $\mu$ is constant along orbits of $\fld$.
Then $\mu\alpha\in\DFunc$ as well.
Thus if $\flow_{\alpha}$ is a diffeomorphism, then $\flow_{\mu\alpha}$ is also a diffeomorphism.

In particular, $t\alpha\in\DFunc$ for all $t\in\UInt$ and so the family of maps $\{\flow_{t\alpha}\}_{t\in\UInt}$, is an orbit preserving isotopy between $\id_{\Mman}=\flow_{0}$ and $\flow_{\alpha}$.
\added{Also if $\alpha,\beta\in\DFunc$ and $\alpha=\beta$ on some subset $\Xman\subset\Mman$, then $\flow_{t\alpha+(1-t)\beta}$ is an isotopy between $\flow_{\alpha}$ and $\flow_{\beta}$ relatively to $\Xman$.}
\end{enumerate}
\end{lemma}
\begin{proof}
\ref{enum:shift:local_existence}
First suppose that $a=0$, whence $\dif(z)=z$.
As $\fld(z)\not=0$, there are local coordinates $(x_1,\ldots,x_m)$ on some neighborhood $\Uman$ of $z$ in which $z=0$, $\fld = \tfrac{\partial}{\partial x_1}$, and $\flow(x_1,\ldots,x_m,t) = (x_1+t, x_2,\ldots,x_m)$ for all sufficiently small $(x_1,\ldots,x_m, t)$.
Let $\Vman\subset\Uman$ be a connected open neighborhood of $z$ such that $\dif(\Vman)\subset\Uman$.
Since $\dif$ preserves orbits of $\flow$ being straight lines parallel to $x_1$-axis, it follows that $\dif(x_1,\ldots,x_m)=(q(x_1,\ldots,x_m), x_2,\ldots,x_m)$, for some \myemph{unique} $C^k$ function $q:\Vman\to\bR$ such that $q(z)=0$.
Define $\alpha:\Vman\to\bR$ by 
\begin{equation}\label{equ:shift_func_formula:0}
\alpha(x) = q(x) - x_1.
\end{equation}
Then $\dif(x) = \flow(x,\alpha(x))$.

Moreover, the Jacobian of $\dif$ at $z$ equals to $\tfrac{\partial q}{\partial x_1}(z) = 1 + \tfrac{\partial \alpha}{\partial x_1}(z) = 1+\fld(\alpha)(z)$ and the latter expression does not depend on a particular choice of local coordinates at $z$.
Hence $\dif$ is a local diffeomorphism at $z$ iff $1+\fld(\alpha)(z)\not=0$.

Also notice that $\dif$ preserves orientation of orbits at $z$ iff $1+\fld(\alpha)(z)>0$.

\smallskip

Suppose $a\not=0$.
Then $\dif' = \flow_{-a}\circ\dif$ satisfies the assumptions of~\ref{enum:shift:local_existence} with $a=0$, and so $\dif'(x)= \flow(x,\alpha'(x))$ for a unique $C^{k}$ function $\alpha':\Vman\to\bR$.
Let 
\begin{equation}\label{equ:shift_func_formula:a}
\alpha = a+\alpha'.
\end{equation}
Then $\dif(x) = \flow_{a}\circ \flow(x,\alpha'(x))=\flow(x,a+\alpha(x))=\flow(x,\alpha(x))$ for all $x\in\Vman$.
Moreover, $\dif$ is a local diffeomorphism at $z$ iff $\dif'$ is so, that is when $1+\fld(\alpha')(z)\not=0$.
But $\fld(\alpha)=\fld(a+\alpha')=\fld(\alpha')$, whence $\dif$ is a local diffeomorphism at $z$ iff $1+\fld(\alpha)(z)\not=0$.

\ref{enum:shift:shift_func_uniqueness}
Let $A = \{x\in\Uman \mid \alpha(x) = \beta(x)\}$.
Then $A$ is a non-empty closed subset of $\Uman$.
Moreover, 
\xchanged{0}{0}{21.2}{$A$ is open due to (i)}
{if $z\in A$, that is $\alpha(z)=\beta(z)$, then due to formulas~\eqref{equ:shift_func_formula:0} and~\eqref{equ:shift_func_formula:a}, $\alpha=\beta$ on some neighborhood $\Vman$ of $z$.
This means that $\Vman\subset A$, and thus $A$ is open.}
Since $\Uman$ is connected, \added{we obtain that} $A=\Uman$.

Statement \ref{enum:shift:composition} can be proved by a straightforward verification, see~\cite[Proposition~3]{Maksymenko:TA:2003}.

\ref{enum:shift:shift_func_existence}
Let $x\in\Mman$.
Since its orbit is non-closed, there exists a unique number $\alpha(x)\in\bR$ such that $\dif(x) = \flow(x,\alpha(x))$.
Moreover, it follows from~\ref{enum:shift:local_existence} that the obtained function $\alpha:\Mman\to\bR$ is of the same differentiability class as $\dif$.

\ref{enum:shift:diffeomorphism}
If $\flow_{\alpha}$ is a diffeomorphism preserving orientation of orbits, then we get from~\ref{enum:shift:local_existence} that $1+\fld(\alpha) >0$, i.e. $\alpha\in\DFunc$.
Conversely, suppose $\alpha\in\DFunc$.
Then by~\ref{enum:shift:local_existence}, $\flow_{\alpha}$ is a local diffeomorphism preserving orientation of orbits, so it suffices to show that $\dif$ is a bijection.
This is trivial for singular points.
If $\orb$ is a closed orbit, then $\flow_{\alpha}:\orb\to\orb$ is a local diffeomorphism (i.e.\! a covering map) homotopic to $\id_{\orb}$, whence $\flow_{\alpha}|_{\orb}$ is a bijection.
Finally, if $\orb$ is a non-closed orbit, then $\flow_{\alpha}:\orb\to\orb$ is a proper (as $\Mman$ is now compact) monotone function, and so it is a bijection as well.

\ref{enum:shift:isotopy_between_shifts}
If $\fld(\alpha)>-1$ and $\fld(\mu)=0$, then $\fld(\mu\alpha) = \mu \fld(\alpha) + \fld(\mu)\alpha  = \mu \fld(\alpha) >-1$, since $0\leq \mu\leq 1$.
\end{proof}

\begin{lemma}\label{lm:shifts_and_covering_maps}
Let $p:\tMman\to\Mman$ be a covering map.
Then 
\begin{enumerate}[label={\rm(\alph*)}]
\item\label{enum:lm:shifts_and_covering_maps:hflow}
 $\flow$ lifts to a flow $\hflow:\tMman\times\bR\to\tMman$ such that $\flow_t\circ p = p\circ \hflow_t$ for all $t\in\bR$;
\item\label{enum:lm:shifts_and_covering_maps:orbits}
for every orbit $\cov{\gamma}$ of $\hflow$ its image $p(\cov{\gamma})$ is an orbit of $\flow$;
\item\label{enum:lm:shifts_and_covering_maps:shifts}
for each function $\alpha:\Mman\to\bR$ the map $\hflow_{\alpha\circ p}$ is a lifting of $\flow_{\alpha}$, i.e. $p\circ \hflow_{\alpha\circ p} = \flow_{\alpha} \circ p$;

\item\label{enum:lm:shifts_and_covering_maps:xi}
\jadded{34.8}{$\hflow$ commutes with each covering transformations $\xi:\tMman\to\tMman$, i.e. $\hflow_t\circ\xi = \xi\circ\hflow_t$ for all $t\in\bR$ which can also we rewritten as follows:}
\begin{align}\label{equ:xi_commutes_with_hflow}
\hflow_t(\xi(z),t) &= \xi\circ\hflow(z,t), \qquad (z,t)\in\tMman\times\bR.
\end{align}
\end{enumerate}
\end{lemma}
\begin{proof}
\ref{enum:lm:shifts_and_covering_maps:hflow}
Define a \myemph{homotopy with ``open ends''} $q:\tMman\times\bR\to\Mman$ by $q(x,t) = \flow(p(x),t)$.
Since $q_0=p\circ\flow_0 = \id_{\tMman}\circ p$, it follows from the homotopy lifting property for covering maps that $q$ lifts to a continuous map $\hflow:\tMman\times\bR\to\tMman$ such that $\hflow_{0}=\id_{\tMman}$ and $q = p \circ\hflow$.
In other words, $\flow(p(x),t) = p\circ \hflow(x,t)$, i.e. $\flow_t\circ p = p\circ \hflow_t$ for all $t\in\bR$.
One easily checks that $\hflow_{t+s} = \hflow_t \circ \hflow_s$ for all $s,t\in\bR$, which means that $\hflow$ is a flow.
	
\ref{enum:lm:shifts_and_covering_maps:orbits}
Let $y\in\cov{\gamma}$.
Then $p\circ\hflow(y,t)=\flow(p(y),t)$ for all $t\in\bR$, whence $p$ maps the $\hflow$-orbit $\cov{\gamma}=\hflow(y\times\bR)$ of $y$ onto the $\flow$-orbit $\flow(p(y)\times\bR)$ of $p(y)$.

\ref{enum:lm:shifts_and_covering_maps:shifts}
Let $z\in\Mman$.
Then
\begin{equation}\label{equ:lifting_of_shift}
\aligned
p\circ\hflow_{\alpha\circ p}(z) &= 
p\circ\hflow\bigl(z,\alpha\circ p(z)\bigr) = 
\flow\bigl(p(z),\alpha\circ p(z)\bigr)= \\ 
&= \flow\bigl(p(z),\alpha\circ p(z)\bigr) = \flow_{\alpha}\circ p(z),
\endaligned
\end{equation}
which means that $\hflow_{\alpha\circ p}$ is a lifting of $\flow_{\alpha}$.

\ref{enum:lm:shifts_and_covering_maps:xi}
Let $\xi:\tMman\to\tMman$ be any covering transformation, i.e.\! a homeomorphism such that $p\circ\xi=p$.
It suffices to show that the map $\hflow':\tMman\times\bR\to\tMman$ defined by $\hflow'_t = \xi\circ\hflow_t\circ\xi^{-1}$, $t\in\bR$, coincides with $\hflow$.
Notice that
\[
p\circ\hflow'_t = p\circ (\xi\circ\hflow_t\circ\xi^{-1}) = p\circ \hflow_t\circ\xi^{-1}= \flow_t\circ p\circ\xi^{-1} =\flow_t\circ p,
\]
i.e. $\hflow'$ is also a lifting of $\flow$.
Moreover, since $\hflow'_0 = \xi\circ\hflow_0\circ\xi^{-1} = \id_{\tMman} =\hflow_0$, it follows from the uniqueness of liftings that $\hflow=\hflow'$.
\end{proof}

\begin{lemma}\label{lm:shifts:covering}
Let $p:\tMman\to\Mman$ is a \myemph{regular} covering map, i.e. there is a properly discontinuous action of some group $G$ on $\tMman$ such that $p$ induces a homeomorphism $\tMman/G \cong \Mman$.
Suppose also that $\tMman$ and $\Mman$ are connected, and that all orbits of the lifted flow $\hflow:\tMman\times\bR\to\tMman$ are non-singular and non-closed.
Let also $\dif\in\Diff(\Mman)$ and $\hdif\in\Diff(\tMman)$ be a lifting of $\dif$, i.e. $p\circ\hdif=\dif\circ p$.
If $\hdif$ preserves each orbit of $\hflow$ then the following statements hold.
\begin{enumerate}[label={\rm{(\alph*)}}]
\item\label{enum:shift:covering:h_pres_orb}
$\dif$ preserves each orbit of $\flow$.

\item\label{enum:shift:covering:th_tfta}
There exists a unique $\Cinfty$ function $\halpha:\tMman\to\bR$ such that $\hdif=\hflow_{\halpha}$.

\item\label{enum:shift:covering:haxi_ha}
\xadded{1}{0}{28.1-28.2}{}%
\xadded{2}{0}{32.1-32.3}{}%
Let $\xi\in G$ and $\ppy\in\tMman$.
Then the following conditions are equivalent
\begin{align*}
\text{\xc{1}} &~ \hdif\circ\xi(\ppy)=\xi\circ\hdif(\ppy), &
\text{\xc{2}} &~ \halpha\circ\xi(\ppy)=\halpha(\ppy), \\
\text{\xc{3}} &~ \hdif\circ\xi=\xi\circ\hdif\ \text{on all of $\tMman$}, &
\text{\xc{4}} &~ \halpha\circ\xi=\halpha \ \text{on all of $\tMman$}.
\end{align*}
In particular, if either of these conditions holds for every $\xi\in G$, then $\halpha$ induces a unique $\Cinfty$ function $\alpha:\Mman\to\bR$ such that $\halpha = \alpha\circ p$ and $\dif = \flow_{\alpha}$.

\item\label{enum:shift:covering:flow_has_nonclosed_orbit}
If $\flow$ has at least one non-closed orbit then any other lifting $\hdif'$ of $\dif$ preserving orbits of $\hflow$ coincides with $\hdif$, and $\hdif$ commutes with $G$.
Hence by~\ref{enum:shift:covering:haxi_ha}, $\dif = \flow_{\alpha}$ for some $\Cinfty$ function $\alpha:\Mman\to\bR$.

\item\label{enum:shift:covering:th_has_fixed_pt}
If $\hdif$ has a fixed point $\ppy$, then again it commutes with $G$, and by~\ref{enum:shift:covering:haxi_ha}, $\dif = \flow_{\alpha}$ for some $\Cinfty$ function $\alpha:\Mman\to\bR$.
Moreover, in this case $\alpha(p(\ppy))=0$.

\item\label{enum:shift:covering:lift}
Suppose each orbit $\omega$ of $\flow$ is closed and its inverse image $p^{-1}(\omega)$ is connected, whence it must be a non-closed orbit of $\hflow$.
Suppose also that there exists an orbit $\gamma$ of $\flow$ and a $\Cinfty$ function $\alpha:\gamma\to\bR$ such that $\dif(x) = \flow(x,\alpha(x))$ for all $x\in\gamma$.
Then $\alpha$ extends to a unique $\Cinfty$ function $\alpha:\Mman\to\bR$ such that $\dif= \flow_{\alpha}$.
\end{enumerate}
\end{lemma}
\begin{proof}
\ref{enum:shift:covering:h_pres_orb}
Let $\gamma$ be an orbit of $\flow$, and $\cov{\gamma}$ be an $\hflow$-orbit with $p(\tilde{\gamma}) = \gamma$.
Then $\dif(\gamma)=\dif\circ p(\cov{\gamma}) = p\circ\hdif(\cov{\gamma}) =p(\cov{\gamma})=\gamma$.

\ref{enum:shift:covering:th_tfta} 
Existence of $\halpha$ is a direct consequence of Lemma~\ref{lm:shift}\ref{enum:shift:shift_func_existence} since by assumption all orbits of $\hflow$ are non-singular and non-closed.

\ref{enum:shift:covering:haxi_ha}
First notice that $\hdif(\ppy)=\hflow(z,\halpha(\ppy))$ and
\begin{align}\label{equ:}
\xi^{-1}\circ\hdif\circ\xi(\ppy) = 
\xi^{-1}\circ\hflow( \xi(\ppy), \halpha\circ\xi(\ppy) ) \stackrel{\eqref{equ:xi_commutes_with_hflow}}{=}
\xi^{-1}\circ\xi\circ\hflow(\ppy, \halpha\circ\xi(\ppy) ) =
\hflow(z, \halpha\circ\xi(\ppy) ). 
\end{align}

\xc{1}$\Leftrightarrow$\xc{2}
The latter identity shows that the relation $\hdif(\ppy)=\xi^{-1}\circ\hdif\circ\xi(\ppy)$ is equivalent to $\hflow(\ppy,\halpha(\ppy)) = \hflow(\ppy, \halpha\circ\xi(\ppy) )$, which in turn is equivalent to the assumption $\halpha\circ\xi(\ppy)=\halpha(\ppy)$ since $\ppy$ is non-periodic and non-fixed for $\hflow$.

The equivalence \xc{3}$\Leftrightarrow$\xc{4} follows from \xc{1}$\Leftrightarrow$\xc{2}.
Moreover, we also have obvious implications \xc{3}$\Rightarrow$\xc{1} and \xc{4}$\Rightarrow$\xc{2}.

\xc{1}$\Rightarrow$\xc{3}
Suppose $\hdif\circ\xi(z) = \xi\circ\hdif(z)$ at some $z\in\Mman$.
Then $\hdif\circ\xi$ and $\xi\circ\hdif$ are two liftings of $\dif$ which coincide at $\ppy\in\tMman$, whence they coincide on all of $\tMman$.

\ref{enum:shift:covering:flow_has_nonclosed_orbit}
Let $\gamma$ be a non-closed orbit of $\flow$, and $\tilde{\gamma}$ an orbit of $\hflow$ with $p(\tilde{\gamma}) = \gamma$.
Then $p$ induces a continuous bijection of $\tilde{\gamma}$ onto $\gamma$, whence $\dif\circ p = p \circ\hdif = p \circ\hdif'$ implies that $\hdif'|_{\tilde{\gamma}} = \hdif|_{\tilde{\gamma}}$.
Since $\tMman$ is connected, it follows that the liftings $\hdif'$ and $\hdif$ of $\dif$ coincide on all of $\tMman$.

In particular, $\xi^{-1}\circ\hdif\circ\xi$ is a lifting of $\dif$ which preserves orbits of $\hflow$, whence $\hdif=\xi^{-1}\circ\hdif\circ\xi$, i.e. $\dif\circ\xi=\hdif\circ\xi$.

\ref{enum:shift:covering:th_has_fixed_pt}
If $\hdif(z)=z$ for some $z\in\tMman$, then for each $\xi\in G$, $\hdif' = \xi^{-1}\circ\hdif\circ\xi$ is also a lifting of $\dif$ satisfying $\hdif'(z) = z = \hdif(z)$, whence $\hdif'\equiv \hdif$, i.e. $\hdif$ commutes with $\xi$.

Moreover, $\hdif(z)= z = \hflow(z,\halpha(z))$ implies that $\halpha(z)=0$.
Hence $\alpha(p(z))=\halpha(z)=0$.

\ref{enum:shift:covering:lift}
Take any point $\ppy\in \cov{\gamma} := p^{-1}(\gamma)$ and put $\ppx=p(\ppy)\in\gamma$, $a=\alpha(\ppx)$, $\ppy'=\hflow_{a}(\ppy) \in\tMman$.
Then 
$p(\ppy') = p\circ \hflow_{a}(\ppy) = \flow_{a}\circ p(\ppy) = \flow_{a}(\ppx) = \dif(\ppx).$
Since $p$ is a covering map, there exists a unique map $\hgdif:\tMman\to\tMman$ of $\dif\circ p:\tMman\to\Mman$ such that $\hgdif(\ppy) = \ppy'$ and $\dif\circ p = p \circ\hgdif$.

We claim that $\hgdif$ preserves each orbit of $\hflow$.
Indeed, by assumption each orbit of $\hflow$ is of the form $p^{-1}(\omega)$ for some orbit $\omega$ of $\flow$.
Moreover, as $\dif(\omega)=\omega$, the identity $\dif\circ p = p \circ\hgdif$ implies that  $\hgdif^{-1}(p^{-1}(\omega)) = p^{-1}(\dif^{-1}(\omega))=p^{-1}(\omega)$, whence the orbit $p^{-1}(\omega)$ is invariant under $\hgdif$.

Therefore by~\ref{enum:shift:covering:th_tfta}, $\hgdif = \hflow_{\halpha}$ for a unique $\Cinfty$ function $\halpha:\tMman\to\bR$.

Let us show that $\halpha = \alpha\circ p$ on $\cov{\gamma}$.
Indeed, notice that $\hgdif|_{\cov{\gamma}} = (\hflow_{\halpha})|_{\cov{\gamma}}$ is a lifting of $\dif|_{\gamma} = \flow_{\alpha}$.
On the other hand due to~\eqref{equ:lifting_of_shift} the restriction $\hflow_{\alpha\circ p}|_{\cov{\gamma}}$ is also a lifting of $\flow_{\alpha}$.
Moreover, $\hflow(\ppy, \halpha(\ppy)) = \hgdif(\ppy) = \hflow(\ppy, a) =\hflow(\ppy, \alpha\circ p(\ppy))$, i.e. those liftings conicide at $\ppy$, and therefore they are equal on all of $\cov{\gamma}$.
Thus $\hflow(\ppz, \halpha(\ppz)) = \hflow(\ppz, \alpha\circ p(\ppz))$ for all $z\in\cov{\gamma}$.
Since $\cov{\gamma}$ is a non-closed orbit, we obtain that $\halpha = \alpha\circ p$ on $\cov{\gamma}$.

Furthermore, as $\cov{\gamma}$ is invariant under each $\xi\in G$, it follows that $\halpha\circ\xi(y) = \alpha\circ p\circ\xi(y) =\halpha\circ\xi(y)$.
Hence by~\ref{enum:shift:covering:haxi_ha}, $\halpha\circ\xi=\halpha$ on all of $\tMman$, and thus it induces a unique $\Cinfty$ function $\alpha':\Mman\to\bR$ such that $\alpha'\circ p = \halpha$ and $\dif = \flow_{\alpha'}$.

Finally, if $\ppx\in\gamma$ and $\ppy\in p^{-1}(x)$, then $\alpha'(\ppx)  = \alpha'\circ p(\ppy) = \halpha(\ppy)=\alpha\circ p(\ppy) = \alpha(\ppx)$.
In other words, $\alpha'=\alpha$ on $\cov{\gamma}$, i.e. $\alpha'$ is the desired extension of $\alpha$.
\end{proof}

\begin{figure}[ht]
\centering
\includegraphics[height=1.7cm]{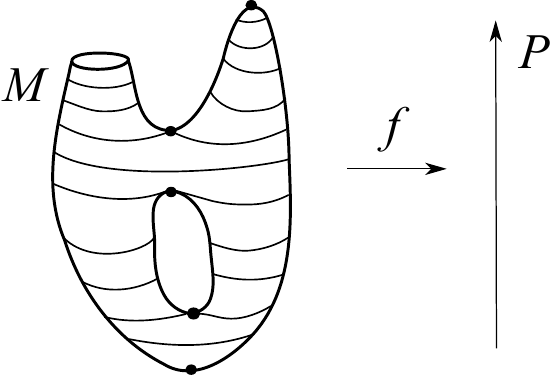}
\caption{Orbits of a Hamiltonian-like vector field}\label{fig:ham_like_vf}
\end{figure}

{\bf Hamiltonian-like flows for $\func\in\FF(\Mman,\Pman)$.}
Till the end of this section assume that $\Mman$ is an orientable \added{possibly non-connected} compact surface, $\func\in\FF(\Mman,\Pman)$, and $\fSing$ is the set of critical points of $\func$.
\begin{definition}\label{def:ham_like_vf}
A smooth vector field $\fld$ on $\Mman$ will be called \myemph{Hamiltonian-like} for $\func$ if the following conditions hold true.
\begin{enumerate}[leftmargin=5ex, topsep=0pt, wide, label={\rm(\alph*)}]
\item\label{enum:HamVF:F_0_crpt}
$\fld(z)=0$ if and only if $z$ is a critical point of $\func$.
\item\label{enum:HamVF:Ff_0}
$\fld(\func)\equiv0$ everywhere on $\Mman$.
\item\label{enum:HamVF:local_form}
Let $z$ be a critical point of $\func$.
Then there exists a local representation of $\func$ at $z$ in the form of a homogeneous polynomial $\gfunc:(\bR^2,0)\to(\bR,0)$ without multiple factors as in Axiom~\AxCrPt, such that in the same coordinates $(x,y)$ near the origin $0$ in $\bR^2$ we have that $\fld = -\gfunc'_{y}\,\tfrac{\partial}{\partial x} + \gfunc'_{x}\,\tfrac{\partial}{\partial y}$.
\end{enumerate}
\end{definition}
Condition~\ref{enum:HamVF:Ff_0} means that $\func$ is constant along orbits of $\fld$.
Then it follows from~\ref{enum:HamVF:F_0_crpt} and Axiom~\AxBd\ that the orbits of $\fld$ are critical points of $\func$, regular components of level sets of $\func$ being closed orbits of $\fld$, and connected components of the sets $\crLev\setminus\fSing$, where $\crLev$ runs over all critical components of level-sets of $\func$\jchanged{22.1}{, see}{. See} Figure~\ref{fig:ham_like_vf}.

By~\cite[Lemma~5.1]{Maksymenko:AGAG:2006} or~\cite[Lemma~16]{Maksymenko:ProcIM:ENG:2010} every $\func\in\FF(\Mman,\Pman)$ admits a Hamiltonian-like vector field.
For the proof take the Hamiltonian vector field $\fld$ for $\func$ with respect to any symplectic form $\orb$ on $\Mman$, and then properly change $\fld$ near each critical point of $\func$ in accordance with~\ref{enum:HamVF:local_form} of Definition~\ref{def:ham_like_vf}.

Let $\fld$ be a Hamiltonian-like vector field for $\func$.
Since $\fld$ is tangent to $\partial\Mman$, $\fld$ generates a flow $\flow:\Mman \times \bR\to \Mman$ which will also be called \myemph{Hamiltonian-like} for $\func$.
Then condition~\ref{enum:HamVF:Ff_0} is equivalent to the assumption that $\func\circ\flow_t = \func$ for all $t\in\bR$, that is $\flow_t\in\fStab$.
In fact, $\func\circ\flow_{\alpha} = \func$ for all $\alpha\in C^{\infty}(\Mman,\bR)$.

\begin{lemma}\label{lm:prop:StabIdfX}
\begin{enumerate}[label={\rm(\alph*)}, wide]
\item\label{enum:lm:hamvf:shift:SIdfX}
Let $\Xman\subset\Mman$ be a closed subset, and $\alpha\in\DFunc$\jremoved{22.2}{be} such that $\alpha=0$ on $\Xman$.
Then $\flow_{\alpha} \in \StabilizerId{\func,\Xman}$.
\item\label{enum:lm:hamvf:shift:saddles}
\jadded{26.3}{Suppose $\Mman$ is connected and $\func$ has a saddle point.
If $\flow_{\alpha}=\flow_{\beta}$ for some $\alpha,\beta\in \Ci{\Mman}{\bR}$, then $\alpha=\beta$.
In particular, in this case for each $\dif\in\StabilizerId{\func}$ there may exist no more than one function $\alpha\in\DFunc$ such that $\dif=\flow_{\alpha}$.}
\end{enumerate}
\end{lemma}
\begin{proof}
\ref{enum:lm:hamvf:shift:SIdfX}
By \ref{enum:shift:isotopy_between_shifts} of Lemma~\ref{lm:shift} $\{\flow_{t\alpha}\}_{t\in\UInt}$ is an isotopy between $\id_{\Mman} = \flow_{0}$ and $\flow_{\alpha}$ in $\Stabilizer{\func}$.
But $t\alpha=0$ on $\Xman$ for all $t\in\UInt$, hence $\flow_{t\alpha}$ is fixed on $\Xman$, and so $\flow_{\alpha} \in \StabilizerId{\func,\Xman}$.

\ref{enum:lm:hamvf:shift:saddles}
It follows from the structure of level-sets of $\func$ at saddle critical points, that $\flow$ has at least one non-closed orbit $\gamma$.
By assumption, $\flow(x,\alpha(x))=\flow(x,\beta(x))$, whence $\alpha(x)=\beta(x)$ for all $x\in\gamma$.
As $\Mman$ is connected, we get from Lemma~\ref{lm:shift}\ref{enum:shift:shift_func_uniqueness} that $\alpha=\beta$ of all of $\Mman$.
\end{proof}

The following theorem explains more precisely the statement~\ref{enum:th:stab_orb_full_info:Sidf} of Theorem~\ref{th:stab_orb:full_info}:
\begin{theorem}\label{th:charact_Stabf}{\rm\jadded{23.1}{(\cite[Theorem~1.3]{Maksymenko:AGAG:2006}, \cite[Theorem~3]{Maksymenko:ProcIM:ENG:2010}).}}
Let $\Mman$ be a connected compact surface, $\func\in\FF(\Mman,\Pman)$, $\flow:\Mman\times\bR\to\Mman$ the flow generated by some Hamiltonian vector field $\fld$, and $\Sh{\fld}:\DFunc \to \StabilizerId{\func}$\jremoved{23.2}{be} the map defined by $\Sh{\fld}(\alpha) = \flow_{\alpha}$.
If $\func$ has at least one \myemph{saddle} or a \myemph{degenerate local extreme}, then $\Sh{\fld}$ is a \myemph{homeomorphism} with respect to $C^{\infty}$ topologies and $\StabilizerId{\func}$ is contractible (because $\DFunc$ is convex).
Otherwise, there is $\theta\in\DFunc$ such that $\Sh{\fld}$ can be represented as a composition
\[
\Sh{\fld}: \DFunc \xrightarrow{~\text{quotient}~}
\DFunc/\{n\theta\}_{n\in\bZ} \xrightarrow{~\text{homeomorphism}~}
\StabilizerId{\func}
\]
of the quotient map by the closed discrete subgroup $Z = \{n\theta\}_{n\in\bZ}$ of $\DFunc$ and a homeomorphism of the quotient of $\DFunc$ by $Z$ onto $\StabilizerId{\func}$.
In particular, $\Sh{\fld}$ is an infinite cyclic covering map and $\StabilizerId{\func}$ is homotopy equivalent to the circle.
\qed
\end{theorem}

\section{Simplification of diffeomorphisms from $\Stabilizer{\func}$}
In what follows we will assume that \added{$\Mman$ is a compact and not necessarily connected smooth surface,} $\func\in\FF(\Mman,\Pman)$, $\fld$ is a Hamiltonian vector field for $\func$, $\flow$ is the flow generated by $\fld$, and $\Xman$ is an $\func$-adapted \myemph{submanifold}.
Thus each connected component of $\Xman$ is either a regular component of some level-set of $\func$ or an $\func$-adapted subsurface.

Let $\Vman$ be another $\func$-adapted submanifold.
The following technical Lemma~\ref{lm:extension_of_shift_functions} implies that if a diffeomorphism $\dif\in\Stabilizer{\func,\Vman}$ has a shift-function $\gamma_{\dif}$ on $\Xman$ which vanishes on $\Vman\cap\Xman$, then $\dif$ can be deformed in $\Stabilizer{\func,\Vman}$ to a diffeomorphism $\dif'$ fixed on $\Vman\cup\Xman$ and even on some $\func$-regular neighborhood of $\Vman\cup\Xman$ as well.
This is a kind of a ``\textit{simplification}'' of $\dif$ in its isotopy class in $\pi_0\Stabilizer{\func,\Vman}$.
In fact such a statement is established in~\cite[Lemma~4.14]{Maksymenko:AGAG:2006} but the presented formulation takes to account a possible continuous dependence of $\gamma_{\dif}$ on $\dif$.

\begin{lemma}\label{lm:extension_of_shift_functions}
{\rm\jadded{23.3}{(cf.~\cite[Lemma~4.14]{Maksymenko:AGAG:2006}).}}
Let $\mathcal{A} \subset \Stabilizer{\func}$ be a subset.
Suppose that there exists a continuous map $\gamma:\mathcal{A} \to C^{\infty}(\Xman,\bR)$ such that
\[
\dif(x) = \flow\bigl(x,\gamma(\dif)(x)\bigr), \qquad  \dif\in\mathcal{A}, \  x\in\Xman.
\]
Then for every pair of $\func$-regular neighborhoods $\Uman$ and $\Nman$ of $\Xman$ such that $\Uman \subset\Int{\Nman}$ there exists a continuous map $\beta: \mathcal{A} \ \to \ \DAFunc{\fld} \  \subset \  C^{\infty}(\Mman,\bR)$ such that
\begin{enumerate}[topsep=0pt, label={\rm(\arabic*)}]
\item\label{enum:cond:extsh:a}
$\beta(\dif) = \gamma(\dif)$ on $\Xman$ for each $\dif\in\mathcal{A}$;
\item\label{enum:cond:extsh:b}
$\dif|_{\Uman} = \flow_{\beta(\dif)}$ for each $\dif\in\mathcal{A}$, i.e. $\dif(x) = \flow\bigl(x,\beta(\dif)(x)\bigr)$, $x\in\Uman$;
\item\label{enum:cond:extsh:c}
$\beta(\dif) = 0$ on $\overline{\Mman\setminus\Nman}$ for each $\dif\in\mathcal{A}$;
\item\label{enum:cond:extsh:d}
if $\Wman \subseteq \Uman$ is an $\func$-regular neighborhood of $\Xman$ or $\Wman=\Xman$, then for each $\dif\in\mathcal{A}$ fixed on $\Wman$ and such that $\gamma(\dif)\equiv0$, we have that $\beta(\dif)=0$ on $\Wman$ as well;
\item\label{enum:cond:extsh:e}
the homotopy
\begin{align}\label{equ:homotopy_H}
H&: \mathcal{A} \times\UInt \to  \Stabilizer{\func}, &
H(\dif,t) = H_t(\dif) = \added{\dif\circ (\flow_{t\beta(\dif)})^{-1}} 
\end{align}
has the following properties:
\begin{enumerate}[label={\rm(\alph*)}, leftmargin=5ex]
\item\label{extsh:enum:Ht_on_U}
\added{$H_t(\dif)|_{\Uman} = \flow_{(1-t)\beta\circ \flow_{t\beta(\dif)}^{-1} }$ for all $\dif\in\mathcal{A}$ and $t\in\UInt$;}

\item\label{extsh:enum:H0_is_id}
 $H_0=\id_{\mathcal{A}}$;

\item\label{extsh:H1_A_SfU}
$H_1(\mathcal{A}) \subset\Stabilizer{\func,\Uman}$;

\item\label{extsh:Ht_ASfU}
if $\dif$ and $\Wman$ are the same as in~\ref{enum:cond:extsh:d}, then $H_t(\dif)$ is fixed on $\Wman$ as well as $\dif$, whence
\xadded{1.4}{0}{24.1}{}
\[ H_t\fremoved{(\dif)}\bigl( \Stabilizer{\func,\Wman}\cap  \mathcal{A}\cap\gamma^{-1}(0)\bigr) \subset \Stabilizer{\func,\Wman}\,;\]

\item\label{extsh:Ht_ADfU}
\added{$H_t\bigl( \mathcal{A}\cap \FolStabilizer{\func} \bigr) \subset \FolStabilizer{\func}$.}
\end{enumerate}
\end{enumerate}
\end{lemma}
\begin{proof}
First notice that there exists a $C^{\infty}$ function $\mu:\Mman\to\UInt$ with the following properties: $\mu=0$ on some neighborhood of $\overline{\Mman\setminus\Nman}$, $\mu=1$ on some neighborhood of $\Uman$, and $\fld(\mu)=0$, see Figure~\ref{fig:reg_nbh_extdif}.
\begin{figure}[!htbp]
\centering
\includegraphics[width=3.5cm]{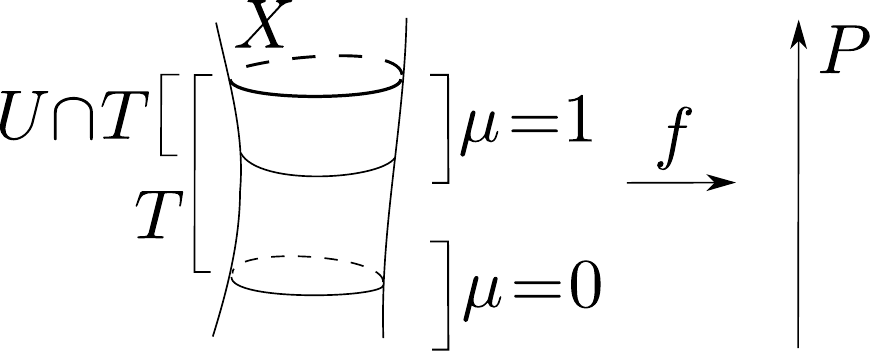}
\protect\caption{}\label{fig:reg_nbh_extdif}
\end{figure}
Indeed, for every connected component ${T}$ of $\Nman\setminus\Xman$, there exists a diffeomorphism $\sigma_{T}:\overline{T}\to \Circle\times\UInt$ such that $\sigma_{T}^{-1}(\Circle\times t)$ is a regular component of some level-set of $\func$, $\sigma_{T}(\partial\Nman\cap{T}) = \Circle\times 0$, $\sigma_{T}(\Xman\cap\overline{T}) = \Circle\times 1$, and $\sigma_T(\Uman\cap\overline{T}) = \Circle\times [\eps_{T}, 1]$ for some $\eps_{T}>0$.

Fix any $C^{\infty}$ function $\mu_{T}:\UInt\to\UInt$ such that $\mu_{T}=0$ on $[0,\eps_{T}/3]$ and $\mu_{T}=1$ on $[2\eps_{T}/3,1]$.
Let also $p:\Circle\times\UInt\to\UInt$ be the projection onto the second coordinate, that is $p(\phi,t)=t$.
Then $\mu:\Mman\to\UInt$ can be defined by the following rule: $\mu=0$ on $\overline{\Mman\setminus\Nman}$, $\mu=1$ on $\Uman$, and $\mu|_{T} = \mu_{T} \circ p\circ \sigma_{T}$ for every connected component ${T}$ of $\Nman\setminus\Xman$.

It follows from Lemma~\ref{lm:shift}\ref{enum:shift:local_existence} that \myemph{for each $\dif\in\mathcal{A}$ the function $\gamma(\dif):\Xman\to\bR$ uniquely extends to a $C^{\infty}$ function $\cov{\gamma}(\dif):\Nman\to\bR$ such that $\dif(x) = \flow(x,\cov{\gamma}(\dif)(x))$ for $x\in\Nman$}.
Moreover, by the formulas for such extension presented in Lemma~\ref{lm:shift}\ref{enum:shift:local_existence} the correspondence $\dif\mapsto\cov{\gamma}(\dif)$ is a \myemph{continuous} map $\cov{\gamma}:\mathcal{A} \to C^{\infty}(\Nman,\bR)$.

Then the map $\beta:\mathcal{A}\to C^{\infty}(\Mman,\bR)$ can be defined as follows: if $\dif\in\mathcal{A}$, then
\begin{equation}\label{equ:beta_homotopy}
\beta(\dif)(x) =
\begin{cases}
\mu(x)\cov{\gamma}(\dif)(x), & x \in \Nman, \\
0, & x\in\Mman\setminus\Nman.
\end{cases}
\end{equation}
Indeed, evidently, $\beta$ is continuous.

Let us check that $\beta(\mathcal{A})\subset\DFunc$.
Indeed, let $\dif\in\mathcal{A}$.
Since $\dif|_{\Nman}$ is a diffeomorphism preserving orientation of orbits of $\fld$ on $\Nman$, it follows Lemma~\ref{lm:shift}\ref{enum:shift:diffeomorphism} that $\fld(\cov{\gamma}(\dif))>-1$ on $\Nman$.
But $\mu$ is constant along orbits of $\fld$, whence by Lemma~\ref{lm:shift}\ref{enum:shift:isotopy_between_shifts}, $\fld(\beta(\dif))>-1$ on all of $\Mman$, that is $\beta(\dif)\in\DFunc$.

Let us check properties~\ref{enum:cond:extsh:a}-\ref{enum:cond:extsh:e} of $\beta$.

\ref{enum:cond:extsh:a} \& \ref{enum:cond:extsh:b}
Since $\mu=1$ on $\Uman \supset \Xman$, it follows from~\eqref{equ:beta_homotopy} that $\beta(\dif) = \cov{\gamma}(\dif)$ on $\Uman$.
In particular,
$\beta(\dif)|_{\Xman} = \cov{\gamma}(\dif)|_{\Xman}=\gamma|_{\Xman}$ and
$\flow(x,\beta(\dif)(x)) = \flow(x,\cov{\gamma}(\dif)(x)) = \dif(x)$ for $x\in\Uman$.

\ref{enum:cond:extsh:c}
$\beta(\dif)=0$ on $\overline{\Mman\setminus\Nman}$ by~\eqref{equ:beta_homotopy}.

\ref{enum:cond:extsh:d}
If $\Wman=\Xman$, the statement follows from~\ref{enum:cond:extsh:a}.
Suppose $\Wman \subseteq \Uman$ is an $\func$-regular neighborhood of $\Xman$, and $\dif\in\mathcal{A}$ is a diffeomorphism fixed on $\Wman$ and such that $\gamma(\dif)\equiv0$.
Then $\id_{\Wman} = \dif|_{\Wman} = \flow_{0}|_{\Wman} = \flow_{\cov{\gamma}(\dif)}|_{\Wman}$.
Thus $0$ and $\cov{\gamma}(\dif)$ are two shift functions for the identity map on $\Wman$ which coincide on $\Xman$.
Since $\Xman$ intersects every connected component of $\Wman$, it follows from Lemma~\ref{lm:shift}\ref{enum:shift:shift_func_uniqueness} that $\cov{\gamma}(\dif)=0$ on $\Wman$.
In particular, $\beta(\dif)|_{\Wman}=\cov{\gamma}(\dif)|_{\Wman}=0$ on $\Wman$ as well.

\ref{enum:cond:extsh:e}
Let $\dif\in\mathcal{A}$.

\added{\ref{extsh:enum:Ht_on_U}
Since $\dif|_{\Uman}=\flow_{\beta(\dif)}$, we get from Lemma~\ref{lm:shift}\ref{enum:shift:composition} that}
\[
\added{H_t(\dif)|_{\Uman} = 
\flow_{\beta(\dif)} \circ \flow_{t\beta(\dif)}^{-1} \stackrel{\eqref{equ:shift_fa_fbinv}}{=}
\flow_{(\beta(\dif) - t\beta(\dif))\circ\flow_{t\beta(\dif)}^{-1}}=
\flow_{(1-t)\beta(\dif) \circ\flow_{t\beta(\dif)}^{-1}}.}
\]

Statements~\ref{extsh:enum:H0_is_id} and~\ref{extsh:H1_A_SfU} are direct consequences of~\ref{extsh:enum:Ht_on_U}. 

\ref{extsh:Ht_ASfU}
Suppose $\dif$ and $\Wman \subset\Uman$ are the same as in~\ref{enum:cond:extsh:d}.
Then $\beta(\dif)=0$ on $\Wman$, and therefore $(1-t)\beta\circ \flow_{t\beta(\dif)}^{-1}=0$ on $\Wman$ as well.
Hence by~\ref{extsh:enum:Ht_on_U}, $H_t(\dif)$ is fixed on $\Wman$.
%
%
%


\ref{extsh:Ht_ADfU}
\added{Notice that $\flow_{t\beta(\dif)} \in \StabilizerId{\func} \subset \FolStabilizer{\func}$ for each $t\in\UInt$.
Therefore if $\dif\in \mathcal{A}\cap \FolStabilizer{\func}$, then the composition $H_t(\dif) = \dif\circ(\flow_{\beta(\dif)})^{-1} \in \FolStabilizer{\func}$ as well.}
\end{proof}

\begin{corollary}\label{cor:SfU_SfX_hom_equiv}
{\rm\jadded{25.1}{(cf.~\cite[Corollary 6.1]{Maksymenko:UMZ:ENG:2012})}}
Let $\Uman$ be an $\func$-regular neighborhood of $\Xman$.
Then the following inclusions of pairs are homotopy equivalences:
\added{
\begin{gather}\label{equ:inclusios:stab}
\bigl( \Stabilizer{\func,\Uman}, \FolStabilizer{\func,\Uman} \bigr) 
\ \subset \ 
  \bigl( \StabilizerNbh{\func,\Xman}, \FolStabilizerNbh{\func,\Xman} \bigr) 
\ \subset \ 
  \bigl( \Stabilizer{\func,\Xman}, \FolStabilizer{\func,\Xman} \bigr), \\
\label{equ:inclusios:stab:isot_id}
\bigl( \StabilizerIsotId{\func,\Uman}, \FolStabilizerIsotId{\func,\Uman} \bigr) 
\ \subset \ 
  \bigl( \StabilizerNbhIsotId{\func,\Xman}, \FolStabilizerNbhIsotId{\func,\Xman} \bigr) 
\ \subset \ 
  \bigl( \StabilizerIsotId{\func,\Xman}, \FolStabilizerIsotId{\func,\Xman} \bigr),
\end{gather}
whence we get isomorphisms of the corresponding sequences
$\seqStab{\func,\Uman} \cong \seqStabNbh{\func,\Xman} \cong \seqStab{\func,\Xman}$ and
$\seqStabIsotId{\func,\Uman} \cong \seqStabNbhIsotId{\func,\Xman} \cong \seqStabIsotId{\func,\Xman}$.
}
\end{corollary}
\begin{proof}
Let $\Nman$ be an $\func$-regular neighborhood of $\Uman$.
Then it is also an $\func$-regular neighborhood of $\Xman$.
Consider the constant (and therefore continuous) map $\gamma:\Stabilizer{\func, \Xman} \to C^{\infty}(\Xman,\bR)$ into the zero function, that is $\gamma(\dif)=0$ for all $\dif\in\Stabilizer{\func, \Xman}$.
Then due to Lemma~\ref{lm:extension_of_shift_functions} applied to the pair $\Uman\subset\Nman$ there exists a homotopy $H:\Stabilizer{\func, \Xman}\times\UInt \to \Stabilizer{\func, \Xman}$ such that $H_0=\id$, $H_t(\Stabilizer{\func,\Wman}) \subset \Stabilizer{\func,\Wman}$ for any $\func$-regular neighborhood $\Wman\subseteq\Uman$ of $\Xman$, and $H_1(\Stabilizer{\func,\Xman})\subset \Stabilizer{\func,\Uman}$.
This means that $H$ is a deformation of $\Stabilizer{\func,\Xman}$ into $\Stabilizer{\func,\Uman}$ which leaves spaces $\Stabilizer{\func,\Uman}$ and $\StabilizerNbh{\func, \Xman}$ invariant.
\added{Moreover, by Lemma~\ref{lm:extension_of_shift_functions}\ref{extsh:Ht_ADfU}, $H$ preserves $\FolStabilizer{\func}$, and therefore deforms $\FolStabilizer{\func,\Xman}$ into $\FolStabilizer{\func,\Uman}$ leaving invariant $\FolStabilizerNbh{\func,\Xman}$ and $\FolStabilizer{\func,\Uman}$.}
Hence the inclusions~\eqref{equ:inclusios:stab} are homotopy equivalences \added{of pairs}.

By Lemma~\ref{lm:Lf_Sf} \added{each group in~\eqref{equ:inclusios:stab:isot_id}} consists of path components of the corresponding group from~\eqref{equ:inclusios:stab}, and therefore must be invariant under the deformation $H$.
Hence it suffices to check that $H_1(\StabilizerIsotId{\func,\Xman})\subset \StabilizerIsotId{\func,\Uman}$ and $H_1(\FolStabilizerIsotId{\func,\Xman})\subset \FolStabilizerIsotId{\func,\Uman}$.

Due to Lemma~\ref{lm:DiffMX_DiffMV} the inclusion $i:\Diff(\Mman,\Uman)\subset\Diff(\Mman,\Xman)$ induces a monomorphism $i_0:\pi_0\Diff(\Mman,\Uman)\to\pi_0\Diff(\Mman,\Xman)$, which means that $\DiffId(\Mman,\Uman) = \Diff(\Mman,\Uman) \cap \DiffId(\Mman,\Xman)$.
Therefore
\begin{align*}
H_1(\StabilizerIsotId{\func,\Xman}) &\ \subset \
\Stabilizer{\func,\Uman} \cap \StabilizerIsotId{\func,\Xman} \ = \ \Stabilizer{\func} \ \cap \ \Diff(\Mman,\Uman) \ \cap \ \DiffId(\Mman,\Xman) = \\
&=
\ \Stabilizer{\func} \cap \DiffId(\Mman,\Uman)  \ =: \ \StabilizerIsotId{\func,\Uman}.
\end{align*}
The proof of the inclusion $H_1(\FolStabilizerIsotId{\func,\Xman})\subset \FolStabilizerIsotId{\func,\Uman}$ is similar.
\end{proof}

{\bf The group $\Stabilizer{\func,\Vman;\Xman}$.}
Given two $\func$-adapted submanifolds $\Vman$ and $\Xman$ of $\Mman$, denote by $\Stabilizer{\func,\Vman;\Xman}$ the subset of $\Stabilizer{\func,\Vman}$ consisting of diffeomorphisms $\dif$ for which there exists a $C^{\infty}$-function $\gamma_{\dif}:\Xman\to\bR$ such that $\gamma_{\dif}=0$ on $\Vman\cap\Xman$ and $\dif|_{\Xman}=\flow_{\gamma_{\dif}}$.
In particular, \jadded{26.2}{and by Lemma~\ref{lm:prop:StabIdfX}\ref{enum:lm:hamvf:shift:SIdfX},
$\dif|_{\Xman} \in \StabilizerId{\func|_{\Xman}, \Vman\cap\Xman}$.}

It follows from formulas in Lemma~\ref{lm:shift}\ref{enum:shift:composition} that $\Stabilizer{\func,\Vman;\Xman}$ is a group.

Let $\FolStabilizer{\func,\Vman;\Xman} := \FolStabilizer{\func}\cap\Stabilizer{\func,\Vman;\Xman}$ \added{and $\GrpKR{\func,\Vman;\Xman}:=\Stabilizer{\func,\Vman;\Xman}/\FolStabilizer{\func,\Vman;\Xman}$.
Then we have the following short exact sequence:
\[
\seqStab{\func,\Vman;\Xman}: \ \FolStabilizer{\func,\Vman;\Xman} \monoArrow \Stabilizer{\func,\Vman;\Xman} \epiArrow \GrpKR{\func,\Vman;\Xman}.
\]}
Notice that there are natural inclusions:
\begin{align}\label{equ:SfVX_SfV_X}
\FolStabilizer{\func,\Vman\cup\Xman} &\subset \FolStabilizer{\func,\Vman;\Xman}, &
\Stabilizer{\func,\Vman\cup\Xman} & \subset \Stabilizer{\func,\Vman;\Xman},
\end{align}
since for $\dif\in\Stabilizer{\func,\Vman\cup\Xman}$ one can take $\gamma_{\dif}$ to be the zero function on $\Xman$.

\begin{corollary}\label{cor:SX_Sidf}
\xadded{1}{0}{26.1}{}
Let $\Xman \subset\Mman$ be a \myemph{connected $\func$-adapted subsurface} containing at least one saddle critical point of $\func$, and $\Vman\subset \Xman$ a possibly empty $\func$-adapted submanifold.
Then the inclusion
\begin{equation}\label{equ:inclusions_DfVX_SfVX}
\bigl(
\Stabilizer{\func,\Xman}, \ \FolStabilizer{\func,\Xman} 
\bigr) 
\ \subset \
\bigl(
\Stabilizer{\func,\Vman; \Xman}, \
\FolStabilizer{\func,\Vman; \Xman}
\bigr)
\end{equation}
is a homotopy equivalence of pairs, whence $\seqStab{\func,\Xman}\cong\seqStab{\func,\Vman;\Xman}$.
\end{corollary}
\begin{proof}
As $\Xman$ contain a saddle point of $\func$, we have a homeomorphism $\Sh{}:\DAFunc{\fld|_{\Xman}} \to \StabilizerId{\func|_{\Xman}}$ from Theorem~\ref{th:charact_Stabf} defined by $\Sh{}(\beta) = \flow_{\beta}$.

Also, by definition for each $\dif\in\Stabilizer{\func,\Vman; \Xman}$ there exists $\gamma_{\dif}\in\Ci{\Xman}{\bR}$ such that $\gamma_{\dif}=0$ on $\Vman$ and $\dif|_{\Xman}=\flow_{\gamma_{\dif}} \in \StabilizerId{\func|_{\Xman},\Vman} \subset \StabilizerId{\func|_{\Xman}} = \Sh{}(\DAFunc{\fld|_{\Xman}})$.
Hence $\dif|_{\Xman} = \flow_{\beta}$ for some $\beta\in\DFunc$.
Since $\Xman$ contain a saddle point of $\func$, it follows from Lemma~\ref{lm:prop:StabIdfX}\ref{enum:lm:hamvf:shift:saddles} that $\beta=\gamma_{\dif}$.
In particular, $\beta=0$ on $\Vman$.

Thus we get a continuous map $\gamma:\Stabilizer{\func,\Vman; \Xman} \to C^{\infty}(\Xman,\bR)$ given by $\gamma(\dif) = \Sh{}^{-1}(\dif|_{\Xman})$, whence by Lemma~\ref{lm:extension_of_shift_functions} there is a deformation $H$ of $\Stabilizer{\func,\Vman;\Xman}$ into $\Stabilizer{\func,\Xman}$.
Moreover, by Lemma~\ref{lm:extension_of_shift_functions}\ref{extsh:Ht_ASfU}, $H$ leaves 
$\Stabilizer{\func,\Xman}$ invariant and by Lemma~\ref{lm:extension_of_shift_functions}\ref{extsh:Ht_ADfU} it also preserves $\FolStabilizer{\func,\Vman;\Xman}$ and $\FolStabilizer{\func,\Xman}$.
This implies that the inclusions~\eqref{equ:inclusions_DfVX_SfVX} are homotopy equivalences.
\end{proof}

\begin{lemma}\label{lm:shift_func_near_crlevel}
{\rm\jadded{27.1}{(\cite{Maksymenko:AGAG:2006, Maksymenko:ProcIM:ENG:2010, Maksymenko:MFAT:2010})}}
Let $\Mman$ be an orientable compact surface, $\func\in\FF(\Mman,\Pman)$, $\crLev$ a non-extremal critical component of some level-set of $\func$, and $\regN{\crLev}$ an $\func$-regular neighborhood of $\crLev$.
Let also \jchanged{27.2}{$\dif\in\Stabilizer{\func}$ be}{$\dif$ be an element of $\Stabilizer{\func}$} such that $\dif(\crLev)=\crLev$.
\begin{enumerate}[wide,label={\rm(\arabic*)}]
\item\label{enum:lm:shift_func_near_crlevel:1-5}
Consider the following conditions on $\dif$:
\begin{enumerate}[label=$\mathrm{(\roman*)}$, leftmargin=11ex]
\item\label{enum:pres:bnd_comp}
every component of $\partial\regN{\crLev}$ is \invhp{\dif};

\item\label{enum:pres:edge}
at least one edge of $\crLev$ is \invhp{\dif};

\item
\label{enum:pres:vertice_with_edges}
$\dif$ fixes every vertex of $\crLev$ and 
each edge of $\crLev$ is \invhp{\dif};

\item
\label{enum:pres:has_shift_func}
$\dif|_{\regN{\crLev}} = \flow_{\alpha}$ for some function $\alpha\in\DAFunc{\fld|_{\regN{\crLev}}}$;

\item
\label{enum:pres:belong_to_StabId}
$\dif|_{\regN{\crLev}} \in \StabilizerId{\func|_{\regN{\crLev}}}$.
\end{enumerate}
Then we have the following implications:
\ref{enum:pres:belong_to_StabId}$\Leftrightarrow$%
\ref{enum:pres:has_shift_func}$\Leftrightarrow$%
\ref{enum:pres:vertice_with_edges}$\Leftrightarrow$%
\ref{enum:pres:edge}$\Rightarrow$%
\ref{enum:pres:bnd_comp}.

\added{In particular, if $\Vman\subset\Mman\setminus\regNK$ is an $\func$-adapted submanifold, then due to~\ref{enum:pres:has_shift_func}, $\Stabilizer{\func,\Vman;\regNK}$ is the subset of $\Stabilizer{\func,\Vman}$ consisting of diffeomorphisms satisfying~\ref{enum:pres:edge}-\ref{enum:pres:belong_to_StabId}.}

\item\label{enum:lm:shift_func_near_crlevel:flat}
If either $\dif\in\DiffIdM$, that is $\dif\in\StabilizerIsotId{\func}=\Stabilizer{\func}\cap\DiffIdM$, or $\regN{\crLev}$ can be embedded into $\bR^2$, then all the above conditions are equivalent.
\end{enumerate}
\end{lemma}
\begin{proof}
\ref{enum:lm:shift_func_near_crlevel:1-5}
\added{Implication~\ref{enum:pres:vertice_with_edges}$\Rightarrow$\ref{enum:pres:edge} is trivial, and
\ref{enum:pres:edge}$\Rightarrow$\ref{enum:pres:vertice_with_edges} is proved in~\cite[Claim~7.1.1]{Maksymenko:AGAG:2006}.}

\xadded{0}{0}{27.3}{
\ref{enum:pres:vertice_with_edges}$\Rightarrow$\ref{enum:pres:bnd_comp}.
Fix a Riemannian metric $\rho$ on $\Mman$ and denote by $\Xman_{-}$ (resp. $\Xman_{+}$) the set of points of $\partial\regNK$ at which gradient lines of $\func$ with respect to $\rho$ are directed inward (resp. outward) $\regNK$.
Then $\Xman_{-}$ and $\Xman_{+}$ are unions of connected components of $\partial\regNK$ and $\dif(\Xman_{\pm})=\Xman_{\pm}$.
Let $\Vman_i$, $i=1,\dots,a$, be all the connected components of $\Xman_{-}$.
We will show that each $\Vman_i$ is \invhp{\dif}.
The proof for $\Xman_{+}$ is similar.

Notice that for every $i$ there exists a unique connected component $\Zman_i$ of $\regNK\setminus\crLev$ containing $\Vman_i$.
Put $\dCyl_i := \overline{\Zman_i} \setminus \fSing$.
Then it is easy to see that 

(a) $\Bman_i:=\dCyl_i\cap\crLev$ is a union of edges of $\crLev$, whence by~\ref{enum:pres:vertice_with_edges}, $\dif(\Bman_i)=\Bman_i$;

(b) $\Bman_i\cap\Bman_j=\varnothing$ for $i\not=j$;

(c) $\dCyl_i$ is diffeomorphic, via some diffeomorphism $\phi_i$, with $\dCyl'_i:=\Circle\times\UInt\setminus \FSet_i$, where $\FSet_i$ is a finite subset of $\Circle\times 1$, so that $\phi_i(\Circle\times0) =\Vman_i$, and $\phi_i\bigl((\Circle\times 1) \setminus \FSet_i\bigr) = \Bman_i$.

Suppose $\dif(\Vman_i)=\Vman_j$ for some $i,j$.
As $\dif$ interchanges connected components of $\regNK\setminus\crLev$, it follows that $\dif(\dCyl_i)=\dCyl_j$, whence
$\dif(\Bman_i)=\dif(\dCyl_i\cap\crLev) = \dCyl_j\cap\crLev = \Bman_j$, and so $i=j$.

Thus~\ref{enum:pres:vertice_with_edges} implies that $\dif(\Vman_i)=\Vman_i$ for all $i$, whence $\dif(\dCyl_i)=\dCyl_i$ as well.
Moreover, again due to~\ref{enum:pres:vertice_with_edges}, $\phi_i\circ\dif\circ\phi_i^{-1}$ is a diffeomorphism of $\dCyl'_i$ preserving the arcs of $(\Circle\times1)\setminus \FSet_i$ with their orientation.
This implies that $\phi_i\circ\dif\circ\phi_i^{-1}$ also preserves orientation of $\Circle\times0$, whence $\dif$ preserves orientation of $\phi_i(\Circle\times0)=\Vman_i$.}

\added{
Equivalence \ref{enum:pres:belong_to_StabId}$\Leftrightarrow$\ref{enum:pres:has_shift_func} follows from Theorem~\ref{th:charact_Stabf}.
Implications \ref{enum:pres:has_shift_func}$\Rightarrow$\ref{enum:pres:bnd_comp}-\ref{enum:pres:vertice_with_edges} hold because $\flow_{\alpha}$ preserves all orbits of $\flow$ in $\regN{\crLev}$ with their orientation, and in particular this is true for vertices and edges of $\crLev$ and connected components of $\partial\regNK$.
}
\removed{Implications (iii)$\Rightarrow$(ii)$\Rightarrow$(i) are trivial, and (ii)$\Rightarrow$(iii) is proved in~[35, Claim~7.1.1].}

The implication \ref{enum:pres:vertice_with_edges}$\Rightarrow$\ref{enum:pres:has_shift_func} is very technical and is established in \cite[Lemma~6.4]{Maksymenko:AGAG:2006} for Morse functions and in~\cite[Claim~1]{Maksymenko:ProcIM:ENG:2010} for functions from $\FF(\Mman,\Pman)$.
\jadded{27.5-28.2}{}%
We will just briefly indicate the arguments.
Since every edge $e$ of $\crLev$ is a non-closed orbit of $\flow$, it follows that the values of $\alpha$ on $e$ are uniquely defined.

We claim that $\alpha$ uniquely extends to $\added{\Xman:=} \regN{\crLev}\setminus\fSing$ so that $\dif(x) = \flow(x,\alpha(x))$ for all \changed{$x\in\regN{\crLev}\setminus\fSing$}{$x\in\Xman$}.
Indeed, let \changed{$p:\widetilde{\regN{\crLev}}\to\regN{\crLev}$}{$p:\tXman\to\Xman$} be the universal cover of \changed{$\regN{\crLev}$}{$\Xman$} and $\hflow$ be the lifting of the flow $\flow$ on 
\changed{$\widetilde{\regN{\crLev}}$}{$\tXman$.
Notice that every closed orbit $\gamma$ of $\flow$ in $\Xman$ is not null-homotopic.
Since $\tXman$ is the universal covering space of $\Xman$, it follows that $p^{-1}(\gamma)$ consists of non-closed orbits of $\hflow$.
In other words, all orbits of $\hflow$ are non-singular and non-closed.
}

One deduces from~\ref{enum:pres:vertice_with_edges} that $\dif|_{\regN{\crLev}}$ is homotopic to $\id_{\regN{\crLev}}$ via a homotopy preserving orbits of $\flow$.
\jchanged{27.4}{This implies that}{Then by homotopy lifting property} $\dif$ lifts to a diffeomorphism $\hdif$ of \changed{$\widetilde{\regN{\crLev}}$}{$\tXman$} preserving orbits of $\hflow$\jadded{27.5}{, so $\dif\circ p = p \circ\hdif$.
Moreover, as $\flow$ has a non-closed orbit, it follows from Lemma~\ref{lm:shifts:covering}\ref{enum:shift:covering:flow_has_nonclosed_orbit} that $\dif|_{\Xman}=\flow_{\alpha'}$ for some $\Cinfty$ function $\alpha':\Xman\to\bR$.
This function must coincide with $\alpha$ on non-closed orbits of $\flow$, and therefore it is the desired extension of $\alpha$.
}

The principal problem is to show that $\alpha$ can be defined on $\crLev\cap\fSing$ to become $C^{\infty}$ on all of $\regN{\crLev}$.
This is proved in~\cite{Maksymenko:TA:2003} for non-degenerate critical points, and extended in~\cite{Maksymenko:CEJM:2009, Maksymenko:hamv2} to singularities of homogeneous polynomials $\bR^2\to\bR$ without multiple factors.

\ref{enum:lm:shift_func_near_crlevel:flat}
The implication \ref{enum:pres:bnd_comp}$\Rightarrow$\ref{enum:pres:edge} is also nontrivial and established in \cite[Theorem~7.1]{Maksymenko:AGAG:2006} under assumption that either $\dif\in\DiffIdM$ or $\regN{\crLev}$ can be embedded into $\bR^2$.
A similar statement is proved in E.~Kudryavtseva~\cite[Lemma~3.3]{Kudryavtseva:MatSb:ENG:2013} under assumption that $\dif$ trivially acts on the homologies of $\Mman$.
\end{proof}

\section{Proof of Lemma~\ref{lm:reduction_Mconn_V_dM}}\label{sect:proof:lm:reduction_Mconn_V_dM}
Let $\Mman$ be a compact possibly non-connected surface, $\func\in\FF(\Mman,\Pman)$, $\Vman$ an $\func$-adapted submanifold, $\regN{}$ an $\func$-regular neighborhood of $\Vman$, $\Xman_1,\ldots,\Xman_{\cnt}$ all the connected components of $\overline{\Mman \setminus \regN{}}$, and $\hXman_i := \Xman_i \cap \regN{}$, $i=1,\ldots,\cnt$.
We should construct isomorphisms of short exact sequences $\seqStab{\func,\Vman}\cong\prod_{i=1}^{\cnt}\seqStab{\func|_{\Xman_i},\hXman_i}$ and $\seqStabIsotId{\func,\Vman}\cong\prod_{i=1}^{\cnt}\seqStabIsotId{\func|_{\Xman_i},\hXman_i}$.
\jadded{29.8-29.9}{}

Let $\Uman$ be an $\func$-regular neighborhood of $\regN{}$.
Then for each $i=1,\ldots,\cnt$, $\hUman_i := \Xman_i\cap \Uman$ is an $\func$-regular neighborhood of $\hXman_i$ in $\Xman_i$.
Hence by Corollary~\ref{cor:SfU_SfX_hom_equiv}, we have isomorphisms 
\begin{align*}
\seqStab{\func,\Vman}&\cong\seqStab{\func,\regN{}}\cong\seqStab{\func,\Uman}, & \seqStab{\func|_{\Xman_i},\hXman_i}&\cong\seqStab{\func|_{\Xman_i},\hUman_i}
\end{align*}
inducing isomorphisms of the corresponding $\DSG$-sequences.
Therefore it suffices to construct isomorphisms $\seqStab{\func,\Uman}\cong\prod_{i=1}^{\cnt}\seqStab{\func|_{\Xman_i},\hUman_i}$ and $\seqStabIsotId{\func,\Uman}\cong\prod_{i=1}^{\cnt}\seqStabIsotId{\func|_{\Xman_i},\hUman_i}$.

Notice that we have an \myemph{isomorphism of topological groups} $\rIso:\Diff(\Mman,\Uman) \to \prod_{i=1}^{\cnt}\Diff(\Xman_i,\hUman_i)$ defined by the formula: $\rIso(\dif) = \bigl(\dif|_{\Xman_1}, \ldots, \dif|_{\Xman_{\cnt}}\bigr)$.
Indeed, it is evident that $\rIso$ is a continuous homomorphism of groups.
Moreover, as $\Uman \cup (\cup_{i=1}^{\cnt}\Xman_i) =\Mman$, every $\dif\in \Diff(\Mman,\Uman)$ is completely determined by its restrictions on each $\Xman_i$, which implies 
\myemph{injectivity} of $\rIso$.
Also, since $\Uman$ contains an open neighborhood of each intersection $\regNK\cap\Xman_i$, every $\cnt$-tuple of diffeomorphisms $(\dif_1,\ldots,\dif_{\cnt}) \in \prod_{i=1}^{\cnt}\Diff(\Xman_i,\hUman_i)$ being fixed on some neighborhood of $\regNK\cap\bigl(\cup_{i=1}^{\cnt}\Xman_i\bigr)$ extends by the identity to a unique diffeomorphism $\dif\in\Diff(\Mman,\Uman)$.
Then $\rIso(\dif)=(\dif_1,\ldots,\dif_{\cnt})$, whence $\rIso$ is a bijection.
Verification continuity of the inverse $\rIso^{-1}$ is left for the reader.

Evidently, $\dif\in\Diff(\Mman,\Uman)$ preserves $\func$ iff each of the restrictions $\dif|_{\Xman_i}$ preserves $\func|_{\Xman_i}$, i.e. $\rIso\bigl(\Stabilizer{\func,\Uman}\bigr)=\prod_{i=1}^{\cnt}\Stabilizer{\func|_{\Xman_i},\hUman_i}$.
Moreover, as each $\Xman_i$ is a union of connected components of level sets of $\func$, it directly follows from the definitions, that $\rIso\bigl(\FolStabilizer{\func,\Uman}\bigr)=\prod_{i=1}^{\cnt}\FolStabilizer{\func|_{\Xman_i},\hUman_i}$.
Hence $\rIso$ induces an isomorphism $\seqStab{\func,\Uman}\cong\prod_{i=1}^{\cnt}\seqStab{\func|_{\Xman_i},\hUman_i}$.

Finally, since $\rIso$ is a \myemph{homeomorphism}, $\rIso\bigl(\DiffId(\Mman,\Uman)\bigr) = \prod_{i=1}^{\cnt}\DiffId(\Xman_i,\hUman_i)$, whence it maps the intersections 
$\StabilizerIsotId{\func,\Uman}=\Stabilizer{\func,\Uman}\cap \DiffId(\Mman,\Uman)$ and 
$\FolStabilizerIsotId{\func,\Uman}=\FolStabilizer{\func,\Uman}\cap \DiffId(\Mman,\Uman)$ onto the corresponding products $\prod_{i=1}^{\cnt}\StabilizerIsotId{\func|_{\Xman_i},\hUman_i}$ and $\prod_{i=1}^{\cnt}\FolStabilizerIsotId{\func|_{\Xman_i},\hUman_i}$ respectively.
Hence $\rIso$ yields as well an isomorphism $\seqStabIsotId{\func,\Uman}\cong\prod_{i=1}^{\cnt}\seqStabIsotId{\func|_{\Xman_i},\hUman_i}$.

The last statement about conjugations of subgroups $\zB{i}$ is a direct consequence of the construction of $\alpha$ and we leave it for the reader.

\section{Sketch of proof of Theorem~\ref{th:stab:chi_neg}}\label{sect:proof:th:pi0SfX:decomp}
This theorem extends \cite[Theorem~1.7]{Maksymenko:MFAT:2010} which in turn extends a technique developed in the paper \cite{JacoShalen:Topology:1977} by W.~Jaco and P.~Shalen on deformations of incompressible subsurfaces.
For the completeness of exposition we will just indicate principal steps of the proof and reformulate them in terms of the present paper. 

Let $\func\in\FF(\Mman,\Pman)$, $\Vman\subset\partial\Mman$ be an $\func$-adapted submanifold, 
$\crLev_i$, $i=1,\ldots,a$, \jremoved{30.1}{be} all non-extremal critical components of level sets of $\func$ whose canonical neighborhoods have negative Euler characteristic.
Let also $\regN{i}\subset\Int{\Mman}$ be an $\func$-regular neighborhood of $\crLev_i$ chosen so that $\regN{i}\cap\regN{j}=\varnothing$ for $i\not=j$, $\canN{i}$ the canonical neighborhood of $\crLev_i$ obtained by attaching to $\regN{i}$ all connected components of $\overline{\Mman\setminus\regN{i}}$ which are $2$-disks, $\regN{} = \cup_{i=1}^{a} \regN{i}$, $\Xman_1,\ldots,\Xman_{\cnt}$ all the connected components of $\overline{\Mman\setminus\added{\Vman\cup\regN{}}}$, and $\hXman_i = \Xman_i\cap(\added{\Vman\cup\regN{}})$.
We should show that $\Xman_1,\ldots,\Xman_{\cnt}$ satisfy statement of Theorem~\ref{th:stab:chi_neg}.

Statement~\ref{enum:pi0SfX:decomp:1} that each $\Xman_i$ is either a $2$-disk or an annulus or a M\"obius strip in established in~\cite[Corollary~6.8]{Maksymenko:MFAT:2010}.
For the proof of~\ref{enum:pi0SfX:decomp:2} we assume that $\Mman$ is orientable.
Non-orientable case can be deduced from orientable by passing to the oriented double covering of $\Mman$ using technique from~\cite[\S4.1]{Maksymenko:AGAG:2006}.
\jadded{30.2-31.2}{}

\begin{lemma}\label{lm:SisotIdf_deform}{\rm(cf. \cite[Claim~7.2]{Maksymenko:MFAT:2010})}
Every $\dif\in\StabilizerIsotId{\func}$ preserves each vertex and every edge of $\crLev_i$ with its orientation for all $i=1,\ldots,a$.
Moreover, if $\dif$ is fixed on $\Vman\cup\regN{}$, then it is isotopic to $\id_{\Mman}$ rel. $\Vman\cup\regN{}$.
\qed
\end{lemma}
In fact this statement is proved for empty $\Vman$, however since $\regN{}$ is disjoint from $\Vman$, the proof can be modified to hold for arbitrary $\Vman$ due to~\cite[Proposition~4.5(B)]{Maksymenko:MFAT:2010}.

Evidently, the first sentence of Lemma~\ref{lm:SisotIdf_deform} says that $\dif$ satisfies condition~\ref{enum:pres:vertice_with_edges} of Lemma~\ref{lm:shift_func_near_crlevel} for every $i=1,\ldots,a$.
Therefore it also satisfies condition~\ref{enum:pres:has_shift_func}, i.e. $\dif|_{\regN{i}} = \flow_{\alpha_i}$ for some $\Cinfty$ function $\alpha_i:\regN{i}\to\bR$.
Since $\{\regN{i}\}$ are mutually disjoint, those functions $\{\alpha_i\}$ determine a $\Cinfty$ function $\alpha:\regN{}=\cup_{i=1}^{a} \regN{i}\to\bR$ such that $\dif|_{\regN{}} = \flow_{\alpha}$.
This means that $\StabilizerIsotId{\func}\subset\Stabilizer{\func;\regN{}}$, and, in particular, $\StabilizerIsotId{\func,\Vman}\subset\Stabilizer{\func,\Vman;\regN{}}$.
Moreover, the second sentence of the lemma means that $\StabilizerIsotId{\func,\Vman}\cap\Stabilizer{\func,\Vman\cup\regN{}} \subset \StabilizerIsotId{\func,\Vman\cup\regN{}}$.
Since the inverse inclusion holds trivially, we obtain that $\StabilizerIsotId{\func,\Vman}\cap\Stabilizer{\func,\Vman\cup\regN{}} = \StabilizerIsotId{\func,\Vman\cup\regN{}}$.

Also, as $\StabilizerIsotId{\func,\Vman}$ consists of path components of $\Stabilizer{\func,\Vman}$, it follows that the deformation of $\Stabilizer{\func,\Vman;\regN{}}$ into $\Stabilizer{\func,\Vman\cup\regN{}}$ constructed in  Corollary~\ref{cor:SX_Sidf} induces a deformation of $\StabilizerIsotId{\func,\Vman}$ into $\StabilizerIsotId{\func,\Vman} \cap \Stabilizer{\func,\Vman\cup\regN{}} = \StabilizerIsotId{\func,\Vman\cup\regN{}}$.
Moreover, by Lemma~\ref{lm:extension_of_shift_functions}\ref{extsh:Ht_ADfU}, that deformation also preserves the intersections with $\FolStabilizer{\func}$.
This implies that the inclusion of pairs
\[ 
\bigl( \StabilizerIsotId{\func,\Vman\cup\regN{}}, \FolStabilizerIsotId{\func,\Vman\cup\regN{}} \bigr) 
\ \subset \
\bigr(\StabilizerIsotId{\func, \Vman}, \FolStabilizerIsotId{\func, \Vman}\bigr)
\]
is a homotopy equivalence inducing isomorphisms of the corresponding $\funcSeq'$-sequences.
\qed

\section{Proof of Theorem~\ref{th:stab:annulus}}\label{sect:proof:th:cyl:func}
First we will define a map~\eqref{equ:map_A_lift} below and describe its properties in Lemma~\ref{lm:DiffIdA}.
This will be used in the proofs of all three cases of the theorem.

Let $(\Cylinder,\Vman,\Wman)=(\Circle\times\UInt, \Circle\times 0, \Circle\times 1)$, $p:\bR\times\UInt \to \Cylinder$ be the universal covering map given by $p(y,s) = (e^{2\pi i y}, s)$, and $\xi:\bR\times\UInt \to\bR\times\UInt$, $\xi(y,s) = (y+1,s)$, the ``shift'' generating the group $\bZ$ of covering transformations.

For a diffeomorphism $\dif\in\Diff(\Cylinder,\Vman)$ we will denote by $\lf{\dif}:\bR\times\UInt \to \bR\times\UInt$ its unique lifting fixed on $\bR\times 0$.
Then $\xi^{-1}\circ\lf{\dif}\circ\xi$ is also a lifting of $\dif$ fixed on $\bR\times0$, whence it coincides with $\lf{\dif}$, which means that $\lf{\dif}$ commutes with $\xi$, c.f.~Lemma~\ref{lm:shifts:covering}\ref{enum:shift:covering:th_has_fixed_pt}.

Now let $\flow:\Cylinder\times\bR\to\Cylinder$ be any flow on $\Cylinder$ such that $\Wman$ is a periodic orbit of $\flow$ of period $1$, and $\hflow:(\bR\times\UInt)\times\bR\to\bR\times\UInt$ be its lifting.
Then $\bR\times1$ a non-closed orbit of $\hflow$.
Moreover, if $\dif\in\Diff(\Cylinder,\Vman)$, then $\dif(\Wman)=\Wman$, whence $\lf{\dif}(\bR\times 1) = \bR\times 1$.
Since $\lf{\dif}$ commutes with $\xi$, we get from Lemma~\ref{lm:shifts:covering}\ref{enum:shift:covering:haxi_ha}
that there exists a \myemph{unique} $\Cinfty$ function $\gamma(\dif):\Wman\to\bR$ such that 
\begin{align}\label{equ:h_lfh}
\dif|_{\Wman} &= \flow_{\gamma(\dif)}, &
\lf{\dif}|_{\bR\times1} &= \hflow_{\gamma(\dif)\circ p}.
\end{align}
Notice that the correspondence $\dif\mapsto\gamma(\dif)$ is a continuous map 
\begin{equation}\label{equ:map_A_lift}
 \gamma:\Diff(\Cylinder,\Vman) \to C^{\infty}(\Wman,\bR).
\end{equation}
\jadded{32.6-32.9}{}

One can assume that 
$\flow(z,1, t) = (ze^{2\pi i t}, 1)$ and $\hflow(y,1,t)=(y+t,1)$ for all $z\in \Circle$ and $y,t\in\bR$.
Then $\xi(y,1)=\flow_{1}(y,1)$ and formulas~\eqref{equ:h_lfh} mean that
\begin{align}\label{equ:h_lfh_formulas}
\dif(z,1) &= (ze^{2\pi i \gamma(\dif)}, 1), &
\hflow(y,1,t)&=(y+\gamma(\dif)(e^{2\pi i y}),1).
\end{align}
\begin{lemma}\label{lm:DiffIdA}
\begin{enumerate}[wide, label={\rm(\alph*)}, itemsep=0.8ex]
\item\label{enum:lm:DiffIdA:lf_composition}
If $\gdif,\dif\in\Diff(\Cylinder,\Vman)$, then $\lf{\dif\circ\gdif} = \lf{\dif} \circ \lf{\gdif}$, $\gamma(\dif\circ\gdif)=\gamma(\dif)\circ\gdif + \gamma(\gdif)$, $\gamma(\gdif^{-1})=-\gamma(\gdif)\circ\gdif^{-1}$, and 
$\gamma(\dif\circ\gdif^{-1}) = (\gamma(\dif)-\gamma(\gdif))\circ\gdif^{-1}$.

\item\label{enum:lm:DiffIdA:Falpha}
Let $\alpha:\Cylinder\to\bR$ be a $\Cinfty$ function such that $\alpha=0$ on $\Vman$ and $\flow_{\alpha}$ is a diffeomorphism, i.e. $\flow_{\alpha}\in\Diff(\Cylinder,\Vman)$.
Then $\gamma(\flow_{\alpha}) = \alpha|_{\Wman}$ and
\begin{equation}\label{equ:gamma_h_fbinv}
\gamma(\dif\circ\flow_{\alpha}^{-1}) = \bigl(\gamma(\dif) - \alpha|_{\Wman}\bigr)\circ\flow_{\alpha}^{-1}|_{\Wman}.
\end{equation}

\item\label{enum:lm:DiffIdA:DidC_charact}
Let $\dif\in\Diff(\Cylinder,\Vman)$.
Then, due to~\eqref{equ:h_lfh_formulas}, $\dif\in \Diff(\Cylinder,\partial\Cylinder)$ iff $\gamma(\dif)$ is constant and takes an integer value.
Moreover, the restriction $\gamma:\Diff(\Cylinder,\partial\Cylinder)\to \bZ$ is a homomorphism whose kernel is $\DiffId(\Cylinder,\partial\Cylinder)$.
In particular, for $\dif\in\Diff(\Cylinder,\Vman)$ the following statements are equivalent:
\begin{align}\label{equ:charact_gamma_0}
\mathrm{(a)}&~\dif\in \DiffId(\Cylinder,\partial\Cylinder); &
\mathrm{(b)}&~\text{$\lf{\dif}$ is fixed on $\bR\times1$}; &
\mathrm{(c)}&~\gamma(\dif)\equiv 0; 
\end{align}
\xadded{-3}{0}{32.10}{}
If $\Uman$ is a regular neighborhood of $\Wman$, then $\dif\in \DiffId(\Cylinder,\Uman\cup\Vman)$ iff $\lf{\dif}$ is fixed on $p^{-1}(\Uman)$.
\end{enumerate}
\end{lemma}
\begin{proof}
\ref{enum:lm:DiffIdA:lf_composition}
Notice that $\lf{\dif\circ\gdif}$ and $\lf{\dif} \circ \lf{\gdif}$ are liftings of $\dif\circ\gdif$ fixed on $\bR\times0$, whence they must coincide due to uniqueness of lifts.
Hence
\begin{align*}
\hflow_{\gamma(\dif\circ\gdif)\circ p} &=
\lf{\dif\circ\gdif}|_{\bR\times 1} = \bigl(\lf{\dif} \circ \lf{\gdif}\bigr)|_{\bR\times 1}  =
 \hflow_{\gamma(\dif)\circ p} \circ  \hflow_{\gamma(\gdif)\circ p} 
 \stackrel{\text{Lemma~\ref{lm:shift}\ref{enum:shift:composition}}}{=} \\ 
&= \hflow_{\gamma(\dif) \circ p \circ \lf{\gdif} + \gamma(\gdif)\circ p} 
= \hflow_{\gamma(\dif) \circ \gdif \circ p + \gamma(\gdif)\circ p}
= \hflow_{\bigl( \gamma(\dif) \circ \gdif + \gamma(\gdif)\bigr)\circ p}.
\end{align*}
Since $\bR\times 1$ is a non-closed orbit of $\hflow$, we get from Lemma~\ref{lm:shift}\ref{enum:shift:shift_func_uniqueness} that $\gamma(\dif\circ\gdif) = \gamma(\dif) \circ \gdif + \gamma(\gdif)$.
The proof of formulas for $\gamma(\gdif^{-1})$ and $\gamma(\dif\circ\gdif^{-1})$ is similar.

\ref{enum:lm:DiffIdA:Falpha}
By~\eqref{equ:lifting_of_shift} the map $\hflow_{\alpha\circ p}$ is a lifting of $\flow_{\alpha}$.
Moreover, $\alpha\circ p(\bR\times0)=\alpha(\Vman)=0$, which implies that $\hflow_{\alpha\circ p}$ is fixed on $\bR\times 0$.
Hence, by the construction, $\gamma(\flow_{\alpha}) = \alpha|_{\Wman}$, and formula~\eqref{equ:gamma_h_fbinv} follows from~\ref{enum:lm:DiffIdA:lf_composition}.

\ref{enum:lm:DiffIdA:DidC_charact}
Let $\tau\in\Diff(\Cylinder,\partial\Cylinder)$ be a Dehn twist along $\Vman$, see~\eqref{equ:Dehn_twist}.
Then it easily follows that $\gamma(\tau)=\pm1$.
Moreover, it is well known that $\pi_0\Diff(\Cylinder,\partial\Cylinder)\cong \bZ$ and this group is generated by the isotopy class of a Dehn twist $\tau$ along $\Vman$.
Hence $\gamma:\Diff(\Cylinder,\partial\Cylinder)\to \bZ$ is surjective and its kernel is the isotopy class of $\tau^0=\id_{\Cylinder}$, i.e it coincides with $\DiffId(\Cylinder,\partial\Cylinder)$.
All other statements are easy and we leave them for the reader.
\end{proof}

{\bf Proof of \ref{eqnu:th:cyl:func:no_cr_pt}.}
Suppose $\func\in\FF(\Cylinder,\bR)$ has no critical points.
Then one can assume that $\func(x,s) = s$, and $\flow(x,s,t) = (xe^{2\pi i t}, s)$, whence $\hflow(y,s,t) = (y+t, s)$ for $x\in \Circle$, $s\in\UInt$, $y,t\in\bR$.
In particular, for each $s\in\UInt$ the circle $\Circle\times s = \func^{-1}(s)$ is a periodic orbit of $\flow$ of period $1$, and $\bR\times s=p^{-1}(\Circle\times s)$ is a non-closed orbit of $\hflow$.

\added{
\ref{enum:th:cyl:func:no_cr_pt:Spr}
We should prove that $\pi_0\Stabilizer{\func,\Vman}=\Stabilizer{\func,\Vman}/\StabilizerId{\func,\Vman}=0$.
It suffices to show that $\Stabilizer{\func,\Vman}\subset\StabilizerId{\func,\Vman}$. 
Evidently, each $\dif\in\Stabilizer{\func,\Vman}$ preserves orbits of $\flow$, while $\lf{\dif}$ preserves the orbits of $\hflow$.
Since $\lf{\dif}$ is fixed on all of $\bR\times 0$, it follows from Lemma~\ref{lm:shifts:covering}\ref{enum:shift:covering:th_has_fixed_pt} that $\dif = \flow_{\alpha_{\dif}}$ and $\hdif=\hflow_{\alpha_{\dif}\circ p}$ for some $\Cinfty$ function $\alpha_{\dif}:\Cylinder\to\bR$ such that $\alpha_{\dif}=0$ on $\Vman = p(\bR\times 0)$.
Hence $\dif = \flow_{\alpha_{\dif}}\in\StabilizerId{\func,\Vman}$ due to Lemma~\ref{lm:prop:StabIdfX}\ref{enum:lm:hamvf:shift:SIdfX}.}

\ref{enum:th:cyl:func:no_cr_pt:S}
Notice that by the construction, $\alpha_{\dif}|_{\Wman}=\gamma(\dif)$ for each $\dif\in\Stabilizer{\func,\Vman}$.
Hence we have a homomorphism $\gamma: \Stabilizer{\func,\partial\Cylinder} \equiv \Stabilizer{\func,\Vman}\cap\Diff(\Cylinder,\partial\Cylinder)\to \bZ$ from~\eqref{equ:h_lfh_formulas}, associating to each $\dif\in\Stabilizer{\func,\partial\Cylinder}$ the (constant and integer) value of $\alpha_{\dif}$ on $\Wman$.
We claim that it is surjective and its kernel is $\StabilizerId{\func,\partial\Cylinder}$.
This will imply that  $\pi_0\Stabilizer{\func,\partial\Cylinder}=\Stabilizer{\func,\partial\Cylinder}/\StabilizerId{\func,\partial\Cylinder}=\bZ$.

Indeed, let $\tau\in\Stabilizer{\func,\partial\Cylinder}$ be an $\func$-adapted Denh twist along $\Vman$.
Then $\gamma(\tau)=\pm1$, whence $\gamma$ is \myemph{surjective}.
Moreover, as $\gamma$ is a continuous map into a discrete group $\bZ$, it takes constant values on connected components of $\Stabilizer{\func,\partial\Cylinder}$, whence $\StabilizerId{\func,\partial\Cylinder}\subset\ker(\gamma)$.
Conversely, if $\dif\in\ker(\gamma)$, that is $\gamma(\dif)=\alpha_{\dif}|_{\Wman}=0$, i.e. $\alpha(\dif)=$ on $\partial\Cylinder$, we obtain from Lemma~\ref{lm:prop:StabIdfX}\ref{enum:lm:hamvf:shift:SIdfX} that $\dif = \flow_{\alpha_{\dif}}\in\StabilizerId{\func,\partial\Cylinder}$.
\qed

\medskip

{\bf Proof of~\ref{eqnu:th:cyl:func:incl}.}
Let $\func\in\FF(\Cylinder,\bR)$.
\jadded{33.3-33.4}{}
Fix two $\func$-regular neighborhoods $\Uman$ and $\Nman$ of $\Wman$ with $\Uman\subset \Int{\Nman}$, and let $H:\Stabilizer{\func, \Vman} \times\UInt \to \Stabilizer{\func, \Vman}$, $H_t(\dif) = \dif\circ(\flow_{t\beta(\dif)})^{-1}$, be the homotopy from Lemma~\ref{lm:extension_of_shift_functions}, where $\beta(\dif):\Mman\to\bR$ is a certain $C^{\infty}$ function such that 
\begin{align*}
\beta(\dif)|_{\Wman}&=\gamma(\dif), &
\beta(\dif)|_{\overline{\Mman\setminus\Nman}}&=0, &
\dif|_{\Uman} &= \flow_{\beta(\dif)}|_{\Uman}, &
H_t(\dif)|_{\Uman} = \flow_{ (1-t)\beta(\dif)\circ\flow^{-1}_{t\beta(\dif)}}.
\end{align*}
Statement~\ref{eqnu:th:cyl:func:incl} of Theorem~\ref{th:stab:annulus} is contained in~\ref{enum:lm:shift_map_on_W_cyl1:SptfUV} of the following Lemma~\ref{lm:shift_map_on_W_cyl1}.
\begin{lemma}\label{lm:shift_map_on_W_cyl1}
For $\dif\in\Stabilizer{\func,\Vman}$ the following statements hold.
\begin{enumerate}[label={\rm(\roman*)}, itemsep=1ex]
\item\label{enum:lm:shift_map_on_W_cyl1:lf_Ht_U}
$\lf{\dif} = \hflow_{\beta(\dif)\circ p}$ on $\cov{\Uman} := p^{-1}(\Uman)$;

\item\label{enum:lm:shift_map_on_W_cyl1:beta_H_t_h__U}
$\gamma( H_t(\dif) ) = (1-t)\gamma(\dif)\circ \flow_{t\gamma(\dif)}^{-1}$ \ and \ 
$\beta( H_t(\dif) )|_{\Uman} = (1-t)\beta(\dif)\circ\flow^{-1}_{t\beta(\dif)}$;

\item\label{enum:lm:shift_map_on_W_cyl1:gamma_h_SidfDC}
$\gamma(\dif) \equiv 0$ 
$\Longleftrightarrow$ 
$\lf{\dif}$ is fixed on $\bR\times1$
$\Longleftrightarrow$ 
$\dif\in\StabilizerIsotId{\func,\partial\Cylinder}\equiv\StabilizerIsotId{\func,\Vman}\cap\DiffId(\Cylinder,\partial\Cylinder)$;

\item\label{enum:lm:shift_map_on_W_cyl1:beta_h_SidfUV}
$\beta(\dif)|_{\Uman} \equiv 0$
$\ \stackrel{\ref{enum:lm:shift_map_on_W_cyl1:lf_Ht_U}}{\Longleftrightarrow} \ $
$\lf{\dif}$ is fixed on $\cov{\Uman}$
$\ \ \stackrel{\text{\rm Lemma~\ref{lm:DiffIdA}\ref{enum:lm:DiffIdA:DidC_charact}}}{\Longleftrightarrow}\ \ $
$\dif\in\StabilizerIsotId{\func,\Uman\cup\Vman}$;

\item\label{enum:lm:shift_map_on_W_cyl1:H_1_SfV}
$H_1(\Stabilizer{\func,\Vman}) \subset \StabilizerIsotId{\func,\Uman\cup\Vman} \subset \StabilizerIsotId{\func,\partial\Cylinder}$, that is $\beta(H_1(\dif))|_{\Uman}\equiv0$ for all $\dif\in\Stabilizer{\func,\Vman}$;
\jadded{33.3-33.4}{}

\item\label{enum:lm:shift_map_on_W_cyl1:Ht_SptfV}
$H\bigl(\StabilizerIsotId{\func,\Uman\cup\Vman}\times\UInt\bigr) \subset \StabilizerIsotId{\func,\Uman\cup\Vman}$
\ and \
$H\bigl(\StabilizerIsotId{\func,\partial\Cylinder}\times\UInt\bigr)\subset \StabilizerIsotId{\func,\partial\Cylinder}$, i.e. if $\beta(\dif)|_{\Uman}\equiv0$ (resp. $\gamma(\dif)\equiv0$), then $\beta(H_t(\dif))|_{\Uman}\equiv0$ (resp. $\gamma(H_t(\dif))\equiv0$) for all $t\in\UInt$;

\item\label{enum:lm:shift_map_on_W_cyl1:SptfUV}
\xadded{1}{0}{33.1}{}
$H$ is a deformation of $\bigl(\Stabilizer{\func,\Vman}, \FolStabilizer{\func,\Vman}\bigr)$ into $\bigl(\StabilizerIsotId{\func,\Uman\cup \Vman}, \FolStabilizerIsotId{\func,\Uman\cup \Vman}\bigr)$ as well as into $\bigl(\StabilizerIsotId{\func,\partial\Cylinder}, \FolStabilizerIsotId{\func,\partial\Cylinder}\bigr)$.
\end{enumerate}
\end{lemma}
\begin{proof}
\ref{enum:lm:shift_map_on_W_cyl1:lf_Ht_U}
Since $\dif|_{\Uman}=\flow_{\beta(\dif)}|_{\Uman}$, it follows from~\eqref{equ:lifting_of_shift} that $\hflow_{\beta(\dif)\circ p}|_{\cov{\Uman}}$ is a lifting of $\dif$.
On the other hand, $\lf{\dif}$ is also a lifting of $\dif$.
Moreover, 
$\lf{\dif}|_{\bR\times 1} = 
\hflow_{\gamma(\dif)\circ p}=
\hflow_{\beta(\dif)\circ p}|_{\bR\times 1}$,
whence those liftings coincide on $\cov{\Uman}$.

\ref{enum:lm:shift_map_on_W_cyl1:beta_H_t_h__U}
Notice that 
\begin{align*}
\gamma(H_t(\dif)) &\stackrel{\eqref{equ:gamma_h_fbinv}}{\ \ = \ \ }
\bigl(\gamma(\dif) - t\beta(\dif)|_{\Wman}\bigr)\circ\flow_{t\beta(\dif)}^{-1}|_{\Wman} =  \\
& \ \ = \ \ \bigl(\gamma(\dif) - t\gamma(\dif)\bigr)\circ\flow_{t\gamma(\dif)}^{-1}=
(1-t)\gamma(\dif)\circ\flow_{t\gamma(\dif)}^{-1}.
\end{align*}

Denote $\delta = (1-t)\beta(\dif)\circ\flow^{-1}_{t\beta(\dif)}$.
Then, $\flow_{\beta(H_t(\dif))}= H_t(\dif)|_{\Uman} = \flow_{\delta}$\,\!, that is $\beta(H_t(\dif))$ and $\delta$ are two shift functions for $H_t(\dif)$ on the connected set $\Uman$ having no singular points of $\flow$.
Moreover, $\beta(H_t(\dif))|_{\Wman} = \gamma(H_t(\dif)) = (1-t)\gamma(\dif)\circ\flow_{t\gamma(\dif)}^{-1} = \delta|_{\Wman}$.
Hence by Lemma~\ref{lm:shift}\ref{enum:shift:shift_func_uniqueness} they must coincide on all of $\Uman$.

Statement~\ref{enum:lm:shift_map_on_W_cyl1:gamma_h_SidfDC} is the same as~\eqref{equ:charact_gamma_0} for $\dif\in\Stabilizer{\func,\Vman}$, and the arguments for~\ref{enum:lm:shift_map_on_W_cyl1:beta_h_SidfUV} are presented in the formulation.
Statements~\ref{enum:lm:shift_map_on_W_cyl1:H_1_SfV} and~\ref{enum:lm:shift_map_on_W_cyl1:Ht_SptfV} directly follow from formulas for $\gamma$ and $\beta$ given in~\ref{enum:lm:shift_map_on_W_cyl1:beta_H_t_h__U} and characterizations of $\StabilizerIsotId{\func,\Uman\cup\Vman}$ and $\StabilizerIsotId{\func,\partial\Cylinder}$ in~\ref{enum:lm:shift_map_on_W_cyl1:gamma_h_SidfDC} and~\ref{enum:lm:shift_map_on_W_cyl1:beta_h_SidfUV}.

\ref{enum:lm:shift_map_on_W_cyl1:SptfUV}
By Lemma~\ref{lm:extension_of_shift_functions}\ref{extsh:Ht_ADfU}, $H_t\bigl(\FolStabilizer{\func} \cap \Stabilizer{\func,\Vman}) \subset \FolStabilizer{\func}$.
Together with~\ref{enum:lm:shift_map_on_W_cyl1:H_1_SfV} and~\ref{enum:lm:shift_map_on_W_cyl1:Ht_SptfV} this implies that the spaces 
$\FolStabilizer{\func,\Vman} := \FolStabilizer{\func} \cap \Stabilizer{\func,\Vman}$, 
$\FolStabilizerIsotId{\func,\Uman\cup\Vman}:=\FolStabilizer{\func} \cap \StabilizerIsotId{\func,\Uman\cup\Vman}$,
and $\FolStabilizerIsotId{\func,\partial\Cylinder}:=\FolStabilizer{\func} \cap \StabilizerIsotId{\func,\partial\Cylinder}$ are also invariant under $H$.
Hence $H$ is a deformation of the mentioned pairs.
\end{proof}

\medskip 

{\bf Proof of~\ref{eqnu:th:cyl:func:subset}.}
Recall that $\func\in\FF(\Mman,\Pman)$, where $\Mman = \Bman\cup\Cylinder$ is a not necessarily orientable connected compact surface, $\Bman$ and $\Cylinder$ are two $\func$-adapted subsurfaces such that $\Wman:=\Bman\cap\Cylinder$ is a regular component of some level set of $\func$ separating $\Mman$, and $\Xman\subset\Bman\setminus\Wman$ is an $\func$-adapted submanifold. 
Moreover, $\Cylinder$ is an annulus and $\dif(\Cylinder)=\Cylinder$ for all $\dif\in\Stabilizer{\func,\Xman\cup\Vman}$, where $\Vman$ is another boundary component of $\Cylinder$ distinct from $\Wman$. 
We have to construct a homotopy equivalence $\Stabilizer{\func, \Xman \cup \Vman} \cong \Stabilizer{\func|_{\Bman}, \Xman \cup \Wman} \times \StabilizerIsotId{\func|_{\Cylinder}, \partial\Cylinder}$.

Fix any $\func$-regular neighborhood $\Uman$ of $\Wman$, put $\UB := \Bman\cap\Uman$, $\UA := \Uman\cap\Cylinder$, and define the following groups:
\begin{align*}
\ggB &= \Stabilizer{\func, \Xman\cup\Uman\cup\Cylinder}, &
\ggA &= \Stabilizer{\func, \Bman\cup\Uman\cup\Vman}, &
\ggC &= \Stabilizer{\func, \Xman\cup\Uman\cup\Vman}, \\
\ggB' &= \StabilizerIsotId{\func, \Xman\cup\Uman\cup\Cylinder}, &
\ggA' &= \StabilizerIsotId{\func, \Bman\cup\Uman\cup\Vman}, &
\ggC' &= \StabilizerIsotId{\func, \Xman\cup\Uman\cup\Vman}.
\end{align*}
Then we have the following homeomorphisms ($\cong$) obtained by restrictions maps as in Lemma~\ref{lm:reduction_Mconn_V_dM}, and homotopy equivalences ($\simeq$) from Corollary~\ref{cor:SfU_SfX_hom_equiv}:
\begin{align*}
\ggB &\cong \Stabilizer{\func|_{\Bman}, \Xman\cup \UB} \simeq 
\Stabilizer{\func|_{\Bman}, \Xman\cup\Wman}, &
\ggA &\cong \Stabilizer{\func|_{\Cylinder}, \UA\cup\Vman} \simeq 
\Stabilizer{\func|_{\Cylinder}, \partial\Cylinder}.
\end{align*}
Evidently, $\ggA$ and $\ggB$ are subgroups of $\ggC$.
Moreover, $\supp(b)\subset \overline{\Bman\setminus\Uman}$ for each $b\in\ggB$ and $\supp(a)\subset \overline{\Cylinder\setminus\Uman}$ for each $a\in\ggA$.
In particular, $\ggA$ and $\ggB$ commute, $\ggA\cap\ggB = \{ \id_{\Mman}\}$ and they generate all of $\ggC$, i.e. the map $\alpha:\ggB\times\ggA\to\ggC$ given by $\alpha(b,a)=b\circ a$ is an isomorphism.

Let also $\Nman$ be an $\func$-regular neighborhoods of $\Uman$ with $\Uman\subset \Int{\Nman}$.
Since $\Nman\cup\Cylinder$ is an annulus, it is orientable (regardless of whether $\Mman$ is orientable or not), and so one can construct a vector field $\fld$ on $\Mman$ such that the restriction $\fld|_{\Nman\cup\Cylinder}$ is Hamiltonian like for $\func|_{\Nman\cup\Cylinder}$.
Let $\flow$ be the flow on $\Nman\cup\Cylinder$ generated by $\fld$.

Since $\Cylinder$ is invariant with respect to all diffeomorphisms from $\Stabilizer{\func,\Vman}$, we still have a well-defined map $\gamma:\Stabilizer{\func,\Vman}\to \Ci{\Wman}{\bR}$ as above satisfying~\eqref{equ:h_lfh}.
Then similarly to the previous case~\ref{eqnu:th:cyl:func:incl} we can define a homotopy  $H:\Stabilizer{\func, \Xman \cup \Vman} \times\UInt \to \Stabilizer{\func, \Xman \cup \Vman}$ such that $H_0$ is the identity and $H_1(\Stabilizer{\func, \Xman \cup \Vman}) \subset \Stabilizer{\func, \Xman\cup \Uman\cup\Vman}=:\ggC \cong \ggB\times\ggA$.

Moreover, it follows from Lemma~\ref{lm:shift_map_on_W_cyl1} that $\beta(\dif)|_{\Uman}\equiv0$ for each $\dif\in\ggB \cup \ggA'$.
Indeed, if $\dif\in\ggB$, then $\dif$ is fixed on $\Cylinder$, whence $\beta(\dif)\equiv0$ even on $\Cylinder$.
On the other hand, if $\dif\in\ggA$, then $\beta(\dif)|_{\Uman}\equiv0$ due to Lemma~\ref{lm:shift_map_on_W_cyl1}\ref{enum:lm:shift_map_on_W_cyl1:beta_h_SidfUV}.
Hence $H$ preserves the subgroup $\ggB \times \ggA'$.
Moreover, by Lemma~\ref{lm:shift_map_on_W_cyl1}\ref{enum:lm:shift_map_on_W_cyl1:Ht_SptfV}, $H_1(\dif)|_{\Cylinder} \in \StabilizerIsotId{\func|_{\Cylinder},\UA\cup\Vman}$, whence in fact $H_1(\Stabilizer{\func, \Xman \cup \Vman}) \subset \ggB\times\ggA'$.
Therefore $H_1$ is a homotopy equivalence 
\[
\Stabilizer{\func,\Xman\cup \Vman} \ \simeq \ \ggB\times\ggA' \ \simeq \ \Stabilizer{\func|_{\Bman}, \Xman\cup\Wman} \times \StabilizerIsotId{\func|_{\Cylinder}, \partial\Cylinder}.
\]

Notice also that by Lemma~\ref{lm:DiffMX_DiffMV}, $\ggB' = \ggB\cap \StabilizerIsotId{\func,\Xman\cup\Vman}$, whence $H$ also preserves $\ggB'$.
Indeed, the inclusion $\ggB' \subset \ggB\cap \StabilizerIsotId{\func,\Xman\cup\Vman}$ is evident.
Conversely, suppose $\dif\in\ggB\cap \StabilizerIsotId{\func,\Xman\cup\Vman}$, i.e. $\dif$ is fixed on the \myemph{annulus} $\Uman\cap\Cylinder$ being also a \myemph{regular neighborhood} of $\Vman$ and is isotopic to $\id_{\Mman}$ rel.~$\Vman$.
Then by Lemma~\ref{lm:DiffMX_DiffMV}, $\dif$ is isotopic to $\id_{\Mman}$ rel.~$\Uman\cap\Cylinder$, i.e. $\dif\in\ggB'$.

Thus $H_1(\StabilizerIsotId{\func,\Xman\cup\Vman}) \subset \StabilizerIsotId{\func,\Xman\cup\Vman} \cap (\ggB\times\ggA') \subset \ggB'\times\ggA'$, whence that $H_1$ induces a homotopy equivalence $\StabilizerIsotId{\func,\Xman\cup \Vman} \ \simeq \ \ggB'\times\ggA'$.
Finally, by Lemma~\ref{lm:extension_of_shift_functions}\ref{extsh:Ht_ADfU}, $H$ preserves the intersections of all the above groups with $\FolStabilizer{\func}$, which implies that $H$ induces the required homotopy equivalences of the corresponding pairs.

{\bf Proof of~\eqref{equ:bseq_two_cyl}.}
Suppose $\Xman=\varnothing$.
Then~\eqref{equ:cyl:split} can be written in the following form: 
\begin{equation}\label{equ:cyl_split_X_empty}
\seqStabIsotId{\func,\Vman} \cong 
\seqStabIsotId{\func|_{\Bman}, \Wman} \times 
\seqStabIsotId{\func|_{\Cylinder}, \partial\Cylinder}.
\end{equation}
Since $\Cylinder$ is a collar of the boundary component $\Vman$, it follows that $\Mman$ and $\Bman$ are diffeomorphic.
If $\Mman$ and $\Bman$ are $2$-disks, then $\Wman=\partial\Bman$, $\Vman=\partial\Mman$, and~\eqref{equ:cyl_split_X_empty} is the same as the required isomorphism~\eqref{equ:bseq_two_cyl}:
$\seqStabIsotId{\func,\partial\Mman} \cong 
\seqStabIsotId{\func|_{\Bman}, \partial\Bman} \times 
\seqStabIsotId{\func|_{\Cylinder}, \partial\Cylinder}$.
On the other hand, if $\Mman$ and $\Bman$ are annuli, then by the case~\ref{eqnu:th:cyl:func:incl}, $\seqStabIsotId{\func,\Vman} \cong \seqStabIsotId{\func,\partial\Mman}$ and
$\seqStabIsotId{\func|_{\Bman},\Wman} \cong \seqStabIsotId{\func|_{\Bman},\partial\Bman}$,
whence again~\eqref{equ:cyl_split_X_empty} implies~\eqref{equ:bseq_two_cyl}.
\qed

\section{Sketch of proof of Theorem~\ref{th:stab:disk:one_crpt}}
\label{sect:proof:th:disk:func_one_crpt}
{\bf Preliminaries.}
Let $\HPlane = \bR\times[0,+\infty)$ and $p:\HPlane\to\bC$ be the map defining polar coordinates, that is $p(\phi,r) = r e^{2\pi i \phi}$.
Then $p$ is a ``branched'' covering over $0\in\bC$ and the map $\xi:\HPlane\to\HPlane$, $\xi(\phi, r)=(\phi+1, r)$, generates the group of covering transformations and satisfies $p\circ \xi = p$.
Moreover, the restriction of $p:\bR\times(0,+\infty) \to \bC\setminus 0$ is the universal covering map.
The following statement seems to be well-known.
\begin{lemma}\label{lm:branched_lift}
Let $\Uman$ be an open neighborhood of $0$ in $\bC$, $\dif=(X,Y):\Uman\to\bC$ a $\Cr{k}$-embedding, $k\geq1$, such that $\dif(0)=0$, and $\hdif=(\Phi,R): p^{-1}(\Uman\setminus 0) \to \bR\times(0,+\infty)$ any continuous (and therefore $\Cr{k}$) lifting of $\dif$, i.e. $p\circ\hdif = \dif\circ p$.
Then $\hdif$ extends on $\bR\times 0$ to become a $\Cr{k-1}$-embedding $\hdif: p^{-1}(\Uman) \to \HPlane$ which also satisfies $p\circ\hdif = \dif\circ p$.

If the Jacobi matrix $A$ of $\dif$ at $0$ is a rotation, i.e. $A(z) = z e^{2\pi i a}$ for some $a\in\bR$, then there exists $k\in\bZ$ such that $\hdif(\phi,0) = (\phi + a + k, 0)$ for all $\phi\in\bR$.
\end{lemma}
\newcommand\mvect[2]{\left(\begin{matrix}#1 \\ #2\end{matrix}\right)}
\newcommand\smvect[2]{\left(\begin{smallmatrix}#1 \\ #2\end{smallmatrix}\right)}
\newcommand\mmatr[4]{\left(\begin{matrix}#1 & #2 \\ #3 & #4 \end{matrix}\right)}
\newcommand\smmatr[4]{\left(\begin{smallmatrix}#1 & #2 \\ #3 & #4 \end{smallmatrix}\right)}

\begin{proof}
We should prove that the coordinate functions $R$ and $\Phi$ of $\hdif$ can be defined on $\bR\times0$ and become of class $\Cr{k-1}$.
As $\dif(0)=0$, it follows from the Hadamard lemma that there exist $\Cr{k-1}$ functions $\alpha_{ij}:\Uman\to\bR$, $i,j\in\{0,1\}$, such that $X(x,y)=x \alpha_{00} + y \alpha_{01}$ and $Y(x,y)=x \alpha_{10} + y \alpha_{11}$.
Moreover, if we put $A = \smmatr{\alpha_{00}}{\alpha_{01}}{\alpha_{10}}{\alpha_{11}}$, then $A(0,0)$ is the Jacobi matrix of $\dif$ at the origin.
In particular, $A$ is non-degenerate at $(0,0)$, whence decreasing $\Uman$ if necessary, one can assume that $A$ is non-degenerate everywhere on $\Uman$.
Hence one can define the following $\Cr{k-1}$-map $q:\HPlane \to \bC\setminus0$ by 
$q(\phi,r) = 
\smmatr{\alpha_{00}(p(\phi, r))}
       {\alpha_{01}(p(\phi, r))}
       {\alpha_{10}(p(\phi, r))}
       {\alpha_{11}(p(\phi, r))}
\smvect{\cos\phi}{\sin\phi}$.

Fix $\cov{q}=(\cov{R},\cov{\Phi}):\HPlane \to \bR\times(0,+\infty)$ any lifting of $q$, i.e. $q = p\circ\cov{q}$, and thus $q(\phi,r) = \cov{R} e^{2\pi i\cov{\Phi}}$.
Then the relation $p\circ\hdif=\dif\circ p: p^{-1}(\Uman\setminus 0) \to\bR\times(0,+\infty)$ can be written as follows:
\begin{equation}\label{equ:pth=hp}
\smvect{R\cos(2\pi \Phi)}{R\sin(2\pi \Phi)} =  
r 
\smmatr{\alpha_{00}(p(\phi, r))}
{\alpha_{01}(p(\phi, r))}
{\alpha_{10}(p(\phi, r))}
{\alpha_{11}(p(\phi, r))}
 = r q(\phi,r).
\end{equation}
Hence $R e^{2\pi i\Phi} \stackrel{\eqref{equ:pth=hp}}{=} r q(\phi,r) = r\cov{R}e^{2\pi i\cov{\Phi}}$ which implies that $R=r\cov{R}$ and $\Phi = \cov{\Phi} + n$ on $p^{-1}(\Uman\setminus0)$ for some $n\in\bZ$.
The right parts of these formulas are $\Cr{k-1}$ on $p^{-1}(\Uman)$ and therefore they give the desired extensions of $R$ and $\Phi$.

The second statement is easy and we leave it for the reader.
\end{proof}

{\bf Proof of Theorem~\ref{th:stab:disk:one_crpt}.}
Not loosing generality and due to Axiom~\AxCrPt\ one can assume that $\func:\bR^2\to\bR$ is a homogeneous polynomial without multiple factors for which the origin $0$ is a unique critical point and this points is a global minimum.
\jadded{33.6-34.3}{}%
Then $\aDisk = \func^{-1}(\UInt)$ is a $2$-disk and $\func|_{\aDisk} \in \FF(\aDisk,\bR)$ has a unique critical point.
Therefore it suffices to prove Theorem~\ref{th:stab:disk:one_crpt} for $\func|_{\aDisk}$.
Denote
\begin{align*}
\cov{\aDisk} &= p^{-1}({\aDisk}),  & 
\dCyl &= {\aDisk}\setminus 0, &  
\hdCyl&=p^{-1}(\dCyl) = \cov{\aDisk}\setminus(\bR\times0),  & 
\partial\hdCyl &= p^{-1}(\partial{\aDisk}).
\end{align*}
Let also \jadded{34.7}{$\fld = -\frac{\partial\func}{\partial y}\, \tfrac{\partial}{\partial x}  +  \frac{\partial\func}{\partial x}\, \tfrac{\partial}{\partial y}$} be the Hamiltonial vector field of $\func$ and $\flow$ the corresponding flow on $\aDisk$.

Since $p:\hdCyl \to \dCyl$ is a covering map, each $\dif\in\Stabilizer{\func,\partial\aDisk}$ lifts to a unique $\Cinfty$ diffeomorphism $\lf{\dif}:\hdCyl\to\hdCyl$ fixed on $\partial\hdCyl$ and by Lemma~\ref{lm:branched_lift}, $\lf{\dif}$ extends to a $\Cinfty$ diffeomorphism of all $\cov{\aDisk}$.
In particular, $\flow$ lifts to a unique flow $\hflow$ on $\cov{\aDisk}$ such that $p \circ\hflow_t = \flow_t\circ p$ and $\xi\circ\hflow_t=\hflow_t\circ\xi$ for all $t\in\bR$.
As each orbit $\omega$ of $\flow$ in $\dCyl$ is closed and its inverse $p^{-1}(\omega)$ is an orbit of $\hflow$, we get from Lemma~\ref{lm:shifts:covering}\ref{enum:shift:covering:lift} that for each $\dif\in\Stabilizer{\func,\Vman}$ there exists a unique $\Cinfty$ function $\alpha_{\dif}:\dCyl\to\bR$ such that 
\begin{align*}
\alpha_{\dif}|_{\partial\dCyl}&=0, &
\dif|_{\dCyl} &= \flow_{\alpha_{\dif}},&
\lf{\dif}|_{\hdCyl} &= \hflow_{\alpha_{\dif}\circ p}.
\end{align*}

\jadded{35.5}{From this point the proof of non-degenerate and degenerate cases are essentially differs.}
Let $\LStab(\func) \subset \GL(\bR^2)$ be the group of Jacobi matrices of all $\dif\in\Stabilizer{\func,\Vman}$ at $0$.
\begin{itemize}[wide]
\item[\ref{enum:2disk:1crpt:nondeg}]
Suppose $0$ is a non-degenerate critical point.
Then by Lemma~\ref{lm:LStabf} one can also assume that $\func(x,y)=x^2+y^2$ and $\LStab(\func) = \SO(2)$.
In this case $\aDisk=\{|z|\leq1\}$, the Hamiltonian flow $\flow:\aDisk\times\bR\to\aDisk$ for $\func$ and its lifting $\hflow:\HPlane\times\bR\to\bR$ can be given by $\flow(z,t)=z e^{2\pi i t}$ and $\hflow(\phi,r,t)=(\phi+t,r)$.
Notice that all orbits of $\hflow$ are non-singular and non-closed.
Therefore $\alpha_{\dif}\circ p$ extends to a unique $\Cinfty$ function on all of $\hdCyl$.

On the other hand, the Jacobi matrix of $\dif$ is a rotation, whence by Lemma~\ref{lm:branched_lift} there exists $a\in\bR$ such that $\lf{\dif}(\phi,0)=(\phi+a,0) = \hflow_{a}(\phi,0)$.
This means that $\alpha_{\dif}\circ p(\phi,0) = a$ for all $\phi\in\bR$, whence if we set $\alpha_{\dif}(0)=a$, then $\alpha_{\dif}$ becomes continuous on $\aDisk$.

It was shown in~\cite[Lemma~31]{Maksymenko:TA:2003} that $\alpha_{\dif}$ is in fact $C^{\infty}$ on all of $\aDisk$, which implies that $\dif =\flow_{\alpha_{\dif}} \in\StabilizerId{\func,\Vman}$. 
Hence $\Stabilizer{\func,\Vman}\subset \StabilizerId{\func,\Vman}$ which means that $\pi_0\Stabilizer{\func,\Vman}=0$, and thus all groups in the sequence $\seqStabIsotId{\func,\Vman}$ are trivial.

\item[\ref{enum:2disk:1crpt:deg}]
Suppose $0$ is a degenerate extreme of $\func$ of symmetry index $m$.
Then $\deg\func\geq4$, and therefore the coordinate functions of the Hamiltonian vector field $\fld$ of $\func$ are polynomials of degree $\geq 3$.
This implies that the Jacobi matrix of each $\flow_t$ at $0$ is the identity, whence its lifting $\hflow_t$ must be fixed on $\bR\times 0$, i.e.\! $\bR\times 0$ consists of fixed orbits of $\hflow$.
Hence if $\alpha_{\dif}\circ p:\hdCyl\to\bR$ (even continuously) extends to $\bR\times0$, then $\lf{\dif}=\hflow_{\alpha_{\dif}\circ p}$ on all of $\cov{\aDisk}$, and therefore $\lf{\dif}$ must be fixed on $\bR\times 0$.
We will show that in fact the behavior of $\lf{\dif}$ on $\bR\times 0$ completely classifies path components of $\Stabilizer{\func,\partial\aDisk}$.

Indeed, by Lemma~\ref{lm:LStabf} one can assume that $\LStab(\func)$ consists of rotations by angles $2\pi k/m$, for $k\in\{0,\ldots,m-1\}$.
\jadded{35.2}{}%
In particular, the Jacobi matrix of $\dif$ is so.
Then again by Lemma~\ref{lm:branched_lift} there exists
$a_{\dif}\in\bR$ such that \jadded{35.3-35.4}{$\lf{\dif}(\phi,0)=(\phi+a_{\dif},0)$} for all $\phi\in\bR$, but now $a_{\dif}$ must be of the form $a_{\dif}=k/m + n$ for some $k\in\{0,1\ldots,m-1\}$ and $n\in\bZ$.
Hence we get a well defined map $\eta:\Stabilizer{\func,\partial\aDisk}\to \bZ$ given by $\eta(\dif)=a_{\dif}m$.

Evidently, $\eta$ is a continuous epimorphism into $\bZ$ with discrete topology.
Hence $\eta$ takes constant values on path components of $\Stabilizer{\func,\partial\aDisk}$.
In particular, $\eta(\StabilizerId{\func,\partial\aDisk}) = \eta(\id_{\aDisk})=0$, i.e. $\StabilizerId{\func,\partial\aDisk} \subset \ker(\eta)$.
The following technical statement implies that the inverse inclusion holds as well.

\begin{lemma}\label{lm:construct_shift_func_deg}
\begin{enumerate}[wide,label={\rm(\alph*)}]
\item\label{enum:resolve:jets}
Suppose the Jacobi matrix of $\dif\in\Stabilizer{\func,\partial\aDisk}$ at $0$ is the identity.
Then there is a $\Cinfty$-function $\gamma_1:\aDisk\to\bR$ supported in arbitrary small neighborhood of $0$ such that $\infty$-jets 
of $\dif$ and $\flow_{\gamma_1}$ at $0$ coincide, whence $\gdif = \flow_{\gamma_1}^{-1}\circ\dif$ is $\infty$-tangent to the identity at $0$, {\rm\cite[Theorem 7.1(1)]{Maksymenko:CEJM:2009}}.
 
\item\label{enum:resolve:flat}
Suppose $\gdif\in\Stabilizer{\func,\partial\aDisk}$ is $\infty$-tangent to the identity at $0$, and $\lf{\gdif}$ is fixed on $\bR\times 0$, i.e. $\eta(\gdif)=0$.
Then $\lf{\gdif}$ is $\infty$-tangent to the identity along $\bR\times 0$.
Moreover, put $\alpha_{\gdif}(0)=0$.
Then $\alpha_{\gdif}$ becomes $\Cinfty$ and \myemph{flat} at $0$.
In particular, $\gdif=\flow_{\alpha_{\gdif}}\in\StabilizerId{\func,\partial\aDisk}$, {\rm\cite[Proposition~3.4]{Maksymenko:hamv2}}.
\end{enumerate}
\end{lemma}
Let $\dif\in\ker(\eta)$.
Then the Jacobi matrix of $\dif$ at $0$ is the identity, whence by  Lemma~\ref{lm:construct_shift_func_deg}\ref{enum:resolve:jets} there exists a $\Cinfty$ function $\gamma$ such that $\gdif = \flow_{\gamma_1}^{-1}\circ\dif$ is $\infty$-tangent to the identity at $0$.
But $\eta(\flow_{\gamma_1})=0$, whence $\eta(\gdif)=-\eta(\flow_{\gamma_1})+\eta(\dif)=0$, and by Lemma~\ref{lm:construct_shift_func_deg}\ref{enum:resolve:flat}, $\gdif\in\StabilizerId{\func,\Vman}$.
Hence $\dif = \flow_{\gamma_1}\circ\gdif \in \StabilizerId{\func,\Vman}$ as well.
Thus $\ker(\eta) = \StabilizerId{\func,\Vman}$.

This implies that $\pi_0\Stabilizer{\func,\Vman} = 
\Stabilizer{\func,\Vman}/\StabilizerId{\func,\Vman}= \Stabilizer{\func,\Vman}/\ker(\eta)=\mathrm{image}(\eta)=\bZ$.

Moreover, notice that every level set of $\func$ if connected, whence each $\dif\in\Stabilizer{\func,\Vman}$ preserves each connected component of each level set of $\func$.
Therefore $\dif\in\FolStabilizer{\func,\Vman}$ iff its Jacobi matrix at $0$ is the identity, which in turn is equivalent to the assumption that $a_{\dif}\in\bZ$, that is $\eta(\dif)=a_{\dif}m$ is a multiple of $m$.
In other words, $\FolStabilizer{\func,\Vman} =\eta^{-1}(m\bZ)$.
This implies that $\seqStab{\func,\Vman} \cong \seqZ{m}:m\bZ\monoArrow\bZ\epiArrow\bZ_m$.
\end{itemize}

\section{Epimorphism $\eta:\Stabilizer{\func,\Vman}\to\bZ$}
Let $\Mman$ be an \myemph{orientable} surface, \jadded{37.1}{$\Vman$ a connected component of $\partial\Mman$,} $\widetilde{K}$ the union of all critical components of all level sets of $\func$, and $Z$ the connected component of $\Mman\setminus\widetilde{K}$ containing $\Vman$.
Then $\crLev = \overline{Z} \setminus Z$ is the ``closest'' to $\Vman$ critical component of $\func$,
similarly to the construction before Theorem~\ref{th:stab:disk_ann:gen_case}.
Let $c = \func(\crLev)$ and $\regNK$ be the connected component of $\func^{-1}[c-\eps,c+\eps]$ containing $\crLev$, where $\eps>0$ is so small that $\regNK\setminus\crLev$ contains no critical points and $\regNK\cap\partial\Mman$ is either empty or consists of some boundary components of $\partial\regNK$.
Then $\regNK$ is an $\func$-regular neighborhood of $\crLev$.
 
\begin{figure}[!htbp]
\centering
\includegraphics[height=2.2cm]{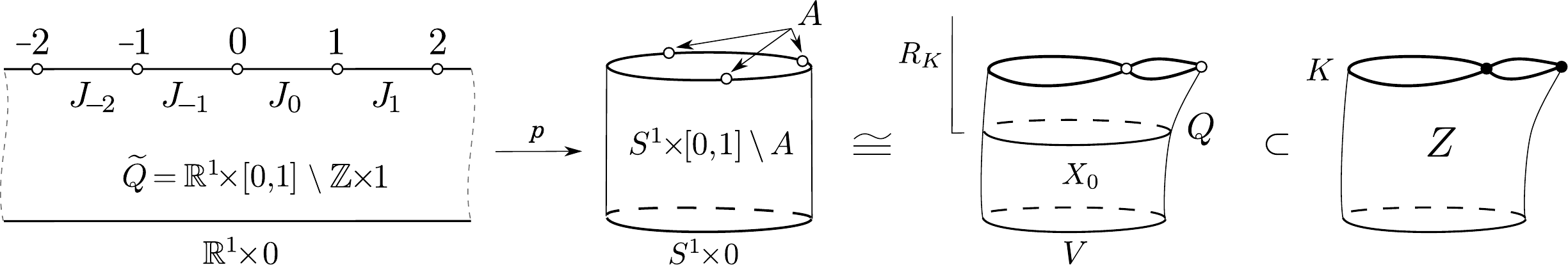}
\caption{}
\protect\label{fig:univ_cov_dCyl}
\end{figure}

Notice that $\dCyl = \overline{Z} \setminus \fSing$ is diffeomorphic with $(\Circle\times\UInt)\setminus\FSet$, where $\FSet$ is a finite subset of $\Circle\times 1$, see Figure~\ref{fig:univ_cov_dCyl}.
Denote by $n$ the number of intervals in $(\Circle\times1) \setminus \FSet$.
Evidently, they correspond to the edges of $\crLev$.

Let also $\hdCyl=(\bR\times\UInt) \setminus (\bZ\times1)$, $p:\hdCyl\to\dCyl$ be the universal covering map for $\dCyl$, and $J_i = (i,i+1) \times 1$ for $i\in\bZ$.
Thus for every edge $e$ of $\dCyl\cap\crLev$ there exists a unique $k\in\{0,\ldots,n-1\}$ such that $p^{-1}(e)=\{ J_{k+an} \}_{a\in\bZ}$.
Let also $\xi:\hdCyl\to\hdCyl$ be the diffeomorphism generating the group $\bZ$ of covering transformations such that $\xi(J_i) = J_{i+n}$, $i\in\bZ$.

Fix a Hamiltonian like flow $\flow$ for $\func$.
Then $\dCyl$ is invariant with respect to $\flow$, whence $\flow$ lifts to a flow $\hflow:\hdCyl\times\bR\to\hdCyl$ on $\dCyl$ commuting with $\xi$.
One can assume that period of the \myemph{closed orbit} $\Vman$ of $\flow$ equals $1$, and 
\begin{align}
\hflow(x,0,t) &= (x-nt,0), &
\xi(x,0) = (x+n, 0) = \hflow(x,0,-1).
\end{align}
In particular, $\flow$ and $\xi$ move points of $\bR\times0$ in opposite directions.

Let $\XFixA = \overline{\dCyl\setminus\regNK}$.
Take any $C^{\infty}$ function $\mu:\Mman\to\UInt$ such that \added{$\mu=0$} on $\Mman\setminus\XFixA$, \jadded{37.2-37.3}{$\mu=1$} near $\Vman$, and $\mu$ is constant along regular components of level sets of $\func$.
Then the map $\tau := \flow_{\mu}$ is an $\func$-adapted Dehn twist around $\Vman$ supported in $\XFixA$, and $\cov{\tau} = \hflow_{\mu\circ p}:\hdCyl\to\hdCyl$ is its lifting fixed on $\bR\times 0$.
Notice that $\cov{\tau} = \xi$ near $(\bR\setminus\bZ)\times 1$.
In particular, 
\begin{equation}\label{equ:htau_Jk_Jkn}
\cov{\tau}(J_k)=J_{k+n}, \qquad k\in\bZ.
\end{equation}

\begin{lemma}\label{lm:eta}
Suppose $\crLev$ is not a local extreme of $\func$ and let $\RB := \regNK \cup \XFixA$.
Then there exists an \myemph{epimorphism} $\eta:\Stabilizer{\func,\Vman} \to \bZ$ with the following properties.
\begin{enumerate}[wide, label={\rm(\arabic*)}, parsep=0.6ex]
\item\label{enum:lm:eta:diag}
We have the following commutative diagram of inclusions of subgroups:
\begin{equation}\label{equ:eta_diagram}
\aligned
\xymatrix@C=1em@R=1.2em{
 \StabilizerId{\func,\Vman} \stackrel{\ref{enum:lm:eta:diag:SidfV}}{=} \StabilizerId{\func}\cap\ker(\eta)     \ar@{^(->}[r]  \ar@{^(->}[d] &
 \FolStabilizer{\func,\Vman;\RB}          \ar@{^(->}[r]  \ar@{^(->}[d] &
 \Stabilizer{\func,\Vman;\RB} \stackrel{\ref{enum:lm:eta:diag:SfVR}}{=} \ker(\eta)  \ar@{->>}[r]  \ar@{^(->}[d] &
 0    \ar@{^(->}[d]  \\
\ \StabilizerId{\func}\cap\eta^{-1}(m\bZ) \   \ar@{^(->}[r]   & 
 \FolStabilizer{\func,\Vman}                               \ar@{^(->}[r]                &
 \Stabilizer{\func,\Vman;\regNK} \stackrel{\ref{enum:lm:eta:diag:inv_mZ}}{=} \eta^{-1}(m\bZ)           \ar@{->>}[r]  \ar@{^(->}[d] &
 m\bZ \ \ar@{^(->}[d]  \\ 
    &&
 \Stabilizer{\func,\Vman}                           \ar@{->>}[r]^-{\eta} &
 \bZ 
}
\endaligned
\end{equation}
where $m=\eta(\tau)$ and the image under $\eta$ of each group in the second row is $m\bZ$.
In fact, diagram~\eqref{equ:eta_diagram} is implied by the following statements.
\begin{enumerate}[label={\rm(\alph*)}, parsep=0.6ex]
\item\label{enum:lm:eta:diag:inv_mZ}
$\Stabilizer{\func,\Vman;\regNK} = \Stabilizer{\func,\Vman} \cap \Stabilizer{\func;\RB}=\eta^{-1}(m\bZ)$.
In other words, for $\dif\in\Stabilizer{\func,\Vman}$ the following conditions are equivalent:
\begin{enumerate}
\item[$\ccInvm$]
$\dif\in\Stabilizer{\func,\Vman;\regNK}$, i.e.\! there is a $\Cinfty$ function $\delta'_{\dif}:\regNK\to\bR$ such that $\dif|_{\regNK} = \flow_{\delta'_{\dif}}$;
\item[$\ccSfR$] 
$\dif\in\Stabilizer{\func,\Vman} \cap \Stabilizer{\func;\RB}$, i.e.\! there is a $\Cinfty$ function $\delta_{\dif}:\RB\to\bR$ such that $\dif|_{\RB} = \flow_{\delta_{\dif}}$;
\item[$\ccSfX$] 
$\dif\in\eta^{-1}(m\bZ)$, i.e. $\eta(\dif)=ma$ for some $a\in\bZ$.
\end{enumerate}
In these cases the functions $\delta'_{\dif}$ and $\delta_{\dif}$ are unique, $\delta_{\dif}|_{\regNK}=\delta'_{\dif}$, and $\delta_{\dif}(\Vman) = \eta(\dif)/m = a$.
 
\item\label{enum:lm:eta:diag:SfVR}
$\Stabilizer{\func,\Vman;\RB} = \ker(\eta)$;

\item\label{enum:lm:eta:diag:DeltafVR}
$\FolStabilizer{\func,\Vman;\RB} :=
\FolStabilizer{\func,\Vman} \cap \Stabilizer{\func,\Vman;\RB}
\stackrel{\ref{enum:lm:eta:diag:SfVR}}{\equiv} \FolStabilizer{\func,\Vman} \cap \ker(\eta)$;

\item\label{enum:lm:eta:diag:SidfV}
$\StabilizerId{\func,\Vman} = \StabilizerId{\func} \cap \ker(\eta)$.
\end{enumerate}

\item\label{enum:lm:eta:diag_pi0}
Every group in~\eqref{equ:eta_diagram} consists of path component of $\Stabilizer{\func,\Vman}$, whence that diagram yields the diagram consisting of the corresponding $\pi_0$-groups:
\begin{equation}\label{equ:eta_diagram_pi0}
\aligned
\xymatrix@C=1.6em@R=1.2em{
\ \pi_0\StabilizerId{\func,\Vman}=0              \ \ar@{^(->}[r]  \ar@{^(->}[d] &
\ \pi_0\FolStabilizer{\func,\RB}         \ \ar@{^(->}[r]  \ar@{^(->}[d] &
\ \pi_0\Stabilizer{\func,\RB} \ \ar@{->>}[r]  \ar@{^(->}[d] &
\ 0  \  \ar@{^(->}[d]  \\
\ \pi_0\bigl(\StabilizerId{\func}\cap\eta^{-1}(m\bZ)\bigr)  \   \ar@{^(->}[r]   & 
\ \pi_0\FolStabilizer{\func,\Vman}                              \ \ar@{^(->}[r]                &
\ \pi_0\Stabilizer{\func,\Vman\cup\regNK}        \ \ar@{->>}[r]  \ar@{^(->}[d] &
\ m\bZ \ \ar@{^(->}[d]  \\ 
 &&
\ \pi_0\Stabilizer{\func,\Vman}                           \ \ar@{->>}[r]^-{\eta} &
\ \bZ \
}
\endaligned
\end{equation}

\item\label{enum:lm:eta:3x3}
Let $\XFixA, \Zman_1,\ldots,\Zman_{\cnt}$ be all the connected components of $\overline{\Mman\setminus\regNK}$ and $\hZman_i = \Zman_i\cap\regNK$, $i=1,\ldots,\cnt$.
Then \eqref{equ:eta_diagram_pi0} implies one more commutative diagram:
\begin{equation}\label{equ:3x3_pi0}
\aligned
\xymatrix@R=3ex{
& 
\prod\limits_{i=1}^{\cnt} \pi_0\FolStabilizer{\func|_{\Zman_i},\hZman_i} \ar@{=}[d]_-{\cong} \ar@{^(->}[r] &
\prod\limits_{i=1}^{\cnt} \pi_0\Stabilizer{\func|_{\Zman_i},\hZman_i}    \ar@{=}[d]^-{\alpha^{-1}}_-{\cong} \ar@{->>}[r]  &
\prod\limits_{i=1}^{\cnt} \GrpKR{\func|_{\Zman_i},\hZman_i}   \ar@{=}[d]_-{\cong}     \\
 \seqStab{\func,\Xman} : \ar@{^(->}[d] &
\pi_0\FolStabilizer{\func,\Xman} \ar@{^(->}[r]                \ar@{^(->}[d] &
\pi_0\Stabilizer{\func,\Xman}    \ar@{->>}[r]                 \ar@{^(->}[d] &
\GrpKR{\func,\Xman}                                           \ar@{^(->}[d] \\
 \seqStab{\func,\Vman}:  \ar@{->>}[d] &
\pi_0\FolStabilizer{\func,\Vman} \ar@{^(->}[r] \ar@{->>}[d] &
\pi_0\Stabilizer{\func,\Vman}    \ar@{->>}[r]^-{\fc}  \ar@{->>}[d]^-{\hb} &
\GrpKR{\func,\Vman}                            \ar@{->>}[d] \\
 \seqZ{m}: &
m\bZ \ar@{^(->}[r] &
\bZ  \ar@{->>}[r] &
\bZ_m
}
\endaligned
\end{equation}
which can be regarded as a sequence homomorphisms of short exact sequences of its rows:  
\begin{equation}\label{equ:3x3_pi0:seq}
\prod\limits_{i=1}^{\cnt}\seqStab{\func|_{\Zman_i},\hZman_i} 
 \ \monoArrow \ \seqStab{\func,\Vman} \ \epiArrow \ \seqZ{m},
\end{equation}
where $\alpha$ is constructed in Lemma~\ref{lm:reduction_Mconn_V_dM}.
 
\item\label{enum:lm:eta:m=1}
If $m=1$, then~\eqref{equ:3x3_pi0:seq} ``\myemph{splits}'' in the sense that we have an isomorphism
\begin{equation}\label{equ:split_bseq}
\seqStab{\func,\Vman}
\ \cong \
\seqStab{\func,\Xman} \times \seqZ{1} 
\ \cong \ \prod\limits_{i=1}^{\cnt}\seqStab{\func|_{\Zman_i},\hZman_i}   \ \times \ \seqZ{1}.
\end{equation}
\end{enumerate}
\end{lemma}
\begin{proof}
\textit{Construction of $\eta$.}
Notice that each $\dif\in\Stabilizer{\func,\Vman}$ lifts to a unique diffeomorphism $\hdif:\hdCyl\to\hdCyl$ fixed on $\bR\times 0$.
Moreover, $\hdif$ preserves the order of intervals $\{J_k\}_{k\in\bZ}$, whence there exists a unique number $\tilde{\eta}(\dif)\in\bZ$ such that $\hdif(J_k) = J_{k+\tilde{\eta}(\dif)}$ for all $k\in\bZ$.
One easily checks that the correspondence $\dif\mapsto\tilde{\eta}(\dif)$ is a continuous homomorphism $\tilde{\eta}:\Stabilizer{\func,\Vman}\to\bZ$.

Since $\cov{\tau} =\xi$ near $\hdCyl\cap(\bR\times 1)$, we obtain that $\tilde{\eta}(\tau) = n>0$, whence $\tilde{\eta}$ is a non-trivial homomorphism, and thus $\tilde{\eta}(\Stabilizer{\func,\Vman}) = l\bZ$ for some $l\in\bN$.
Hence one can define an \myemph{epimorphism} $\eta:\Stabilizer{\func,\Vman} \to \bZ$ by $\eta = \tilde{\eta} / l$.
Then $m:=\eta(\tau) = n/l$.

\ref{enum:lm:eta:diag}
We need to prove \ref{enum:lm:eta:diag:inv_mZ}-\ref{enum:lm:eta:diag:SidfV} which will imply existence of diagram~\eqref{equ:eta_diagram}.
\xadded{1}{0}{38.2}{}

\begin{itemize}[wide]
\item[\ref{enum:lm:eta:diag:inv_mZ}]
In fact, for $\dif\in\Stabilizer{\func,\Vman}$ and $a\in\bZ$ even the following conditions are equivalent:
\begin{align*}
\ccInvm   &~\eta(\dif) = ma, \ \text{i.e.} \ \dif\in\eta^{-1}(m\bZ); &
\ccTildeEta &~\tilde{\eta}(\dif) = na = lma;   \\
\ccSfR     &~\dif\in\Stabilizer{\func,\Vman;\regNK}; &
\ccIntervals &~\hdif(J_i) = J_{i+na} = \xi^a(J_i) \ \text{for all} \ i\in\bZ; \\
\ccSfX      &~\dif\in\Stabilizer{\func;\RB};  &
\ccEdges &~\text{each edge of $\crLev$ is \invhp{\dif}}.
\end{align*}
Indeed, the implications $\ccInvm$  $\Leftrightarrow$ $\ccTildeEta$ $\Leftrightarrow$ $\ccIntervals$ $\Leftrightarrow$ $\ccEdges$ and $\ccSfX$ $\Rightarrow$ $\ccSfR$  are direct consequences of the definitions.
The equivalence $\ccEdges$ $\Leftrightarrow$  $\ccSfR$  follows from the equivalence of conditions \ref{enum:pres:edge}$\Leftrightarrow$\ref{enum:pres:has_shift_func} of Lemma~\ref{lm:shift_func_near_crlevel}.

$\ccSfR$ $\Rightarrow$ $\ccSfX$
Suppose $\dif\in\Stabilizer{\func,\Vman;\regNK}$, i.e. there exists a $\Cinfty$ function $\delta'_{\dif}:\regNK\to\bR$ such that $\dif|_{\regNK}=\flow_{\delta'_{\dif}}$.
Then by Lemma~\ref{lm:shifts:covering}\ref{enum:shift:covering:lift}, the restriction of $\delta'_{\dif}$ to the curve $\Wman = \regNK \cap \XFixA$ extends to $\XFixA$ to a continuous function $\delta_{\dif}: \RB = \regNK \cup \XFixA \to \bR$ such that $\dif|_{\RB} = \flow_{\delta_{\dif}}$ and $\delta_{\dif}$ is $\Cinfty$ on $\regNK$ and $\XFixA$.
Then by Lemma~\ref{lm:shift}\ref{enum:shift:local_existence}, $\delta_{\dif}$ is $\Cinfty$ near $\Wman$, and therefore it is $\Cinfty$ on $\RB$.
This means that $\dif\in\Stabilizer{\func;\RB}$ and proves $\ccSfX$.

Now let us show that $\delta_{\dif}(\Vman)=\eta(\dif)/m = a$.
By~\eqref{equ:lifting_of_shift} the map $\qdif = \hflow_{\delta_{\dif}\circ p}$ is a lifting of $\dif|_{\dCyl}$ preserving each orbit of $\hflow$.
Then every other lifting of $\dif$ has the form $\xi^{k}\circ\qdif$, $k\in\bZ$.
Since $\xi$ shifts intervals $J_i$, it follows that $\qdif$ is a unique lifting of $\dif$ preserving orbits of $\hflow$.
On the other hand, due to $\ccIntervals$, $\xi^{-a}\circ\lf{\dif}$ preserves each interval $J_k$ and therefore preserves each orbit of $\hflow$ as well.
Hence $\xi^{-a}\circ\lf{\dif} = \qdif$, 
Therefore for each $x\in\bR$ we have that
\[ \qdif(x,0) = \hflow\bigl(x,0; \delta_{\dif}(p(x,0)) \bigr) = \xi^{-a}\circ\lf{\dif}(x,0) = \xi^{-a}(x,0)=(x-a,0)=\hflow(x,0; a).\]
Hence $\delta_{\dif}\circ p$ takes the same value $a$ on all of $\bR\times0$, and thus $\delta_{\dif}(\Vman)=\delta_{\dif}\circ p(\bR\times 0)=a$.

Finally, uniqueness of $\delta_{\dif}$ and $\delta'_{\dif}$ follows from the assumption that $\regNK \subset \RB$ are connected and contain non-closed orbits of $\flow$, see Lemma~\ref{lm:prop:StabIdfX}\ref{enum:lm:hamvf:shift:saddles}.

\xadded{1}{0}{38.1}{}
\item[\ref{enum:lm:eta:diag:SfVR}]
We need to prove that $\Stabilizer{\func,\Vman;\RB} = \ker(\eta)$.
Notice that both these sets are contained in $\Stabilizer{\func,\Vman} \cap \Stabilizer{\func,\Xman} \stackrel{\ref{enum:lm:eta:diag:inv_mZ}}{=} \eta^{-1}(m\bZ)$.
Then $\dif\in\Stabilizer{\func,\Vman;\RB}$ iff $\delta_{\dif}(\Vman)= 0 = \eta(\dif)/m$, i.e. when $\dif\in\ker(\eta)$.

\item[\ref{enum:lm:eta:diag:DeltafVR}]
This statement is trivial.

\item[\ref{enum:lm:eta:diag:SidfV}]
We should prove that $\StabilizerId{\func,\Vman} = \StabilizerId{\func} \cap \ker(\eta)$.
The inclusion $\StabilizerId{\func,\Vman} \subset  \StabilizerId{\func}$ is evident.
Moreover, as $\eta$ is a continuous map into the discrete space $\bZ$, it take constant values of connected components of $\Stabilizer{\func,\Vman}$, whence $\eta(\StabilizerId{\func,\Vman})= \eta(\id_{\Mman}) = 0$.
Thus $\StabilizerId{\func,\Vman} \subset \StabilizerId{\func} \cap \ker(\eta)$.
Conversely, let $\dif\in\StabilizerId{\func} \cap \ker(\eta)$.
Then $\dif=\flow_{\alpha}$ for some $\Cinfty$ function $\alpha:\Mman\to\bR$.
Moreover, by~\ref{enum:lm:eta:diag:inv_mZ}, $\alpha|_{\RB} = \delta_{\dif}$, whence $\alpha(\Vman)=\delta_{\dif(\Vman)} = m\eta(\dif) = 0$.
Hence $\dif\in\StabilizerId{\func,\Vman}$ due to Lemma~\ref{lm:prop:StabIdfX}\ref{enum:lm:hamvf:shift:SIdfX}.
This proves~\ref{enum:lm:eta:diag:SidfV}.
\end{itemize}

It remains to note that $\tau =\flow_{\mu} \in \StabilizerId{\func}$ and $\eta(\tau)=m$, i.e. $\tau\in\StabilizerId{\func}\cap\eta^{-1}(m\bZ)$, whence the image of each group in the second row of~\eqref{equ:eta_diagram} is $m\bZ$.

\ref{enum:lm:eta:diag_pi0}
Since all the groups in~\eqref{equ:eta_diagram} contain the identity path component $\StabilizerId{\func,\Vman}$, it follows that they consist of path component of $\Stabilizer{\func,\Vman}$.
Therefore one can put $\pi_0$ at each of the groups in~\eqref{equ:eta_diagram} and the obtained diagram will still consist of inclusions.

Moreover, due to Corollary~\ref{cor:SX_Sidf}, we have that
$\pi_0\FolStabilizer{\func,\Vman;\RB}\cong\pi_0\FolStabilizer{\func,\Vman\cup\RB} \equiv \pi_0\FolStabilizer{\func,\RB}$,
$\pi_0\Stabilizer{\func,\Vman;\RB}\cong\pi_0\Stabilizer{\func,\Vman\cup\RB} =\pi_0\Stabilizer{\func,\RB}$, and 
$\pi_0\FolStabilizer{\func,\Vman;\regNK}\cong\pi_0\FolStabilizer{\func,\Vman\cup\regNK}$, which give the diagram~\eqref{equ:eta_diagram_pi0}. 

Statement~\ref{enum:lm:eta:3x3} is evident.

\ref{enum:lm:eta:m=1}
Let $\tau\in\FolStabilizer{\func,\Vman}$ be an $\func$-adapted Dehn twist supported in $\Int{\XFixA} \subset \RB$ and therefore commuting with each $\dif\in \Stabilizer{\func,\Xman}$.
Suppose $\eta(\tau)=m=1$.
We claim that then the mapping $\beta: \pi_0\Stabilizer{\func,\Xman} \times \bZ\to \pi_0\Stabilizer{\func,\Vman}$ by $\beta([\dif], l) = [\dif\circ  \tau^{l}]$ is an isomorphism of groups.
Indeed, 
\begin{align*}
\beta\bigl( ([\gdif], a) \bigr) \beta\bigl( ([\dif],b) \bigr) &= 
[\gdif \circ \tau^a]\, [\dif \circ \tau^b] = 
[\gdif \circ \tau^a \circ  \dif\circ \tau^b] \\
&=
[\gdif \circ  \dif \circ \tau^a \circ \tau^b] = 
\beta\bigl( ([\gdif\circ\dif], a+b) \bigr).
\end{align*}
Hence $\beta$ is a homomorphism.
Moreover, since $\pi_0\Stabilizer{\func,\Xman}=\ker(\eta)$ and $\eta(\tau)=m=1$, it easily follows that $\beta$ is a bijection, and therefore an isomorphism.

Finally, as $\tau\in\FolStabilizer{\func,\Vman}$, we have that $\dif\circ \tau^{l}\in\FolStabilizer{\func,\Vman}$ for all $\dif\in\FolStabilizer{\func,\Vman}$, whence  $\beta\bigl(\pi_0\FolStabilizer{\func,\Xman} \times \bZ\bigr) = \pi_0\FolStabilizer{\func,\Vman}$, which implies an isomorphism of sequences~\eqref{equ:split_bseq}.
\end{proof}

\section{Proof of Theorem~\ref{th:stab:disk_ann:gen_case}}
Suppose $(\Mman,\Vman)$ is either $(\Disk,\partial\Disk)$ or $(\Circle\times\UInt,\Circle\times 0)$ and $\func\in\FF(\Mman,\Pman)$. 
Notice that notations of Theorem~\ref{th:stab:disk_ann:gen_case} and Lemma~\ref{lm:eta} agree.
We will now deduce several consequences from that latter lemma for our special case of $(\Mman,\Vman)$.

\begin{lemma}\label{lm:MV_DdD_AdA}
\begin{enumerate}[wide, label={\rm(\arabic*)}]
\item\label{enum:cor:MV_DdD_AdA:SZ}
$\Stabilizer{\func,\Vman;\regNK} =\Stabilizer{\bZman}$, whence $\frac{\Stabilizer{\func,\Vman}}{\Stabilizer{\bZman}} \cong
\frac{\Stabilizer{\func,\Vman}}{\Stabilizer{\func,\Vman;\regNK}} \cong \bZ_m$ due to~\eqref{equ:eta_diagram}.

\item\label{enum:cor:MV_DdD_AdA:eff_action}
If $m\geq2$, then the action of $\frac{\Stabilizer{\func,\Vman}}{\Stabilizer{\bZman}}\cong \bZ_m$ on $\bZman$ is \myemph{semifree} and has precisely either one or two fixed elements.

\item\label{enum:cor:MV_DdD_AdA:gm_tau}
There exists $\gdif\in\Stabilizer{\func,\Vman}$ such that $\eta(\gdif) = 1$ and $\gdif^{m}\in\StabilizerId{\func}$.
Moreover, for every such $\gdif$ we have that $[\gdif^{m}] = [\tau]$ in $\pi_0\Stabilizer{\func,\Vman}$, and therefore $[\gdif^m]$ commutes with the subgroup $\pi_0\Stabilizer{\func,\RB}$ of $\pi_0\Stabilizer{\func,\Vman}$.
\end{enumerate}
\end{lemma}
\begin{proof}
\ref{enum:cor:MV_DdD_AdA:SZ}
\jadded{38.3}{}
Recall that $\hZman_i := \Zman_i \cap\regNK = \partial\Zman_i \cap\partial\regNK$ is a common boundary component of $\Zman_i$ and $\regNK$, and it splits $\Mman$.
It suffices to show that for $\dif\in\Stabilizer{\func,\Vman}$ the following conditions are equivalent:
\begin{enumerate}[label={\rm(\alph*)}]
\item\label{enum:cor:MV_DdD_AdA:SZ:SfVR}
$\dif\in\Stabilizer{\func,\Vman;\regNK}$;
\item\label{enum:cor:MV_DdD_AdA:SZ:hbdRk}
$\dif$ preserves boundary component of $\regNK$ with its orientation;
\item\label{enum:cor:MV_DdD_AdA:SZ:hdZ_or}
$\dif$ preserves each $\hZman_i$ with its orientation;
\item\label{enum:cor:MV_DdD_AdA:SZ:hSdZ}
$\dif(\hZman_i)=\hZman_i$ for all $i=1,\ldots,\cnt$;
\item\label{enum:cor:MV_DdD_AdA:SZ:hSZ}
$\dif\in\Stabilizer{\bZman}$, i.e. $\dif(\Zman_i)=\Zman_i$ for all $i=1,\ldots,\cnt$.
\end{enumerate}
The equivalence \ref{enum:cor:MV_DdD_AdA:SZ:SfVR}$\Leftrightarrow$\ref{enum:cor:MV_DdD_AdA:SZ:hbdRk} follows from the equivalence of conditions \ref{enum:pres:has_shift_func} and \ref{enum:pres:bnd_comp} of Lemma~\ref{lm:shift_func_near_crlevel}.
The implications \ref{enum:cor:MV_DdD_AdA:SZ:hbdRk}$\Rightarrow$\ref{enum:cor:MV_DdD_AdA:SZ:hdZ_or}$\Rightarrow$\ref{enum:cor:MV_DdD_AdA:SZ:hSdZ} are trivial.

\ref{enum:cor:MV_DdD_AdA:SZ:hSdZ}$\Rightarrow$\ref{enum:cor:MV_DdD_AdA:SZ:hdZ_or}.
Suppose $\dif(\hZman_i)=\hZman_i$.
By assumption $\dif$ is fixed on $\Vman$, whence it preserves orientation of $\Mman$.
Hence $\regNK$ and its boundary component $\hZman_i$ are \invhp{\dif}.

\ref{enum:cor:MV_DdD_AdA:SZ:hdZ_or}$\Rightarrow$\ref{enum:cor:MV_DdD_AdA:SZ:hbdRk}
By the construction $\regNK \cap \partial\Mman \subset \partial\regNK$.
If $\regNK \cap \partial\Mman=\varnothing$ (which always hold when $\Mman=\Disk$), then 
$\partial\regNK = \cup_{i=1}^{\cnt}\hZman_i$ and the required implication is trivial.
On the other hand, if $\regNK \cap \partial\Mman\not=\varnothing$, then $\Mman=\Circle\times\UInt$, $\regNK \cap \partial\Mman=\Circle\times 1$, and $\partial\regNK = (\Circle\times 1)\cup (\cup_{i=1}^{\cnt}\hZman_i)$.
Then each $\hZman_i$ is \invhp{\dif} by assumption~\ref{enum:cor:MV_DdD_AdA:SZ:hdZ_or}, while $\Circle\times 1$ is \invhp{\dif} since $\dif$ is fixed on $\Vman$.

\ref{enum:cor:MV_DdD_AdA:SZ:hSdZ}$\Leftrightarrow$\ref{enum:cor:MV_DdD_AdA:SZ:hSZ}
Suppose $\dif(\hZman_i)=\hZman_i$ for some $i=1,\ldots,\cnt$.
Let also $\Zman_j=\dif(\Zman_i)$ for some $j$.
Then $\varnothing \not= \dif(\Zman_i \cap\regNK) \cap \Zman_i \supset \dif(\Zman_i)\cap \Zman_i = \Zman_j\cap \Zman_i$, which is possible only if $i=j$.

Conversely, if $\dif(\Zman_i)=\Zman_i$, then $\dif(\hZman_i) = \dif(\Zman_i \cap\regNK) = \dif(\Zman_i)\cap\dif(\regNK)=\Zman_i \cap\regNK=\hZman_i$.

\medskip 

\ref{enum:cor:MV_DdD_AdA:eff_action}
Suppose that $m\geq 2$.
It suffices to show that \myemph{all non-trivial elements of $\frac{\Stabilizer{\func,\Vman}}{\Stabilizer{\bZman}}\cong \bZ_m$ have the same set of fixed elements in $\bZman$ consisting of either one or two elements.}

Pinch every boundary component $\Wman$ of $\Mman$ into a point $z_{\Wman}$ and denote the obtained surface (being a $2$-sphere) by $\hMman$.
Notice that there is a CW-partition $\Theta$ of $\hMman$ whose $0$-cells are critical points of $\func$ belonging to $\crLev$, $1$-cells are connected components of $\crLev\setminus\fSing$, and the $2$-cells of $\Theta$ are connected components of $\hMman\setminus\crLev$.
These $2$-cells are in 1-1 correspondence with elements of $\bZman$.

\newcommand\sgmh{\sigma}

Let $\dif\in\Stabilizer{\func,\Vman}$.
As $\dif(\crLev)=\crLev$, it follows that $\dif$ yields a homeomorphism $\hdif$ of $\hMman$ permuting cells of $\Theta$.
\jremoved{38.4}{Suppose $\hdif(e)=e$ for some cell $e$ of $\Theta$.
Then $e$ is \emph{$(+)$-invaraint} with respect to $\hdif$ if either $\dim e =0$ or $\hdif$ preserves orientaion of $e$.
Otherwise, $e$ is \emph{$(-)$-invariant}.}
Denote by $\sgmh$ the corresponding permutation of cells of $\Theta$ and let $p$ be its period.

\newcommand\trace{\mathrm{tr}}
\begin{sublemma}\label{lm:number_of_invariant_components}{\rm(cf.~\cite[Proposition~5.4]{Maksymenko:MFAT:2010}).}
Suppose $p>1$.
Then the sets of fixed cells of $\sgmh,\sgmh^2,\ldots,\sgmh^{p-1}$ are the same and consist of two cells $e_1$ and $e_2$ each having even dimension.
In particular, the action of the cyclic group $\langle \sgmh \rangle \cong \bZ_p$ on $\Theta$ is \myemph{semifree}.
\end{sublemma}
\begin{proof}
For $i=0,1,2$, let $c_i$ be the number of $i$-cells in $\Theta$, $C_i =\bZ^{c_i}$ the group of integral $i$-chains, $\sgmh_i:C_i\to C_i$ the chain homomorphism induced by $\sgmh$, and $\phi_i:H_i(\Mman)\to H_i(\Mman)$ the homomorphism induced on integral $i$-th homologies.
Then the \myemph{Lefschetz number} $L(\hdif)$ of $\hdif$ equals
\begin{align*}
L(\hdif) &:= \trace(\phi_0)-\trace(\phi_1)+\trace(\phi_2) =  \trace(\sgmh_0)-\trace(\sgmh_1)+\trace(\sgmh_2) = 
L(\id_{\hMman})=\chi(\hMman) = 2,
\end{align*}
where the first equality is just the definition of $L(\hdif)$, the second one is a principal property of $L(\hdif)$ allowing to compute it as well via induced chain homomorphisms, the third holds because $\hdif$ is an orientation preserving homeomorphism of $S^2$, and so it is isotopic to $\id_{S^2}$, and the forth follows from invariantness of $L(\hdif)$ with respect to homotopies.

As $\sgmh$ is a permutation of cells, it follows that the matrices of chain homomorphisms $\sgmh_i$ consist of numbers $-1$, $0$ and $1$ only.

Suppose $\hdif(e)=e$ for some cell $e\in\Theta$.
We claim that then \myemph{$e$ is \invhp{\hdif} and $\dim e = 0$ or $2$.}
Indeed, for $\dim e=2$, $(\dif,+)$-invariantness of $e$ follows from the assumption that $\hdif$ preserves orientation of $\Mman$; for $\dim e=1$ from the assumptions that $\hdif$ preserves orientation of $\Mman$ and $\func$; and for $\dim e = 0$ this holds by definition.
Moreover, if $\dim e=1$, then the implications~\ref{enum:pres:edge} $\Rightarrow$ \ref{enum:pres:bnd_comp} \& \ref{enum:pres:vertice_with_edges} of Lemma~\ref{lm:shift_func_near_crlevel} show that $\hdif$ preserves all cells of $\Theta$ with their orientations, and so $\sgmh$ is the identity permutation, which contradict to the assumption.

The above statement means that diagonal elements of matrices of $\sgmh_0$ and $\sgmh_2$ are either $0$ and $1$, while all diagonal elements of the matrix of $\sgmh_1$ are zeros.
Hence the number of fixed cells under $\sgmh$ equals 
$\trace(\sgmh_0) + \trace(\sgmh_2) = \trace(\sgmh_0) - \trace(\sgmh_1) + \trace(\sgmh_2) = L(\hdif) = 2$,
and these cells are of even degrees.

Let $e_1, e_2$ be the fixed cells of $\sgmh$.
Then $\sgmh^j$, ($j=1,\ldots,p-1$), also fixes these cells.
On the other hand, by the same arguments as above, $\sgmh^j$ must have exactly $2$ fixed cells, whence those cells are $e_1, e_2$ as well.
\end{proof}

Now let $\gdif\in\Stabilizer{\func,\Vman}$ be any element with $\eta(\gdif)=1$.
Then $\gdif$ generates the cyclic group $\frac{\Stabilizer{\func,\Vman}}{\Stabilizer{\bZman}}$.
Hence applying Lemma~\ref{lm:number_of_invariant_components} to $\gdif$ we obtain that all elements of that group have the same fixed cells $e_1$ and $e_2$ of even dimensions.
Notice that $\XFixA$ is \invhp{\gdif} and thus it corresponds, say to $e_1$.
Therefore $e_2 \in \Theta$ is either a $0$-cell i.e. a critical point of $\func$ in $\crLev$ or a $2$-cell corresponding to an element of $\bZman$.

\ref{enum:cor:MV_DdD_AdA:gm_tau}
Fix any $\qdif\in\Stabilizer{\func,\Vman}$ with $\eta(\qdif)=1$.
Then $\eta(\qdif^{m}) = m$, i.e.\! $\qdif^{m}\in\eta^{-1}(m\bZ)$, whence by Lemma~\ref{lm:eta}\ref{enum:lm:eta:diag:inv_mZ} there exists a $\Cinfty$ function $\alpha:\RB\to\bR$ such that $\qdif^{m}|_{\RB} = \flow_{\alpha}$.

Applying Lemma~\ref{lm:extension_of_shift_functions} to the one-element set $\mathcal{A} = \{\qdif^{m}\}$%
\xadded{-0.8}{0}{41.1}{,} we get that $\qdif^m$ is isotopic in $\Stabilizer{\func}$ to a diffeomorphism $\dif$ fixed on some neighborhood of $\RB$, whence $\dif^{-1}\circ\qdif^{m} \in \StabilizerId{\func}$.
Since $\StabilizerId{\func}$ is a normal subgroup of $\Stabilizer{\func}$, it follows that $\qdif^{i}\circ\dif^{-1}\circ\qdif^{m-i}\in\StabilizerId{\func}$ for each $i=0,\ldots,m-1$.
Hence by Theorem~\ref{th:charact_Stabf} there exists a unique $C^{\infty}$ function $\beta_i:\Mman\to\bR$ such that
$\qdif^{i}\circ\dif^{-1}\circ\qdif^{m-i} = \flow_{\beta_i}$.
Moreover, as \jadded{41.2}{$\dif$ is fixed on $\RB$, we see that} $\qdif^{i}\circ\dif^{-1}\circ\qdif^{m-i} = \qdif^{m}$ on $\RB$, \changed{it follows from}{whence by} Lemma~\ref{lm:shift}\ref{enum:shift:shift_func_uniqueness}\added{,} \removed{that} $\alpha = \beta_i$ on $\RB$ for all $i$.
Hence we get a $C^{\infty}$ function $\gamma:\Mman\to\bR$ given by $\gamma = \alpha$ on $\RB$ and $\gamma=\beta_i$ on $\Yi{i}{j}$.
Now define $\gdif\in\Stabilizer{\func,\Vman}$ by
\[
\gdif(x) =
\begin{cases}
\qdif(x), & x\in \RB \cup
\bigcup\limits_{i=0}^{m-2} \bigcup\limits_{j=1}^{c} \Yi{i}{j}, \\
\dif^{-1}\circ\qdif(x), & x\in \bigcup\limits_{j=1}^{c} \Yi{m-1}{j},
\end{cases}
\]
see Figure~\ref{fig:diff_g}.

\begin{figure}[ht]
	\centering\includegraphics[width=5.3cm]{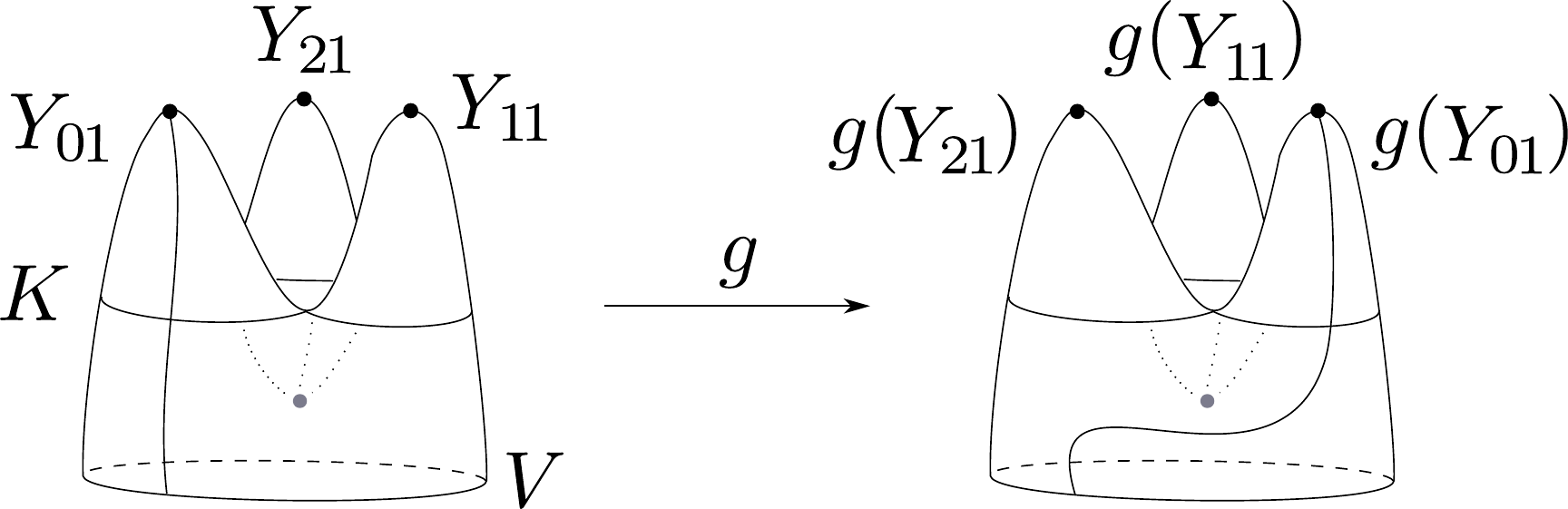}
	\caption{Diffeomorphism $\gdif$}\label{fig:diff_g}
\end{figure}

Then the restriction $\gdif^m|_{\Yi{i}{j}}$ is the composition of the following maps:
\[
\gdif^m|_{\Yi{i}{j}}:
\underbrace{\Yi{i}{j} \xrightarrow{\,\qdif\,} \Yi{i+1}{j} \xrightarrow{\,\qdif\,} \cdots \xrightarrow{\,\qdif\,}
\Yi{m-1}{j} }_{\qdif^{\added{m-i-1}}} \xrightarrow{\,\dif^{-1}\circ\qdif\,}
\underbrace{ \Yi{0}{j} \xrightarrow{\,\qdif\,} \Yi{1}{j} \xrightarrow{\,\qdif\,} \cdots \xrightarrow{\,\qdif\,} \Yi{i}{j} }_{\qdif^{i}},
\]
whence\xadded{-2.5}{0}{41.3}{}
\[
\gdif^m(x) =
\begin{cases}
\qdif^m(x) = \flow_{\alpha}(x) = \flow_{\gamma}(x), & x\in  \RB, \\
\qdif^{i}\circ \dif^{-1} \circ \qdif \circ\qdif^{\added{m-i-1}}(x) = \flow_{\beta_{i}}(x) = \flow_{\gamma}(x), & x\in \Yi{i}{j}.
\end{cases}
\]
\xadded{-3}{0}{41.4}{}
Thus $\gdif^{m} = \flow_{\gamma}\in\StabilizerId{\func}$.
Moreover, as $\gdif=\qdif$ on $\RB$, we see that $\eta(\gdif)=\eta(\qdif)=1$.

\medskip

\xadded{1}{0}{41.6}{}
It remains to show that $[\gdif^{m}] = [\tau]$ in $\pi_0\Stabilizer{\func,\Vman}$.
We have that $\gdif^m = \flow_{\gamma}$ and $\tau = \flow_{\mu}$.
Then by Lemma~\ref{lm:eta}\ref{enum:lm:eta:diag:inv_mZ}, $\gamma(\Vman)=\eta(\gdif^m)/m = \eta(\tau)/m = \mu(\Vman)=1$, whence by Lemma~\ref{lm:shift}\ref{enum:shift:isotopy_between_shifts}, the family of maps $\flow_{(1-t)\gamma + t \mu}$, $t\in\UInt$, is an isotopy between $\gdif^m$ and $\tau$ in $\Stabilizer{\func,\Vman}$.
\end{proof}

{\bf Proof of Theorem~\ref{th:stab:disk_ann:gen_case}.}
\ref{enum:SS:A}
Suppose $\bZmanY=\varnothing$, so $\bZman=\bZmanX=\{\XFixA,\Xman_1,\ldots,\Xman_{a}\}$ and $\Stabilizer{\func,\Vman}=\Stabilizer{\bZman}$.
Then by Lemma~\ref{lm:MV_DdD_AdA}\ref{enum:cor:MV_DdD_AdA:SZ}, $m=1$, whence we have an isomorphism $\seqStab{\func,\Vman} \stackrel{\eqref{equ:split_bseq}}{\cong} \prod\limits_{i=1}^{a}\seqStab{\func|_{\Zman_i},\hZman_i} \times \seqZ{1}$ from Lemma~\ref{lm:eta}\ref{enum:lm:eta:m=1}.
Since $\Mman$ and each $\Zman_i$ is either a $2$-disk or an annulus, it follows from~\eqref{equ:StabIsotIdD2} that one can replace in the latter isomorphism each $\funcSeq$ with $\funcSeq'$, $\Vman$ with $\partial\Mman$, and each $\hZman_i$ with $\partial\Zman_i$.
This gives the required isomorphism 
$\seqStabIsotId{\func,\partial\Mman} \cong \prod\limits_{i=1}^{a}\seqStab{\func|_{\Zman_i},\partial\Zman_i} \times \seqZ{1}$.

\ref{enum:SS:B}
Suppose $\bZmanY\not=\varnothing$, so $\Stabilizer{\bZman} \subsetneq \Stabilizer{\func,\Vman}$.
Then by Lemma~\ref{lm:MV_DdD_AdA}\ref{enum:cor:MV_DdD_AdA:eff_action} the action of $\frac{\Stabilizer{\func,\Vman}}{\Stabilizer{\bZman}} \cong \bZ_m$ on $\bZman$ is semifree.
This implies that each non-fixed orbit consists of $m$ elements, whence $b$ is divided by $m$ and we have $c:=b/m$ orbits.
Let $\Yman_1,\ldots,\Yman_c$ be any elements from $\bZmanY$ belonging to mutually distinct orbits and $\gdif\in\Stabilizer{\func,\Vman}$ be any diffeomorphism such that $\eta(\gdif)=1$.
Then $\gdif$ generates the cyclic group $\frac{\Stabilizer{\func,\Vman}}{\Stabilizer{\bZman}}$.
Denote $\Yi{i}{j} = \gdif^j(\Yman_i)$, $i\in\bZ$, $j=1,\ldots,c$.
Then $\gdif(\Yi{i}{j})=\Yi{i+1}{j}$ and $\Yi{i}{j} = \Yi{i+m}{j}$ for all $i,j$.
This means that
\begin{equation}\label{equ:non_fixed_vertexes}
\begin{matrix}
\Yi{0}{1} & \Yi{1}{1} & \cdots & \Yi{m-1}{1} \\
\Yi{0}{2} & \Yi{1}{2} & \cdots & \Yi{m-1}{2} \\
\cdots    & \cdots   & \cdots & \cdots \\
\Yi{0}{c} & \Yi{1}{c} & \cdots & \Yi{m-1}{c}
\end{matrix}
\end{equation}
is the set $\bZmanY$ of all non-fixed elements of $\bZman$ and $\gdif$ cyclically shifts columns in~\eqref{equ:non_fixed_vertexes}.

\ref{enum:SS:B:B0}
Suppose $\bZmanX = \{\XFixA\}$, so $\bZmanY$ contains all connected components of $\Mman\setminus(\regNK\cup\XFixA)$.
Then by~\eqref{equ:3x3_pi0:seq} we have an exact $(3\times3)$-diagram
$\prod\limits_{i=0}^{m-1} \prod\limits_{j=1}^{c} \seqStab{\func|_{\Yi{i}{j}},\partial\Yi{i}{j}} \monoArrow  \seqStab{\func,\Vman} \epiArrow  \seqZ{m}$ which, in particular, contains an inclusion $\prod\limits_{i=0}^{m-1}
\prod\limits_{j=1}^{c}\pi_0\Stabilizer{\func|_{\Yi{i}{j}},\partial\Yi{i}{j}} \monoArrow\pi_0\Stabilizer{\func,\Vman}$.
Denote by $\zB{i}$, $i=0,\ldots,m-1$, the image of the group $\prod\limits_{j=1}^{c}\pi_0\Stabilizer{\func|_{\Yi{i}{j}},\partial\Yi{i}{j}}$ in 
$\pi_0\Stabilizer{\func,\Vman}$ corresponding to the subsurfaces from the $i$-th column of~\eqref{equ:non_fixed_vertexes}.

By Lemma~\ref{lm:MV_DdD_AdA}\ref{enum:cor:MV_DdD_AdA:gm_tau} one can choose the above  $\gdif\in\Stabilizer{\func,\Vman}$ so that $\gdif^{m}\in\StabilizerId{\func}$, and $[\gdif^m]=[\tau]\in\pi_0\Stabilizer{\func,\Vman}$.
\jadded{41.5-41.6}{}
Since $\gdif$ cyclically shifts columns in~\eqref{equ:non_fixed_vertexes}, it follows from Lemma~\ref{lm:reduction_Mconn_V_dM} that $[\gdif] L_i [\gdif^{-1}] = L_{i+1}$ for each $i$, where $[\gdif]$ is the isotopy class of $\gdif$ in $\pi_0\Stabilizer{\func,\Vman}$ and the indices are taken modulo $m$.
Moreover, as $[\gdif^m]=[\tau]\in\pi_0\Stabilizer{\func,\Vman}$, and $\tau$ is supported in the set $\XFixA$ which is disjoint from all $\Yi{i}{j}$, it follows that $[\gdif^m]$ commutes with each group $\zB{i}$.
Hence by Lemma~\ref{lm:charact_seq_wr} we have an isomorphism 
$\seqStab{\func,\Vman} \cong 
\seqWrm{\bigl(\prod\limits_{j=1}^{c}\seqStab{\func|_{\Yi{0}{j}},\partial\Yi{0}{j}}\bigr)}{m}$, which is the same as $\seqStabIsotId{\func,\partial\Mman} \cong 
\seqWrm{\bigl(\prod\limits_{j=1}^{c}\seqStabIsotId{\func|_{\Yman_j},\partial\Yman_j}\bigr)}{m} \stackrel{\eqref{equ:DOS__jSdS_case_B}}{\equiv} \aSeq$ since each $\Yman_j = \Yi{0}{j}$ is a $2$-disk.

Proof of statement~\ref{enum:SS:B:B0_X1} is explained in its formulation.
Theorem~\ref{th:stab:disk_ann:gen_case} is completed.
\qed

\medskip

{\bf Acknowledgement.}
The author is grateful to Bogdan Feshchenko for useful discussions.
\added{The author is also indebted to the anonymous Referee for very careful reading of the initial version of the manuscript and a lot of suggestions that allowed to clarify all the exposition and correct misprints and several wrong formulations.}


\def\cprime{$'$}
\providecommand{\bysame}{\leavevmode\hbox to3em{\hrulefill}\thinspace}
\providecommand{\MR}{\relax\ifhmode\unskip\space\fi MR }
\providecommand{\MRhref}[2]{%
  \href{http://www.ams.org/mathscinet-getitem?mr=#1}{#2}
}
\providecommand{\href}[2]{#2}

\end{document}